\documentclass[11pt,a4paper,reqno]{amsart}
\usepackage{graphicx} 

\usepackage{amsmath}
\usepackage[foot]{amsaddr}

\usepackage{enumitem}
\usepackage{xcolor}
\usepackage{amssymb}
\usepackage{amsthm}
\usepackage{array}
\newcolumntype{M}[1]{>{\centering\arraybackslash}m{#1}}

\usepackage[backend=bibtex]{biblatex}
\usepackage{float}

\usepackage{cuted}
\usepackage{color}
\usepackage{graphicx}

\usepackage{tikz}
\usetikzlibrary{shapes.geometric}
\usetikzlibrary{arrows.meta,arrows}

\usepackage{tabularx}

\usepackage{comment}

\usepackage[hidelinks, colorlinks=false]{hyperref}

\usepackage{cleveref}

\usepackage[normalem]{ulem}

\newtheorem{thm}{Theorem}
\newtheorem{remark}[thm]{Remark}
\newtheorem{proposition}[thm]{Proposition}
\newtheorem{corollary}[thm]{Corollary}
\newtheorem{lemma}[thm]{Lemma}
\newtheorem{assumption}[thm]{Assumption}
\newtheorem{definition}[thm]{Definition}
\newtheorem{example}[thm]{Example}

\newcommand{\R}{\mathbb{R}}
\newcommand{\N}{\mathbb{N}}
\newcommand{\Z}{\mathbb{Z}}

\newcommand{\scp}[2]{ \left \langle #1, #2 \right\rangle }
\newcommand{\n}[1]{\left\Vert #1\right\Vert}

\newcommand{\eqb}{\begin{flalign}}
\newcommand{\eqe}{\end{flalign}}

\newcommand{\Lra}{\Leftrightarrow}

\newcommand{\platz}{\;\;\;\;\;}

\newcommand{\M}{\mathcal{M}}
\newcommand{\Benedikt}[1]{\textcolor{black}{#1}}

\allowdisplaybreaks

\definecolor{mygrey}{RGB}{220, 220, 220}
\definecolor{concrete}{RGB}{136, 150, 136}
\definecolor{steel}{RGB}{200, 200, 200}
\definecolor{copper}{RGB}{220, 127, 100}
\definecolor{redd}{RGB}{200, 50, 50}
\definecolor{yellow}{RGB}{255, 254, 235}
\definecolor{green}{RGB}{213, 243, 231}
\definecolor{green2}{RGB}{169, 209, 142}
\definecolor{lightgrey}{RGB}{240, 240, 240}
\definecolor{lightred}{RGB}{255, 244, 239}
\definecolor{lightlilac}{RGB}{239, 239, 255}
\definecolor{wood}{RGB}{133, 94, 66}
\definecolor{BackgroundGrey}{RGB}{210, 210, 210}

\begin{document}

\title[Spatial exponential decay of perturbations in optimal control]{Spatial exponential decay of perturbations in optimal control of general evolution equations}
\thanks{B.O.\ extends his gratitude to the Thüringer Graduiertenförderung for funding his scholarship. K.W.\ gratefully acknowledges funding by the German Research Foundation (DFG) – Project-ID 507037103.}
\author{Simone Göttlich$^1$}\address{$^1$Chair of Scientific Computing,  School of Business Informatics and Mathematics, University of Mannheim, Germany.\\ Mail: \textsc{goettlich@uni-mannheim.de}}
\author{Benedikt Oppeneiger$^2$}
\address{$^2$Optimization-based Control Group, Institute of Mathematics, Technische Universität Ilmenau, Germany.\\ Mail: \textsc{\{benedikt-florian.oppeneiger, karl.worthmann\}@tu-ilmenau.de}}
\author{Manuel Schaller$^3$}\address{$^3$Faculty of Mathematics, Chemnitz University of Technology, Germany.\\ Mail: \textsc{manuel.schaller@math.tu-chemnitz.de}}
\author{Karl Worthmann$^2$}

\date{\today}
\begin{abstract} 
    We analyze the robustness of optimally controlled evolution equations with respect to spatially localized perturbations. We prove that if the involved operators are domain-uniformly stabilizable and detectable, then these localized perturbations only have a local effect on the optimal solution. We characterize this domain-uniform stabilizability and detectability for the transport equation with constant transport velocity, showing that even for unitary semigroups, optimality implies exponential damping. \Benedikt{We extend this result to the case of a space-dependent transport velocity. Finally we leverage the results for the transport equation to characterize domain-uniform stabilizability of the wave equation.} Numerical examples in one space dimension complement the theoretical results.
\end{abstract}
\maketitle
\section{Introduction}
The robustness of models with regard to perturbations is a fundamental aspect of many fields of applied mathematics ranging from functional analysis or operator theory (e.g.\ bounded perturbation theorems)
to control theory (synthesis of robust controllers) and scientific computing (analysis of error propagation in numerical algorithms).

When considering optimal control problems, the robustness of solutions is a vivid topic of research. In case of dynamic optimal control, e.g.\ problems governed by evolution equations, the long-term behavior of optimal solutions for increasing time horizon may be described by the turnpike property, cf.~the recent overview articles~\cite{Faulwasser2022Turnpike,GeshZuaz22}. This qualitative property states that optimal solutions reside close to an optimal steady state for the majority of the time. Therefore optimal solutions are - in a certain sense - robust with regard to changes in the initial and terminal conditions. To quantify the closeness to the optimal steady state, exponential versions of the turnpike property have been considered - even for problems governed by PDEs - for example via dissipativity-based ana\-lysis~\cite{Damm2014} or by considering the first-order optimality conditions~\cite{TrelZuaz15,Breiten2020}. Further, the underlying stability of the optimal control problem leading to the turnpike property also implies a particular robustness with respect to numerical errors, cf.~\cite{Schaller2019,Schaller2020} and may be used for efficient predictive control of PDEs via goal-oriented mesh refinement~\cite{Schaller2022}.

Only recently, an exponential decay of sensitivities with regard to perturbations in space was considered and proven, e.g., in \cite{Shin2022} for finite-dimensional nonlinear optimization problems on graphs under second order sufficient conditions that are uniform in the network size. 
Further, such an exponential sensitivity was leveraged in \cite{SperSalu23} for efficient neural-network-based approximations of separable optimal value functions.
For linear-quadratic optimal control of stationary (such as Poisson and Helmholtz equation) and parabolic PDEs, an exponential decay of perturbations in space was shown in~\cite{Göttlich2024} under stabilizability and detectability assumptions, linking the decay in optimal control to system-theoretic properties. 

In this work, we extend the result of the previous work \cite{Göttlich2024} to general evolution equations governed by semigroups of operators. 
This class includes a large number of highly relevant systems (e.g.\ wave, transport, beam equations), which do not exhibit high regularity or a parabolic structure. 
Furthermore, our analysis features a weaker stabilizability/detectability assumption in the sense that we allow for the stabilized solution to exhibit an overshoot behaviour. 
On the contrary, in~\cite{Göttlich2024}, such an overshoot behaviour is not admissible. The essence of the stabilizability/detectability assumption is, that the constants in the exponential estimate of the stabilized semigroups are uniform in the size of the spatial domain.

As a second contribution, we thoroughly investigate this assumption for several hyperbolic equations on a one-dimensional domain, namely the transport equation with constant and space-dependent velocity, the continuity equation and the wave equation. For all of these examples, we propose an easily-verifiable novel necessary and sufficient condition on the control and observation domain to ensure domain-uniform stabilizability/detectability.
    
This work is structured as follows. In~\Cref{Sec: ProblemStatement} the optimal control problem under consideration is specified and the optimality system is derived. In~\Cref{Sec: MainResults} we provide our first main results on the exponential decay of perturbations in optimal control of hyperbolic equations under a domain-uniform stabilizability/detectability assumption. In~\Cref{Sec: TransportConstant} we characterize this stabilizability/detectability assumption for the transport equation with constant coefficients and in~\Cref{Sec: TransportNonConstant} also for space-dependent transport velocities. In~\Cref{Sec: WaveEq} we extend these results to the wave equation with Dirichlet boundary conditions. Finally, we illustrate our findings by numerical experiments in~\Cref{Sec: Numerics}.

\subsection{Nomenclature}
Let~$\N$ denote the natural numbers, $\N_0 = \N \cup\{0\}$, $\R_{\ge 0}=[0,\infty)$ and $\R_{> 0}=(0,\infty)$. The Lebesgue measure of a measurable set $\Omega$ denoted by $|\Omega|$ and $L^2(\Omega)$ denotes the space of square integrable Lebesgue measurable functions. For a Hilbert space $X$ and $T>0$, we denote by $L^p(0,T;X), p\in [1,\infty)$, the space of $X$-valued $p$-Bochner-integrable 
functions endowed with the norm $\n{f}_{L^p(0,T;X)}^p:= \int _0^T \n{f(t)}_X^p\mathrm{d}t$. Further, $L^\infty(0,T;X)$ stands for the space of $X$-valued Bochner-integrable, essentially bounded functions with norm $\|\cdot\|_{L^\infty(0,T;X)} = \text{esssup}_{t \in [0,T]} \|f(t)\|$.
The space 
of $X$-valued continuous functions over $[0,T]$ is denoted by $C(0,T;X)$  with norm $\|\cdot\|_{C(0,T;X)} := \max_{t\in [0,T]} \|f(t)\|$.

\section{Problem statement}
\label{Sec: ProblemStatement}

In this work, we analyse the sensitivity of optimal control problems in view of increasing domain sizes. To this end, we consider a family of domains, and we parameterize the optimal control problem by means of the spatial domain and the optimization horizon.

\begin{definition}
    Let $d\in \N$ and $\mathcal{O}\subset \{\Omega \subset \R^d : \Omega \textrm{ open, bounded with Lipschitz boundary}\}$. We consider a family $\left( \textrm{OCP}_{\Omega} ^T \right)_{\Omega \in \mathcal{O}, T>0}$  of optimal control problems (OCPs) given by
\begin{flalign}
\tag{$\mathrm{OCP}_{\Omega} ^T$}
\label{HyperbolicOCP}
    \begin{split}
    \underset{(x,u)}{\min }\,\,\, \frac{1}{2} \int _0^T \n{C_\Omega \left(x(t)-x^\mathrm{ref}_\Omega\right)}_{Y_\Omega}^2 \, + \, & \n{R_\Omega\left(u(t)-u^\mathrm{ref} _\Omega\right)}_{U_\Omega}^2 \, \mathrm{d}t\\
    \textrm{s.t.}: \platz \dot{x}-A_\Omega x - B_\Omega u &= f_\Omega, \,\, x(0)=x_\Omega^0,
    \end{split}
\end{flalign}
where for each $\Omega \in \mathcal{O}$
\begin{itemize}
    \item $X_\Omega=L^2(\Omega)$ and $A_\Omega: X_\Omega \supset D(A_\Omega)\rightarrow X_\Omega$ generates a strongly-continuous semigroup $(\mathcal{T}_\Omega(t))_{t\geq 0}$ on~$X_\Omega$.
    \item $B_\Omega \in L(U_\Omega,X_\Omega)$ with Hilbert space $U_\Omega$.
    \item $C_\Omega \in L(X_\Omega,Y_ \Omega)$ with Hilbert space $Y_\Omega$ and $R_\Omega \in L(U _\Omega,U _\Omega)$ with $\alpha > 0$ such that \Benedikt{$\forall\, \Omega \in \mathcal{O}: \n{R_\Omega u}_{U_\Omega}^2\geq \alpha \n{u}^2_{U_\Omega}$}, $u \in U_\Omega$.
    \item $f_\Omega \in L^1(0,T;X _\Omega)$, $x^0_\Omega \in X_\Omega$, $x^\mathrm{ref}_\Omega \in X_\Omega$ and $u^\mathrm{ref} _\Omega\in U_\Omega$.
\end{itemize}
\end{definition}
If $A_\Omega$ models a differential operator, such as in partial differential equations, boundary conditions are usually included in its domain $D(A_\Omega)$. For applications of this abstract framework, we refer to the example of the transport equation with periodic boundary conditions that is extensively studied in Sections~\ref{Sec: TransportConstant} and \ref{Sec: TransportNonConstant}.
Here, the dynamics of \eqref{HyperbolicOCP} is meant in a mild sense, that is, $(x,u)$ is admissible for \eqref{HyperbolicOCP} if for all $t\in [0,T]$
\begin{align*}
    x(t) = \mathcal{T}_\Omega(t) x_\Omega ^0 + \int_0^t \mathcal{T}_\Omega(t-s)B_\Omega u(s)\,\mathrm{d}s.
\end{align*}
Correspondingly, the constraint and the state may be eliminated by above variation of constants formula of semigroup theory. Thus, via standard methods cf.~\cite{Tröltzsch2010}, one may conclude that \eqref{HyperbolicOCP} has a unique solution $x_\Omega^T$ and $u_\Omega^T$, see, e.g., \cite[Theorem 2.24]{Schaller2021}. 
If clear from context, we will mostly omit the indices~$T$ and~$\Omega$ for readability. 

In the following, we will always identify the Hilbert spaces $X_\Omega = L^2(\Omega)$, $U_\Omega$, $Y_\Omega$  with their respective dual spaces via the Riesz isomorphism. 
\\

\noindent \textbf{Optimality conditions.} Our analysis will be based on the first-order optimality conditions, i.e.~Pontryagin's Maximum Principle~\cite{YongLi1995}. If $(x,u)\in C(0,T;X_\Omega)\times L^2(0,T;U_\Omega)$ is optimal for~\eqref{HyperbolicOCP}, there exists a Lagrange multiplier 
$\lambda \in C(0,T;X_\Omega)$ such that the optimality system
\begin{equation*}
    \begin{pmatrix}
        C^*C & 0 & -\frac{\mathrm{d}}{\mathrm{d}t} - A ^*\\
        0 & 0 & E^T\\
        0 & R^*R & -B^*\\
        \frac{\mathrm{d}}{\mathrm{d}t} - A & -B & 0\\
        E^0 & 0 & 0
    \end{pmatrix}
    \begin{pmatrix}
        x\\ u\\ \lambda 
    \end{pmatrix}
    =
    \begin{pmatrix}
        C^*Cx^\mathrm{ref}\\
        0\\
        R^*R u^\mathrm{ref}\\
        f\\
        x^0
    \end{pmatrix}
\end{equation*}
is fulfilled in a mild sense.
The operators $E^0$ and $E^T$ capture the initial and terminal conditions, that is,
\begin{equation*}
    E^0: C(0,T;X_\Omega) \rightarrow X_\Omega, \platz E^0 x := x(0)
    \qquad
    \mathrm{and}
    \qquad
    E^T: C(0,T;X_\Omega) \rightarrow X_\Omega, \platz E^T \lambda := \lambda(T).
\end{equation*}
Setting $Q:=R^*R$, we may eliminate the control via $u=Q^{-1}B^*\lambda + u^\mathrm{ref}$. Note that $Q$ is continuously invertible due to the ellipticity of $R$. This leads to the condensed optimality system
\begin{equation}
\label{CondensedOptimalitySystem}
    \begin{pmatrix}
        C^*C & -\frac{\mathrm{d}}{\mathrm{d}t} - A ^*\\
        0 & E^T\\
        \frac{\mathrm{d}}{\mathrm{d}t} - A & -BQ^{-1}B^*\\
        E^0 & 0
    \end{pmatrix}
    \begin{pmatrix}
        x\\ \lambda 
    \end{pmatrix}
    =
    \underset{:=h}{\underbrace{\begin{pmatrix}
        C^*Cx^\mathrm{ref}\\
        0\\
        Bu^\mathrm{ref} + f\\
        x^0
    \end{pmatrix}}},
\end{equation}
where the first two rows correspond to the adjoint equation and the last two rows are correspond to the state equation. 
In the following, we will abbreviate the optimality system by means of the unbounded operator
\begin{align*}
    \M: D(\M)\subset C(0,T;X_\Omega)^2 \rightarrow \left( L^1 (0,T;X_\Omega) \times X_\Omega \right)^2, \qquad \begin{pmatrix}
        x\\ \lambda 
    \end{pmatrix}
    \mapsto \begin{pmatrix}
        C^*C & -\frac{\mathrm{d}}{\mathrm{d}t} - A ^*\\
        0 & E^T\\
        \frac{\mathrm{d}}{\mathrm{d}t} - A & -BQ^{-1}B^*\\
        E^0 & 0
    \end{pmatrix}
    \begin{pmatrix}
        x\\ \lambda 
    \end{pmatrix}
\end{align*}
with $D(\M) = C^1(0,T;X_\Omega)\cap C(0,T;D(A)) \times C^1(0,T;X_\Omega)\cap C(0,T;D(A^*))$. Note that, when endowed with this domain, $M$ is not closed. However, as we only use it to concisely denote the optimality system, we do not rely on analytic properties of $M$.

\section{Exponential decay of optimal solutions in space}
\label{Sec: MainResults}

In this section, the main result of this paper is presented. 
We show, that under the assumption of domain-uniform stabilizability and detectability, the influence of spatially-localized perturbations decays exponentially in optimally-controlled systems. 
Here, domain-uniform stabilizability refers to constants in the closed-loop semigroup that are uniform w.r.t.\ the family of domains. 
This result provides robustness w.r.t.\ disturbances that could, e.g., be caused by discretizing the equation system in~\eqref{CondensedOptimalitySystem} to compute an approximate solution $(\hat{x},\hat{\lambda})$.

We first introduce the notion of exponential localization, which serves to quantify the locality of perturbations and their effect. 
We consider the perturbed version of~\eqref{CondensedOptimalitySystem}
\begin{equation}
\label{DisturbedOptimalitySystem}
    \begin{pmatrix}
        C^*C & -\frac{\mathrm{d}}{\mathrm{d}t} - A ^*\\
        0 & E^T\\
        \frac{\mathrm{d}}{\mathrm{d}t} - A & -BQ^{-1}B^*\\
        E^0 & 0
    \end{pmatrix}
    \begin{pmatrix}
        x^d\\ \lambda ^d 
    \end{pmatrix}
    =
    \begin{pmatrix}
        C^*Cx^\mathrm{ref}\\
        0\\
        Bu^\mathrm{ref} + f\\
        x^0
    \end{pmatrix}
    +
    \begin{pmatrix}
        \varepsilon _1\\
        \varepsilon _2\\
        \varepsilon _3\\
        \varepsilon _4
    \end{pmatrix}
\end{equation}
with disturbance $\varepsilon = \begin{pmatrix}
    \varepsilon _1 & \varepsilon_2 & \varepsilon _3 & \varepsilon _4
\end{pmatrix}^\top \in (L^1 (\R _{\geq 0};L^2(\R^d))\times L^2(\R^d))^2$, where in the above equation, $\varepsilon$ is considered to be restricted to $[0,T]$ and $\Omega$. 
By defining the error variables 
\[
    \delta x := x^d - x \qquad\text{ and }\qquad \delta \lambda := \lambda ^d - \lambda
\]
and subtracting~\eqref{CondensedOptimalitySystem} from~\eqref{DisturbedOptimalitySystem} we find the error system
\begin{equation}
\label{ErrorOptimalitySystem}
   \M \delta z = \begin{pmatrix}
        C^*C & -\frac{\mathrm{d}}{\mathrm{d}t} - A ^*\\
        0 & E^T\\
        \frac{\mathrm{d}}{\mathrm{d}t} - A & -BQ^{-1}B^*\\
        E^0 & 0
    \end{pmatrix}
    \begin{pmatrix}
        \delta x\\ \delta \lambda 
    \end{pmatrix}
    =
    \begin{pmatrix}
        \varepsilon _1\\
        \varepsilon _2\\
        \varepsilon _3\\
        \varepsilon _4
    \end{pmatrix} = \varepsilon.
\end{equation}

\Benedikt{
\begin{remark}
    In optimize-then-discretize fashion, the condensed optimality system~\eqref{CondensedOptimalitySystem} can be solved using finite element methods. This leads to an approximate solution $z_\mathcal{T} \in V_\mathcal{T}$, where in the simplest case $V_\mathcal{T}$ is a space of piecewise affine-linear polynomials on a triangular mesh $\mathcal{T}=\{T_1,\ldots ,T_N\}, T_i \subset \R^2$. However, since $z_\mathcal{T}$ only solves the Galerkin-projected version of~\eqref{CondensedOptimalitySystem}, there is a residual $\varepsilon_\mathcal{T}$ such that
    \begin{equation*}
        \mathcal{M} z_\mathcal{T} = h + \varepsilon _\mathcal{T}.
    \end{equation*}
    This corresponds to the disturbed optimality system~\eqref{DisturbedOptimalitySystem}, i.e. the residual can be interpreted as a perturbation of~\eqref{CondensedOptimalitySystem}. Note that if the resolution of the discretization mesh is fine in certain areas, then the residual in these areas also becomes locally (in space) small (see for example~\cite[Section 10]{Ern2004}). Therefore, the residual can be interpreted as a local perturbation induced by the discretization.
\end{remark}
}

\begin{definition}[Exponential localization]
    \label{ExponentialLocalization}
    \Benedikt{%
    The family of functions $(g_\Omega ^T)_{\Omega \in \mathcal{O}, T>0}$, $g_\Omega ^T \in F_\Omega ^T := (E_\Omega ^T)^m\times (X_\Omega)^n$, $m,n \in \mathbb{N}_0$, with $E_\Omega ^T \in \{L^1(0,T;X_\Omega), L^2(0,T;X_\Omega),C(0,T;X_\Omega)\}$ and Hilbert space of functions~$X_\Omega$ on $\Omega \in \mathcal{O}$, where $\mathcal{O}$ is a set of measurable spatial domains in $\R^d$, 
    is called $\mathcal{O}$-domain-uniformly exponentially localized around $P\in \R ^d$ if there exist constants $\mu _g>0$ and $C_g >0$ such that 
    \begin{equation*}
        \left\| e^{\mu _g \n{P-\cdot}}g_\Omega ^T \right\|_{F_\Omega ^T}< C_g < \infty \qquad\forall\,T > 0, \Omega \in \mathcal{O}.
    \end{equation*}}
\end{definition}

\noindent The above estimate provides a bound in an exponentially weighted norm as considered in~\cite{Breiten2020} and~\cite{Lykina2017}. In this way, it quantifies the decay of a family of functions around a certain point $P$ in space.

The main result of the paper is to show the following sensitivity result of~\eqref{ErrorOptimalitySystem}: \Benedikt{If the family of perturbations $(\varepsilon_\Omega ^T)_{\Omega \in \mathcal{O}, T>0}$ is exponentially localized around $P \in \R^d$ either in $((L^1(0,T;X_\Omega)^2 \times (X_\Omega)^2)_{\Omega \in \mathcal{O},T>0}$ or in $((L^2(0,T;X_\Omega)^2 \times (X_\Omega)^2)_{\Omega \in \mathcal{O},T>0}$ then the family of error variables $(\delta z_\Omega ^T)$ is exponentially localized around $P$ in both 
$(L^2(0,T;X_\Omega)^2)_{\Omega \in \mathcal{O},T>0}$ and  
$(C(0,T;X_\Omega)^2)_{\Omega \in \mathcal{O},T>0}$. 
In particular this means, that the influence of the perturbations on areas in space, which are \textit{far away} from the fixed point $P$ is negligible.}

 Overall we want to find four estimates which correspond to the four elements of the cartesian product given by $\{(L^1 (0,T;X)\times X)^2,(L^2 (0,T;X)\times X)^2\}\times\{(L^2 (0,T;X)\times X)^2,C(0,T;X)^2\}$. To avoid complicated case distinctions and in order to simplify the presentation of our results the following abbreviations are used for estimations:
    \begin{equation*}
        \n{v}_{2 \land \infty}: = \max \left\{\n{v}_{L^2(0,T;X)}, \n{v}_{C(0,T;X)}\right\} \textrm{  and  }
        \n{v}_{1 \lor 2}: = \min \left\{\n{v}_{L^1(0,T;X)}, \n{v}_{L^2(0,T;X)}\right\}.
    \end{equation*}
Using these shorthands we can prove estimates in various norms while also avoiding lengthy case distinctions.
In addition we define the linear spaces
\begin{align*}
    W^{i}:=\left(\left(L^i(0,T;X),\n{\cdot}_{L^i}\right)\times \left(X,\n{\cdot}_{X}\right)\right)^2, \quad W^{2 \land \infty}:=\left(C(0,T;X),\n{\cdot}_{2\land \infty}\right)^2\\
    W^{1\lor2}:=\left(\left(L^1(0,T;X),\n{\cdot}_{1\lor 2}\right)\times \left(X,\n{\cdot}_{X}\right)\right)^2 \hspace{2cm} 
\end{align*}
for $i\in \{1,2\}$. 
To prove the main result of this section we will replace the error variables $\delta x$ and $\delta \lambda$ in~\eqref{ErrorOptimalitySystem} by their scaled versions $\delta \tilde{x}:= e^{\mu \n{P-\omega}} \delta x$ and $\delta \tilde{\lambda}:= e^{\mu \n{P-\omega}} \delta \lambda$ and rewrite the equation system accordingly. 
However to succeed with this strategy we need some information on how the scaling function $\rho (\omega) := e^{\mu \n{P-\omega}}$ from~\Cref{ExponentialLocalization} influences the action of the operators $A$, $B$ and $C$. 
For the input and output operator we 
assume that they are homogeneous as stated in the following assumption.
\begin{assumption}[Uniformity and homogeneity of input, output and weighting]
\label{Ass: ActHom}
    For all $\Omega \in \mathcal{O}$, the operators $B_\Omega$ and $C_\Omega$ are homogeneous with regard to the exponential decay function in the sense that for all $\mu >0$ for all $P \in \R ^d$ we have
    \begin{align}\label{eq:hom}
         \left \langle B_\Omega(e^{\mu \n{P-\cdot}} u), v\right\rangle _{X_\Omega}=\left \langle B_\Omega u, e^{\mu \n{P-\cdot}} v\right\rangle _{X_\Omega} \quad \textrm{ and } \quad C_\Omega(e^{\mu \n{P-\cdot}} x(\cdot)) = e^{\mu \n{P-\cdot}}  C_\Omega x(\cdot).
    \end{align}
Further, let $B_\Omega$, $C_\Omega$ and $R_\Omega$ be bounded uniformly in $\mathcal{O}$, i.e. there exist constants $C_B>0$, $C_C>0$ and $C_R>0$ such that
    \begin{equation*}
        \forall \, \Omega \in \mathcal{O}: \n{B_\Omega}_{L(U_ \Omega,X_ \Omega)} \leq C_B,\quad \n{C_\Omega}_{L(X_\Omega,Y_ \Omega)}\leq C_C \quad\text{ and }\quad \n{R_\Omega}_{L(U_\Omega,U _\Omega)}\leq C_R
    \end{equation*}
and let the ellipticity constant of $R_\Omega \in L(U _\Omega,U _\Omega)$ be uniform in $\mathcal{O}$, i.e., there is $\alpha > 0$ such that $\n{R_\Omega u}_{U_\Omega}^2\geq \alpha \n{u}^2_{U_\Omega}$ for all $u \in U_\Omega$ and all $\Omega \in \mathcal{O}$.
    
\end{assumption}
\begin{example}[Distributed control]
	    Let $\Omega_c \subset \R^d$ be measurable. Then the input operator
	\begin{equation*}
		B_\Omega: L^2(\Omega_c \cap \Omega) \rightarrow X_ \Omega, \platz     		(B_\Omega u)(\omega):= \left\{ \begin{array}{cc}
	      u(\omega), & \omega \in \Omega_c\\
	      0,    & \mathrm{else}
	     \end{array}
	     \right..
	\end{equation*}
    fulfills the homogeneity condition from~\Cref{Ass: ActHom}. Also we find
    \begin{equation*}
        \forall u \in  L^2(\Omega_c \cap \Omega): \n{B_\Omega u}_{X_\Omega}^2 = \int _{\Omega} \n{(B_\Omega u)(\omega )}^2 \mathrm{d}\omega = \int _{\Omega_c \cap \Omega} \n{u(\omega )}^2 \mathrm{d}\omega = \n{u}_{L^2(\Omega_c \cap \Omega)}^2
    \end{equation*}
    which shows uniform boundedness of $B_\Omega$.
\end{example}

\noindent Differential operators do not satisfy the above homogeneity assumption \eqref{eq:hom} but instead yield a perturbation of the original differential operator. Exemplarily for the unbounded gradient operator $L^2(\Omega)$ with domain $H^1(\Omega)$, we get for $x\in H^1(\Omega)$ and $v\in L^2(\Omega;\R^d)$
\begin{equation*}
    \left \langle \nabla(e^{\mu \n{P-\cdot}} x),v\right\rangle_{X_\Omega} = \left \langle e^{\mu  \n{P-\cdot}} \nabla x,v\right\rangle_{X_\Omega} + \left \langle \mu \, \mathrm{sgn}(P-\cdot) \, e^{\mu \n{P-\cdot}} x,v\right\rangle_{X_\Omega}
\end{equation*}
which is a bounded perturbation of the unbounded gradient operator. 

This observation motivates the following structural assumption capturing a wide range of differential operators.
\begin{assumption}[Generators are differential-operator-like]
\label{Ass: Scaling}
    For each $\Omega \in \mathcal{O}$ there exist operator families $\left( S_{i,\Omega}^\mu\right)_{\mu >0} \subset L(X_\Omega), i\in \{1,2\}$ with the following properties:
    \begin{itemize}

    \item[(i)] If $x\in \mathrm{dom}(A _\Omega)$ solves
        \begin{equation*}
         A_\Omega x=f
    \end{equation*}
    for some $f \in X_\Omega$, then  the scaled quantity $\tilde{x}:= e^{\mu\n{P-\cdot}_1}x$ satisfies $\tilde x \in \mathrm{dom}(A_\Omega)$ and solves
    \begin{flalign*}
        (A_\Omega+S_{1,\Omega}^\mu)\tilde{x} = e^{\mu\n{P-\cdot}_1}f=:\tilde{f}.
    \end{flalign*}
    \item[(ii)] If $\lambda \in \mathrm{dom}(A^*_\Omega)$ solves
    \begin{equation*}
         A_\Omega^*\lambda=g,
    \end{equation*}
    then the scaled quantity $\tilde{\lambda}:= e^{\mu\n{P-\cdot}_1}\lambda$ satisfies $\tilde \lambda \in \mathrm{dom}(A^*)$ and solves
    \begin{flalign*}
        (A_\Omega^*+S_{2,\Omega}^\mu)\tilde{\lambda} = e^{\mu\n{P-\cdot}._1}g =: \tilde{g}.
    \end{flalign*}
        \item[(iii)] The operator norms converge towards the zero operator in the uniform operator topology uniformly in the domain size, that is, setting $\mathcal{O}_{>r_0}:=\{\Omega \in \mathcal{O}: |\Omega|>r_0\}$, 
        \begin{equation*}
            \forall \, r_0 >0 \, \forall \varepsilon >0 \, \exists \, \delta >0 \, \forall \, \Omega \in \mathcal{O}_{>r_0} : \mu < \delta \implies \n{S_{i,\Omega}^\mu}_{L(X, X)}< \varepsilon. 
        \end{equation*}
   \end{itemize}
           
\end{assumption}

The following example illustrates Assumption~\ref{Ass: Scaling}.

\begin{example}[First-order problems: One-dimensional transport equation]\label{ex:firstorder}
    We first consider the example of the transport equation that is extensively studied in Sections~\ref{Sec: TransportConstant} and \ref{Sec: TransportNonConstant}. The generator of this partial differential equation on $[0,T]\times \Omega_L$ with $\Omega_{L}:=[0,L]$ with periodic boundary conditions is given by
\begin{equation*}
     A_{\Omega _L}: D(A_{\Omega _L}):= \{x \in H^1(\Omega _L): x(0)=x(L)\}\subset L^2(\Omega_L)  \rightarrow L^2(\Omega _L), \platz A_{\Omega _L}x = -c \frac{\partial}{\partial \omega}x
 \end{equation*}
 with $c\in L^\infty(\Omega_L)$ essentially bounded from below and above by positive constants.
 Correspondingly, we may compute by means of the chain rule for weak derivatives (cf.\ also \cite[Lemma 3]{Göttlich2024}),
 \begin{align*}
     (A_{\Omega_L} \tilde x)(\omega) = c(\omega)\frac{\partial}{\partial \omega} \left(e^{\mu \|P-\omega\|_1} x(\omega)\right) &= c(\omega) \left( e^{\mu |P-\omega|}\frac{\partial}{\partial \omega} x(\omega)-\mu \operatorname{sgn} (P-\omega) e^{\mu |P-\omega|} x(\omega)\right)\\
     &= e^{\mu |P-\omega|} c(\omega) \frac{\partial}{\partial \omega} x(\omega) - c(\omega)\mu \operatorname{sgn} (P-\omega) e^{\mu |P-\omega|} x(\omega).
 \end{align*}
Hence,
\begin{align*}
     (A_{\Omega_L} + S^\mu_{1,\Omega_L}) \tilde x =  e^{\mu |P-\omega|} A_{\Omega_L} x(\omega)
\end{align*}
 with the bounded and linear multiplication operator 
 \begin{align*}
     S^\mu_{1,\Omega_L} : L^2(\Omega_L) \to L^2(\Omega_L),\qquad 
     x \mapsto c\mu \operatorname{sgn} (P-\cdot) x.
 \end{align*}
As $\|S^\mu_{1,\Omega_L}\|_{L(L^2(\Omega_L),L^2(\Omega_L))} \leq \mu \|c\|_{L^\infty(\Omega_L)}$, the uniform convergence property of Assumption~\ref{Ass: Scaling}(iii) is satisfied if $\|c\|_{L^\infty(\Omega_L)}$ is bounded uniformly in $L$.
\end{example}

\begin{example}[Second-order problems: Multi-dimensional wave equation]
\label{ex:vs}
As a second example, we consider distributed control of a wave equation on $\Omega \subset \mathbb{R}^d$, $d\in \N$, given by
\begin{flalign}
\label{eq: Waveex}
\begin{split}
 \rho(\omega) \tfrac{\partial^2}{\partial t^2}{w}(\omega,t) &= \operatorname{div}(\mathcal{T}(\omega) \nabla {w}(\omega,t)) + \chi_{\Omega_c}(\omega) u(t,\omega),  \quad\omega \in \Omega, t \ge 0\\
 w(\gamma,t) &= 0  \quad \gamma \in \partial \Omega, t\geq 0,
\end{split}
\end{flalign}
where $w:\Omega \times [0,T]\to \mathbb{R}$ models the displacement of a membrane with respect to the rest position, where $\rho,\mathcal{T}\in L^\infty(\Omega)$ with $\rho^{-1},\mathcal{T}^{-1}\in L^\infty(\Omega)$ are mass density and Young's modulus, respectively. A naive reformulation as a first order equation by introducing a velocity variable yields the generator
\begin{align*}
    A_{\Omega,1} = \begin{pmatrix}
        0 & I \\
        \tfrac1\rho\operatorname{div}(\mathcal{T} \nabla \cdot) & 0 
    \end{pmatrix},
\end{align*}
e.g.\ on $L^2(\Omega)\times H^{-1}(\Omega)$ with $D(A_{\Omega,1}) = H^1_0(\Omega)\times L^2(\Omega)$.
Proceeding analogously to Example~\ref{ex:firstorder}, however, yields a perturbation $S_{\Omega,1}$ that includes a first derivative, i.e., that is not bounded on $L^2(\Omega)$.

A remedy is a reformulation in momentum and strain variables $(p,q) = (\rho \partial_t w, \mathcal{T}^{-1}\nabla w)$ that yields the generator
\begin{align*}
    A_{\Omega,2} = \begin{pmatrix}
        0 & \operatorname{div} \\
        \nabla & 0 
    \end{pmatrix}
\end{align*}
on the state space $L^2(\Omega)\times L^2(\Omega;\R^d)$ defined on a suitable subset of $H(\operatorname{div},\Omega)\times H^1(\Omega)$, in which the Dirichlet boundary condition on the displacement may be encoded by means of a unique reconstruction, cf.~\cite[Section 7.2]{reis2024linear}. Here, $H(\operatorname{div},\Omega)$ denotes the space of vector-valued functions for which the weak divergence exists in $L^2(\Omega)$. Importantly, this reformulation via $A_{\Omega,2}$ only includes first-order differential operators such that an analogous computation to Example~\ref{ex:firstorder} yields that $A_{\Omega,2}$ satisfies the assertions (i) and (iii) of Assumption~\ref{Ass: Scaling}.
\end{example}
\Benedikt{
\begin{remark}
If $A_\Omega$ is given by a differential operator of order $r\geq 2$ the perturbation operators $S_{\Omega,i}^\mu$ are not bounded in $X_\Omega$ anymore, as they also involve derivative operators. Consequently, \Cref{Ass: Scaling} is not satisfied in this case. 
In the previous~\Cref{ex:vs} we were able to work around this issue via a change of variables. 
Nonetheless the natural question arises, if it is possible to relax~\Cref{Ass: Scaling} to the case of unbounded perturbation operators $S_{\Omega,i}^\mu: D(S_{\Omega,i}^\mu) \rightarrow X_\Omega$. 
However, this leads to some heavy technical obstacles, the most important among which may be that the operators $A_\Omega + S_{1,\Omega}^\mu$ respectively $(A_\Omega + S_{2,\Omega}^\mu)^*$ are not necessarily generators of strongly continuous semigroups anymore, see e.g.~\cite{Chernoff1968}.
This generator property is, however, strictly necessary for the proof of our main result in order to ensure, that the disturbed optimality system is well-posed when written in scaled variables $\tilde{x}(t):=e^{\mu \n{P-\cdot}_1}x(t)$, $\tilde{u}(t):=e^{\mu \n{P-\cdot}_1}u(t)$, $\tilde{\lambda}(t):=e^{\mu \n{P-\cdot}_1}\lambda(t)$. The question of what kind of perturbation operators $S_\Omega$ the $A_ \Omega+S_\Omega$ remains a generator of a strongly continuous semigroup has been extensively studied in the literature. For bounded perturbations $S \in L(X_\Omega)$ this is always true (see~\cite[Theorem 1.3]{Engel2000}). 
For unbounded perturbation operators, there exists a variety of approaches: If $S_\Omega$ is bounded on $(D(A_\Omega), \n{\cdot}_1)$ where $\n{x}_1:=\n{(sI-A_\Omega)}{X_\Omega}$ for some $s\in \rho (A)$ then $A_\Omega + S_\Omega$ is also an infinitesimal generator~\cite[Corollary 1.5]{Engel2000}. If $S_\Omega$ is relatively $A_\Omega$-bounded, i.e.
\begin{equation*}
	D(A_\Omega) \subset D(S_\Omega) \quad \textrm{and} \quad \exists \, a, b >0 \, \forall x \in D(A_\Omega): \n{S_\Omega x}_{X_\Omega}\leq a\n{A_\Omega x}_{X_\Omega} + b \n{x}_{X_\Omega}
\end{equation*}
then this is also true, but only for contraction~\cite[Theorem 2.7]{Engel2000} and analytic~\cite[Theorem 2.10]{Engel2000} semigroups. Another approach is the Lie-Trotter product formula which comes with a sufficient stability condition on the operators $A_\Omega$ and $S_\Omega$ under which their sum generates a semigroup~\cite{Trotter1959}. Finally, it is possible to make regularity assumptions on the abstract Volterra integral operator associated with $(T_\Omega(t))_{t\geq 0}$~\cite[Theorem 3.1]{Engel2000} respectively its adjoint~\cite[Theorem 3.14]{Engel2000}.
We are confident, that it is possible to relax~\Cref{Ass: Scaling} by using operators $S_{i,\Omega}^\mu $ from a suitable class of unbounded operators inspired by the perturbation results cited above. However, this requires a deep dive into pertubation theory, which goes beyond the scope of this work. Also, as was demonstrated in~\Cref{ex:vs}, most linear-higher order systems can be written as a first-order system by defining new state variables from the spatial derivatives. Therefore, our results can be applied to a wide variety of examples, even though~\Cref{Ass: Scaling} may seem restrictive at first glance.
\end{remark}
}

Having formulated suitable assumptions on the generator in Assumption~\ref{Ass: Scaling} and on the control and observation operator in Assumption~\ref{Ass: ActHom} of the optimal control problem \eqref{HyperbolicOCP}, we may now state our main result of this part.
Therein, we prove the spatial decay of perturbations of the optimality system \eqref{CondensedOptimalitySystem}. In its formulation, we assume a uniform bound on the solution operator that will be verified under domain-uniform stabilizability and detectability assumptions in the subsequent Theorem~\ref{Thm: Boundedness}. The formulation and the proof of both results is similar to the results considering exponential decay in time~\cite{Schaller2020} or spatial decay in elliptic and parabolic equations~\cite{Göttlich2024}.

\begin{thm}
\label{Thm: SensitivityResult}
    Let the Assumptions~\ref{Ass: ActHom} and~\ref{Ass: Scaling} be fulfilled. Assume that there exists a constant $c>0$ such that the solution operator of the optimality system~\eqref{CondensedOptimalitySystem} exists and satisfies the bound
    \begin{equation}
    \label{eq: SolOperatorBoundedness}
        \forall \, T>0 \, \forall \, \Omega \in \mathcal{O}: \n{\M^{-1}}_{L(W^{1\lor2}, W^{2\land \infty})}\leq c.
    \end{equation}
    Let $\mu >0$ be such that
    \begin{equation*}
       \forall \, \Omega \in \mathcal{O}: \n{S_{1,\Omega}^\mu}_{L(X)}+ \n{S_{2,\Omega}^\mu}_{L(X)} \leq \frac{1}{2c}.
    \end{equation*}

Let $\varepsilon \in W_{\R^d}^{1,\infty}$ be a disturbance 
for which the family of restrictions $\left(\varepsilon _\Omega ^T\right)_{\Omega \in \mathcal{O}, T>0}\in (W_\Omega^1)_{\Omega \in \mathcal{O}}$ is exponentially localized in either
 $(F_\Omega)_{\Omega \in \mathcal{O}}=(W_\Omega^1)_{\Omega \in \mathcal{O}}$ or $(F_\Omega)_{\Omega \in \mathcal{O}}=(W_\Omega^2)_{\Omega \in \mathcal{O}}$ with $\n{e^{\mu \n{P-\cdot}}\varepsilon _\Omega  ^T}_{F_\Omega ^T}<C _\varepsilon < \infty$.
Then, there exists a constant $K>0$ such that for all $T>0$ and for all $\Omega \in \mathcal{O}$ the inequality
\begin{equation}
    \label{SpatialDecayTheorem}
    \n{e^{\mu \n{P-\cdot}_1} \delta x_\Omega^T}_{2\land \infty} + 
    \n{e^{\mu \n{P-\cdot}_1} \delta \lambda_\Omega^T}_{2\land \infty} + \n{e^{\mu \n{P-\cdot}_1} \delta u_\Omega^T}_{V_{U_\Omega}} \leq K \n{e^{\mu \n{P-\cdot}_1} \varepsilon_{\textcolor{black}{\Omega}}^T}_{1\lor 2} \leq K\cdot C_\varepsilon
\end{equation}
holds both for $V_{U_\Omega} \!= \!L^\infty (0,T;U_\Omega)$ and $V_{U_\Omega} \! = \! L^2 (0,T;U_\Omega)$.
In particular, the family of states and adjoint states $\left(\delta z _\Omega ^T\right)_{\Omega \in \mathcal{O},T>0} = \left(
    \delta x _\Omega ^T, \,\delta \lambda _\Omega ^T 
\right)_{\Omega \in \mathcal{O},T>0}$ 
is exponentially localized in $\left(W^{2\land \infty}_\Omega\right)_{\Omega \in \mathcal{O}}$. Furthermore, the family of controls $\left(\delta u _\Omega ^T\right)_{\Omega \in \mathcal{O},T>0}$ is exponentially localized in $\left(L^\infty (0,T;U_\Omega)\right)_{\Omega \in \mathcal{O}}$ and $\left(L^2 (0,T;U_\Omega)\right)_{\Omega \in \mathcal{O}}$. 
\end{thm}
\begin{proof}
Let $\Omega \in \mathcal{O}$ and $T>0$ be arbitrary but fixed. In the following we will again leave out the indices $\Omega$ and $T$ for readability. Define $\delta \tilde{x}:=e^{\mu \n{P-\cdot}_1} \delta x $, $\delta \tilde{\lambda}:=e^{\mu \n{P-\cdot}_1} \delta \lambda$ and $\tilde{\varepsilon}:=e^{\mu \n{P-\cdot}_1} \varepsilon$.
Due to~\Cref{Ass: Scaling}, the optimality system \eqref{ErrorOptimalitySystem} implies
\begin{equation*}
\left(\M+\mathcal{F}^\mu \right)
    \underset{=:\delta \tilde{z}}{\underbrace{\begin{pmatrix}
        \delta \tilde{x}\\ \delta \tilde{\lambda}
    \end{pmatrix}}}
    =
    \begin{pmatrix}
        \tilde{\varepsilon} _1\\
        \tilde{\varepsilon} _2\\
        \tilde{\varepsilon} _3\\
        \tilde{\varepsilon} _4
    \end{pmatrix},
\end{equation*}
where $\mathcal{F}^\mu: L^2(0,T;X)^2 \rightarrow  W^2$ acts pointwise in time
\begin{equation*}
    \forall z \in L^2(0,T;X)^2: (\mathcal{F}^\mu z)(t):=\begin{pmatrix}
        0 & -S_2^\mu\\
        0 & 0\\
        -S_1^\mu & 0\\
        0 & 0
    \end{pmatrix}
    z(t).
\end{equation*}
This leads to
\begin{equation*}
    (\M+\mathcal{F}^\mu) \delta \tilde{z} = (I+\mathcal{F}^\mu \M^{-1}) \M \delta \tilde{z} = \tilde{\varepsilon} \,\, \implies \,\, \delta \tilde{z} = \M^{-1} (I+\mathcal{F}^\mu \M^{-1} )^{-1} \tilde{\varepsilon}.
\end{equation*}

We first show the inequality
\begin{equation}
\label{L2Inequality}
    \n{e^{\mu \n{P-\cdot}_1} \delta x}_{L^2(0,T;X)} 
    + \n{e^{\mu \n{P-\cdot}_1} \delta \lambda}_{L^2(0,T;X)} \leq K_2 \n{e^{\mu \n{P-\cdot}_1} \varepsilon}_{W^2}
\end{equation}
for a constant $K_2>0$ which does not depend on $T$ or $\Omega$.
From~\Cref{eq: SolOperatorBoundedness} we know, that there exists a constant $c>0$ independent of $T$ and $\Omega$ such that
    \begin{equation*}
        \n{\M^{-1}}_{L(W^2, L^2(0,T;X)^2)}\leq \n{\M^{-1}}_{L(W^{1\lor 2}, W^{2\land \infty})} \leq c.
    \end{equation*}
~\Cref{Ass: Scaling} (iii) implies for $i\in \{1,2\}$ convergence in the sense of$\n{S_i^\mu }_{L(L^2(\Omega ), L^2(\Omega ))} \overset{\mu \rightarrow 0}{\rightarrow} 0$. This leads to
\begin{flalign*}    
\forall z \in L^2(0,T;X)^2: \n{\mathcal{F}^\mu z}_{W^2}^2
     \leq \left(\n{S_1^\mu }_{L(X)}+\n{S_2^\mu }_{L(X)}\right)^2 \n{z}_{L^2(0,T;X)^2} \overset{\mu \rightarrow 0}{\rightarrow} 0.
\end{flalign*}
Therefore we can choose $\mu >0$ such that 
\begin{equation*}
    \n{\mathcal{F}^\mu}_{L(L^2(0,T;X)^2,W^2)} \leq \frac{1}{2\n{\M^{-1}}_{L(W^{1\lor 2},W^{2\land \infty})}} \leq \frac{1}{2\n{\M^{-1}}_{L(W^2, L^2(0,T;X)^2))}}.
\end{equation*}
Using the Neumann series this implies
\begin{equation*}
    \n{(I+\mathcal{F}^\mu \M^{-1})^{-1}}_{(W^2,W^2)}\leq \sum _{k=0}^\infty \left(\frac{1}{2}\right)^k = 2
\end{equation*}
since $\n{\mathcal{F}^\mu \M^{-1}}_{L(W^2,W^2)}\leq \n{\mathcal{F}^\mu}_{L(L^2(0,T;X)^2,W^2)}\n{\M^{-1}}_{L(W^2, L^2(0,T;X)^2)}\leq \frac{1}{2}$.
Therefore we have
\begin{flalign*}
    \n{\delta \tilde{z}}_{L^2(0,T;X)^2} 
    \!\leq\! \n{\M^{-1}}_{L(W^2, L^2(0,T;X)^2)}
    \n{(I+\mathcal{F}^\mu \M^{-1} )^{-1}}_{L(W^2,W^2)} \n{\tilde{\varepsilon}}_{W^2}\leq K_2\n{\tilde{\varepsilon}}_{W^2}
\end{flalign*}
where $K_2:=2c$. This shows~\eqref{L2Inequality}. For the other estimates we use the inversion formula
\begin{equation*}
    (I+F_\mu \M^{-1})^{-1}=I-(I+F_\mu \M^{-1})^{-1}F_\mu \M^{-1}
\end{equation*}
which gives us
\begin{flalign*}
    &\n{\M^{-1}(I+\mathcal{F}^\mu \M^{-1})^{-1}}_{L(W^{1\lor 2},W^{2\land \infty})}=\n{\M^{-1}-\M^{-1}(I+\mathcal{F}^\mu \M^{-1})^{-1}\mathcal{F}^\mu \M^{-1}}_{L(W^{1\lor 2},W^{2\land \infty})}\\
    \leq &\resizebox{0.97\hsize}{!}{$\,\,\n{\M^{-1}}_{L(W^{1\lor 2},W^{2\land \infty})}+
    \n{\M^{-1}}_{L(W^2,W^{2\land \infty})}
    \n{(I+\mathcal{F}^\mu \M^{-1})^{-1}}_{L(W^2,W^2)}
    \n{\mathcal{F}^\mu}_{L(L^2(0,T;X)^2,W^{2})}
    \n{\M^{-1}}_{L(W^{1\lor 2},L^2(0,T;X)^2)}$}\\
    \leq &\,\,2c.
\end{flalign*}
\noindent The last inequality follows from the particular choice of $\mu$. Overall we have
\begin{flalign*}
    \n{\delta \tilde{z}}_{2\land \infty} \!&=\! \n{\M^{-1} (I+\mathcal{F}^\mu \M^{-1} )^{-1} \tilde{\varepsilon}}_{2\land \infty}\! \leq \! \n{\M^{-1} (I+\mathcal{F}^\mu \M^{-1} )^{-1}}_{L(W^{1\lor 2},W^{2\land \infty})}\n{\tilde{\varepsilon}}_{1\lor 2}\! \leq 2c \! \n{\tilde{\varepsilon}}_{1\lor 2}.
\end{flalign*}
Using $\delta \tilde{u}=(R^*R)^{-1}B^*\delta \tilde{\lambda}$ we find
\begin{equation*}
    \n{\delta \tilde{u}}_{L^\infty (0,T;U)} \leq \n{(R^*R)^{-1}}_{L(U,U)}\n{B^*}_{L(X,U)} \n{\delta \tilde{\lambda}}_{C(0,T;X)} \leq \frac{1}{\alpha} C_B 2c \n{\tilde{\varepsilon}}_{1\lor 2}
\end{equation*}
and
\begin{equation*}
    \n{\delta \tilde{u}}_{L^2 (0,T;U)} \leq \n{(R^*R)^{-1}}_{L(U,U)}\n{B^*}_{L(X,U)} \n{\delta \tilde{\lambda}}_{L^2 (0,T;X)} \leq \frac{1}{\alpha} C_B 2c \n{\tilde{\varepsilon}}_{1\lor 2}.
\end{equation*}
For both $V_{U_\Omega} = L^\infty (0,T;U_\Omega)$ and $V_{U_\Omega} = L^2 (0,T;U_\Omega)$ this leads to
\begin{equation*}
    \n{e^{\mu \n{z-\cdot}_1} \delta x_\Omega^T}_{2\land \infty} + 
    \n{e^{\mu \n{z-\cdot}_1} \delta \lambda_\Omega^T}_{2\land \infty} + \n{e^{\mu \n{z-\cdot}_1} \delta u_\Omega^T}_{V_{U_\Omega}} \leq \left(1+\frac{1}{\alpha}C_B\right)2c\n{\tilde{\varepsilon}}_{1\lor 2}=:K\n{\tilde{\varepsilon}}_{1\lor 2}.
\end{equation*}
Since the constant $K$ does not depend on $\Omega$ or $T$ this shows the exponential localization of $\left(\delta z _\Omega ^T\right)_{\Omega \in \mathcal{O},T>0}$ in $W^{2\land \infty}_\mathcal{O}$.
\end{proof}

\Benedikt{
\begin{remark}
    Note that the meaning of the exponential decay from~\Cref{Thm: SensitivityResult} may depend very much on the state space formulation chosen for the problem under consideration. 
    To demonstrate this, we revisit the wave equation in~\Cref{ex:vs} for $d=1$ on $\Omega _L:=[0,L]$, i.e, we have the boundary conditions $w(0,t)=w(L,t)=0$. 
    Let $\mathcal{O}:=\{[0,L]: L>0\}$. 
    Consider the formulation in deflection variables
    \begin{equation*}
        \dot{x}(t)=A_{\Omega,1}x(t),\quad x(0) = x_{\Omega, 1}
    \end{equation*}
    and in strain variables
    \begin{equation*}
        \dot{x}(t)=A_{\Omega,2}x(t), \quad x(0) = x_{\Omega, 2}
    \end{equation*}
    both with Dirichlet boundary conditions $(x_1(t))(0) = (x_1(t))(L) = 0$. Let $x_L^i, i\in \{1,2\}$ be classical (in particular continuous in space) solutions of these two formulations.
    We write $x_L^i(\omega ,t)$ instead of $(x_L^i(t))(\omega )$ where $\omega$ is the spatial variable of the solution. Assume domain-uniform exponential decay of the state of this equation around some point in space, i.e., there exist $P \in \R$ and $M,\mu>0$ such that
    \begin{equation}
    \label{eq: WaveDecay}
        \forall L>0: \n{x_L^i(\omega,t)}\leq M e^{-\mu \n{P-\omega}_1}.
    \end{equation}
    For $i=1$ this corresponds to the decay of the deflection and its time derivative. For $i=2$, \eqref{eq: WaveDecay} is equivalent to the decay of the spatial and time derivative of the deflection. For Dirichlet boundary conditions, the second decay property implies the first since
    \begin{flalign*}
        \forall \omega \in [0,P]: \n{w(\omega,t)} &= \n{\int _0^\omega \frac{\partial}{\partial s}w(s,t)\mathrm{d}s} \leq \int _0^\omega M e^{-\mu \n{P-s}_1}\mathrm{d}s \\
        &= \frac{M}{\mu} \left(e^{-\mu \n{P}_1}-e^{-\mu \n{P-\omega}_1}\right) \leq \frac{M}{\mu}\\
        \forall \omega \in (P,L]: \n{w(\omega,t)} &\leq 
        \int _\omega ^L M e^{-\mu \n{P-s}_1}\mathrm{d}s = \frac{M}{\mu} \left(e^{-\mu \n{P-\omega}_1}-e^{-\mu \n{P-L}_1}\right) \leq \frac{M}{\mu}e^{-\mu \n{P-\omega}_1}.
    \end{flalign*}
    For Neumann boundary conditions $\frac{\partial}{\partial \omega}w(0,t)=\frac{\partial}{\partial \omega}w(L,t)=0$ this implication does not hold anymore. For example a constant initial value $x_{\Omega, 1}\equiv (c,0), c>0$ corresponding to $x_{\Omega, 2}\equiv (0,0), c>0$ leads to the constant solutions $x_{\Omega,1}(\omega, t)=(c,0)$, $x_{\Omega,2}(\omega, t)=(0,0)$. Obviously, the second solution is exponentially decaying while the first does not.
\end{remark}
}
\Cref{Thm: SensitivityResult} shows that exponential localization of the family of perturbations implies exponential localization of the error variables. To this end we assumed, that the solution operator $\M^{-1}$ is bounded from $W^{1\lor2}$ to $W^{2\land\infty}$. In the following, we show, that such a boundedness can be achieved under a domain-uniform stabilizability/detectability assumption. 

\begin{definition}[Domain-uniform exponential stability/stabilizability/detectability]\hfill
\label{Def: DomUniStabDet}
\begin{itemize}
    \item[(i)] We call a family $((\mathcal{T}_\Omega(t))_{t\geq 0})_{\Omega \in \mathcal{O}}$ of strongly continuous semigroups domain-uniformly exponentially stable in $\mathcal{O}$, if and only if there exist constants $M,k>0$ such that
\begin{equation*}
    \forall \, \Omega\in \mathcal{O} \, \forall \, t \geq 0: \n{\mathcal{T}_\Omega(t)}_{L(X_\Omega)} \leq Me^{-kt}.
\end{equation*}
    \item[(ii)] We call a pair $(A_\Omega, B_\Omega)_{\Omega \in \mathcal{O}}$ with $A_\Omega: X_\Omega \supset D(A_\Omega) \rightarrow X_\Omega$ and $B_\Omega \in L(U_\Omega, X_\Omega )$, $\Omega \in \mathcal{O}$
 domain-uniformly exponentially stabilizable in $\mathcal{O}$ if and only if there exists a family of \Benedikt{uniformly (in $\Omega \in \mathcal{O}$) bounded} feedback operators $(K_\Omega^B)_{\Omega \in \mathcal{O}}$, $K_\Omega^B \in L(X_ \Omega,U_\Omega)$, $\Omega\in \mathcal{O}$, such that $\left(A_\Omega + B_\Omega K_\Omega^B\right)_{\Omega\in \mathcal{O}}$ generates a domain-uniformly exponentially stable family of semigroups in $\mathcal{O}$.
    \item[(iii)] We call $(A_\Omega, C_\Omega)_{\Omega \in \mathcal{O}}$, $C_\Omega \in L(X_\Omega, Y_\Omega )$, $\Omega\in \mathcal{O}$ domain-uniformly exponentially detectable in $\mathcal{O}$ if and only if $(A_\Omega^*, C_\Omega^*)$ is domain-uniformly exponentially stabilizable.
\end{itemize}
\end{definition}

\noindent The key property of this domain-uniform stabilizability/detectability property is, to make sure, that the system under consideration can be uniformly stabilized via the optimal control and that it is possible to observe unstable dynamics inside the system via the cost functional. 
In an unstable system small perturbations may cause very large changes in the systems behaviour which means that exponential localization of the disturbance would not necessarily imply exponential localization of the error variables anymore. 

\begin{thm}
\label{Thm: Boundedness}
    Let $(A_\Omega, B_\Omega)_{\Omega \in \mathcal{O}}$ and $(A_\Omega, C_\Omega)_{\Omega \in \mathcal{O}}$ be domain-uniformly stabilizable respectively detectable families in the sense of~\Cref{Def: DomUniStabDet}. 
    Then there exists a constant $c>0$ such that the norm of the solution operator $\M^{-1}$ can be estimated by
    \begin{equation}
        \label{Boundedness}
        \forall\, T>0 \, \forall \, \Omega \in \mathcal{O}: \n{\M^{-1}}_{L(W^{1\lor 2},W^{2\land \infty})}\leq c.
    \end{equation}
\end{thm}
\begin{proof}
    The proof follows along the lines of~\cite[Theorem 10]{Schaller2020} and is provided for completeness in~\Cref{App: Boundedness}.
\end{proof}

\section{Domain-uniform stabilizability of the transport equation on a scalar domain}
\label{Sec: TransportConstant} 
In the previous section, domain-uniform stabilizability and detectability of the underlying operator fa\-milies $(A_\Omega, B_\Omega)_{\Omega \in \mathcal{O}}$ and $(A_\Omega, C_\Omega)_{\Omega \in \mathcal{O}}$ has been the main assumption we needed, to show that exponential localization of the perturbation implies the same for the error variables (see~\Cref{Thm: Boundedness} and~\Cref{Thm: SensitivityResult}). 
In the present and the next section we will present a characterization for this assumption for a controlled transport equation with periodic boundary conditions (see Theorems~\ref{Thm: Stabilizability} and~\ref{Thm: StabilizabilityNonConstantCase}). Thereby we show that domain-uniform stabilizabili\-ty/detectability is indeed a reasonable assumption which is fulfilled for a relevant example if the control domain fulfills some mild requirements. We note that, in the following we only consider the case of stabilizability, and provide a result for detectability using duality arguments in the proof of~\Cref{Cor: TransportDecay}.

We stress that stabilizability of hyperbolic equations on scalar domains is a well-studied subject, see e.g.\ the monograph \cite{bastin2016stability}, and many aspects of the subsequent deductions are standard tools in the study of hyperbolic equations. However, in this work, we are particularly interested in \textit{domain-uniform} stabilizability being the central ingredient of the exponential decay estimates, such that we meticulously will track the dependence of all involved constants on the domain size.

\noindent In this section we will always consider the family of domains $\mathcal{O}:=\{\Omega _L:=[0,L]: L>0\}$. On such a domain $\Omega_L \in \mathcal{O}$, the controlled transport equation with periodic boundary conditions is given by
\begin{subequations}
\label{eq: Transport}
\begin{align}
\label{TransportEquationControlled}
    \forall \, (\omega,t) \in \Omega_L \times [0,T]: \frac{\partial}{\partial t}{x}(\omega,t)&=-c(\omega)\frac{\partial}{\partial \omega}x(\omega,t) + \chi _{\Omega_L^c}(\omega) u(\omega,t)\\
    \label{BoundaryCondition}
    \forall \, t \in [0,T]: x(0,t) &= x(L,t) \\
    \label{InitialCondition}
    \forall \, \omega \in \Omega_L: x(\omega,0) &= x_{\Omega _L}^0
\end{align}
\end{subequations}
with time horizon $T>0$, an initial distribution $x_{\Omega _L}^0\in L^2(\Omega_L)=:X_{\Omega_L}$ and a control $u \in L^2(\Omega _L^c\times [0,T])=:U_{\Omega _L}$. For brevity we often write $x^0$ instead of $x_{\Omega _L}^0$ in the following. The control domain $\Omega_c \subset \R_{\geq 0}$ is assumed to be Lebesgue-measurable and to have positive measure. For $L>0$ we define $\Omega _L^c := \Omega_c \cap [0,L]$. By $\chi _{\Omega_L^c}$ we denote the characteristic function of $\Omega _L^c$.
In this section we will only consider constant transport velocities $c>0$. The case of space-dependent transport velocities will be treated in~\Cref{Sec: TransportNonConstant}.
The operator describing the evolution of \eqref{eq: Transport} is given by
\begin{equation}
\label{Eq: TransportOperator}
     A_{\Omega _L}: D(A_{\Omega _L}):= \{x \in H^1(\Omega _L): x(0)=x(L)\}\subset L^2(\Omega_L)  \rightarrow L^2(\Omega _L), \platz A_{\Omega _L}x = -c \frac{\partial}{\partial \omega}x.
 \end{equation}
 Note that the boundary evaluations in $D(A_{\Omega _L})$ are well-defined as weakly differentiable functions on scalar domains are absolutely continuous. The controlled transport equation with periodic boundary conditions~\eqref{eq: Transport} leads to the inhomogeneous Cauchy problem
\begin{equation}
\label{eq: CauchyTransport}
    \dot{x}(t)= A_{\Omega _L}x(t) + B_{\Omega _L}u(t), \, x(0) = x_{\Omega _L}^0,
\end{equation}    
where the input operator $B_{\Omega _L}$ is defined by
\begin{equation}
\label{eq: ControlOperator}
     B_{\Omega _L}: L^2(\Omega _L^c) \rightarrow L^2(\Omega _L), \platz (B_{\Omega _L}v)(\omega)=\left\{ \begin{array}{cc}
      v(\omega), & \omega \in \Omega _L^c\\
      0,    & \mathrm{else}
     \end{array}
     \right..
\end{equation}
We also define the output operator
\begin{equation}
    \label{eq: OutputOperator}
    C_{\Omega _L}: L^2(\Omega _L) \rightarrow L^2(\Omega _L^o), \platz C_{\Omega _L}v = v|_{\Omega _L^o},
\end{equation}
where the measurement domain $\Omega_o \subset \R_{\geq 0}$ is assumed to be Lebesgue-measurable with positive measure and for $L>0$ we define $\Omega _L^o := \Omega_o \cap [0,L]$. 
For $\Omega _L^o=\Omega _L^c$ we have $B_{\Omega _L}^* = C_{\Omega _L}$. The theorem of Stone~\cite[Theorem 3.8.6]{Weiss2009} implies that the operator $A_{\Omega _L}$ generates a unitary group since it is skew-adjoint. Note that due to duality domain-uniform detectability of $(A,C)$ is the same as domain-uniform stabilizability of $(A^*,C^*)$.

 To analyze the domain-uniform exponential stabilizability of the transport equation we will use solution formulas as presented in the following.
\begin{lemma} 
    \label{Lem: SolConstant}
    For the Cauchy problem~\eqref{eq: CauchyTransport} the following three statements hold:
    \begin{itemize}
        \item [(i)] The mild solution in the uncontrolled case, i.e. $u\equiv 0$, is given by
        \begin{equation*}
            \forall t\in [0,T]: x_{x_0}^L(\cdot , t) = P_{\Omega _L}(x^0)(\cdot - ct),
        \end{equation*}
        where $P_{\Omega _L}: L^2(\Omega _L,\R) \rightarrow L^2((-\infty, L],\R)$ defined by
    \begin{equation*}
        \forall k \in \N: P_{\Omega _L}(x^0) (\omega) = x^0(\omega+(k-1)L) \, \textrm{for a.a. }\omega \in ((1-k)L, (2-k)L]
    \end{equation*}
    is the periodization operator with period $L>0$.
        \item [(ii)] The operator $A_{\Omega _L}$ generates the unitary group $(\mathcal{T}_{\Omega _L}(t))_{t\in \R}$ with
        \begin{equation*}
        \forall t \in \R: \mathcal{T}_{\Omega _L}(t): L^2(\Omega _L) \rightarrow L^2(\Omega _L), \platz \mathcal{T}_{\Omega _L}(t)x^0 := P_{\Omega _L}(x^0)(\cdot - ct).
        \end{equation*}
        \item[(iii)] For all $t \in [0,T]$ the solution in the controlled case is given by
        \begin{flalign*}
            x_{x^0,u}^L(t) &= \mathcal{T}_{\Omega _L}(t)x_0 + \int _0^t \mathcal{T}_{\Omega _L}(t-\tau)B_{\Omega _L}u(\tau)\mathrm{d}\tau\\
            &= P_{\Omega _L}(x^0)(\cdot - ct) + \int _0^t P_{\Omega _L}(B_{\Omega _L}u(\tau))(\cdot -c (t-\tau)) \mathrm{d}\tau.
        \end{flalign*}
    \end{itemize}
\end{lemma}
\begin{proof}
   It is easy to check, that, for $x^0_{\Omega_L}\in D(A_{\Omega_L})$, the formula given in (i) is indeed a solution of the uncontrolled transport equation~\eqref{eq: Transport} and hence a classical solution of the uncontrolled Cauchy problem. Correspondingly, (i) and (ii) follow by a standard density argument in combination with uniqueness of solutions as we already know that $A_{\Omega _L}$ generates a strongly continuous semigroup. 
    The last formula (iii) directly follows from the variation of constants formula~\cite[Corollary 1.7]{Engel2000}.
\end{proof}

\subsection{Motivational examples and negative results}
\label{Sec: NegativeResults}
In this part we will discuss some instructive examples before proving our main result in the subsequent subsection in~\Cref{Thm: Stabilizability}.
We show, that the family $(A_{\Omega _L},B_{\Omega _L})_{\Omega_L \in \mathcal{O}}$ is not domain-uniformly exponentially stabilizable, if the control can only be applied on a single finite interval, i.e. $\Omega_c=[a,b]$ for some $0\leq a < b < \infty$. Finally we show, that both, in case of a control on the full domain $\Omega_c = [0,L]$ as well as in case of a control domain consisting of equidistantly distributed intervals of equal size the condition of domain-uniform stabilizability is fulfilled.
 
To show that it is not stabilizable in case of a single finite control interval we will use the finite propagation velocity $c>0$ which implies that a control on an interval $[a,b]$ cannot instantaneously influence the whole length of the domain $[0,L]$. This is the essence of the following result.
\begin{lemma}
    \label{LemmaFiniteSpeed}
    Consider the controlled transport equation~\eqref{TransportEquationControlled} with control domain $\Omega_c = [a,b], 0 \leq a < b < \infty$ and domain size $L>b$. For $t \in \left[0,\frac{L-b}{c}\right]$ the mild solution is given by
    \begin{equation}
    \label{xRewrittenLocal}
        x_{x_0}^u(\omega,t)=\left\{ \begin{array}{cc}
        P_{\Omega _L}(x^0)(\omega-ct)+ \int _0^t P_{\Omega _L}(B_{\Omega _L}u(\tau))(\omega -c (t-\tau)) \mathrm{d}\tau, & \omega \in [a,b+ct] \\
        P_{\Omega _L}(x^0)(\omega-ct),    & \mathrm{else}
    \end{array}
    \right..   
    \end{equation}
    In particular, for any time $t \in \left[0,\frac{L-b}{c}\right]$ the solution only depends on the control $u$ on the interval $[a, b+ct]$.
\end{lemma}

\begin{proof}
    First we consider the integral term
    \begin{equation}
    \label{IntegralTerm}
        I_L: [0,L]\times \left[0,\frac{L-b}{c}\right] \to \mathbb{R}, \platz I_L(\omega,t):=\int _0^t P_{\Omega _L}(B_{\Omega _L}u(\tau))(\omega -c (t-\tau)) \mathrm{d}\tau
    \end{equation}
    from the solution formula in~\Cref{Lem: SolConstant} (iii). Define a mapping $\alpha_{\omega,t}: \R \rightarrow \R,\, \alpha_{\omega,t} (\tau ):= \omega-c(t-\tau)$.
    Note that for all $(\omega,t) \in \mathrm{dom} (I_L)$ and for all $\tau \in [0,t]$ we have $\alpha _{\omega,t}(\tau) \in [b-L,L]$. Since $P_{\Omega _L}(B_{\Omega _L}u(\tau))(y)=0$ for $y\in (b-L,0)$ and $P_{\Omega _L}(B_{\Omega _L}u(\tau))(y)=(B_{\Omega _L}u(\tau))(y)$ for $y\in [0,L]$ we find the equality
    \begin{flalign*}
        I_L(\omega,t) &= \int _0^t P_{\Omega _L}(B_{\Omega _L}u(\tau))(\alpha _{\omega,t} (\tau )) \mathrm{d}\tau
        = \frac{1}{c}\int _{\omega-ct}^\omega P_{\Omega _L}(B_{\Omega _L}u(\tau))(y) \mathrm{d}y\\
        &= \frac{1}{c} \int _{\mathrm{max}(0,\omega-ct)}^\omega (B_{\Omega _L}u(\tau))(y) \mathrm{d}y.
    \end{flalign*}
    Since for any $\tau \geq 0$ the map $B_{\Omega _L}\,u(\tau)$ only takes non-zero values in $[a,b]$, i.e. the support of $B_{\Omega _L}u(\tau)$ is a subset of $[a,b]$ we can further rewrite $I_L(\omega,t)$ into the form
    \begin{equation*}
        I_L(\omega,t) = \int _{\max (0,\frac{a-\omega}{c}+t)}^{\min (\frac{b-\omega}{c}+t,t)} (B_{\Omega _L}u(\tau))(\omega -c (t-\tau))\mathrm{d}\tau.
    \end{equation*}
    Therefore, it vanishes if $b+ct\leq \omega$ or $a \geq \omega$.
    This shows that for arbitrary but fixed $t \in \left[0,\frac{L-b}{c}\right]$ the solution $x_{x_0}^u(\omega,t)$ does only depend on the control on the interval $[a,b+ct]$, i.e. it can be rewritten as in~\eqref{xRewrittenLocal}.
\end{proof}
\noindent Since the operators $A_{\Omega _L}$ generate unitary (and therefore norm-preserving) groups $(\mathcal{T}_{\Omega _L})_{\Omega _L\in \mathcal{O}}$ the uncontrolled transport equation is not exponentially stable and thus also not domain-uniformly exponentially stable. Using Euclidean division to find a representation $ct=mL+L_0, m\in \N, L_0 \in [0,L)$ this can also be seen via the equations
\begin{flalign}
\label{Eq: NormPreservation}
\begin{split}
    \n{\mathcal{T}_{\Omega _L}(t)x^0}_{L^2([0,L])}^2 &= \int _0^L \n{P_{\Omega _L}(x^0)(\omega - ct)}^2 \mathrm{d}\omega = \int _{-ct}^{L-ct} \n{P_{\Omega _L}(x^0)(y)}^2 \mathrm{d}y\\
    &=\int _{-L_0}^{0} \n{P_{\Omega _L}(x^0)(y)}^2 \mathrm{d}y + \int _{0}^{L-L_0} \n{P_{\Omega _L}(x^0)(y)}^2 \mathrm{d}y\\
    &=\int _{L-L_0}^{L} \n{x^0(y)}^2 \mathrm{d}y + \int _{0}^{L-L_0} \n{x^0(y)}^2 \mathrm{d}y = \n{x^0}_{L^2([0,L])}^2
\end{split}
\end{flalign}
showing that the periodization operator preserves the $L^2(\Omega _L)$-norm.
This observations lead to the following negative result.
\begin{proposition}
\label{Prop: NegativeResult}
    Consider the transport equation~\eqref{TransportEquationControlled} with $\Omega_c = [a,b], 0 \leq a < b < \infty$. 
    Then, the family $(A_{\Omega _L},B_{\Omega _L})_{\Omega_L \in \mathcal{O}}$ is \textit{not} domain-uniformly stabilizable.
\end{proposition}
\begin{proof}
    For each $L>0$ let $K_{\Omega _L}^B\in L(X_{\Omega _L}, L^2(\Omega_c^L))$ be a given feedback and denote by $\mathcal{T}_{L}^\varphi$ the closed-loop semigroup generated by $A_{\Omega _L}+B_{\Omega _L}K_{\Omega _L}^B$. 
    Let $(L_k)_{k\in \N} \subset (b,\infty)^\N$ be a sequence of domain sizes which is defined by $\forall k\in \N: L_k:=b+3 ck$ such that $k \in \left[0,\frac{L_k-b}{c}\right]=\left[0,3 k\right]$. Furthermore, for each $k\in \N$ let $g_k \in L^2(\Omega _{L_k})$ be such that $\mathrm{supp}(g_k) \subset (b+c k ,b+2c k]$.
    Then $\mathcal{T}_{L_k}^\varphi g_k$ solves the controlled transport equation with domain $\Omega _{L_k}$, initial value $g_k$ and control 
    \begin{equation*}
        u_k(\omega,t)=K^{\Omega_{L_k}}_B(\mathcal{T}_{L_k}^\varphi (t) g_k)(\omega)=(K^{\Omega_{L_k}}_B x_{x_0}^{u_k}(\cdot ,t))(\omega).
    \end{equation*}
    Following~\Cref{LemmaFiniteSpeed} we therefore know, that it can represented as in~\eqref{xRewrittenLocal}. This implies
    \begin{flalign*}
        \forall k \in \N: \n{\mathcal{T}_{L_k}^\varphi(k)g_k}_{L^2([0,L_k])}^2 &= \int _{0}^{L_k}\n{(\mathcal{T}_{\Omega_{L_k}}^\varphi(k)g_k)(\omega)}^2 \mathrm{d}\omega = \int _{0}^{L_k}\n{x_{g_k}^{u_k}(\omega,k)}^2 \mathrm{d}\omega\\
        &\geq \int _{b+2ck}^{L_k}\n{x_{g_k}^{u_k}(\omega,k)}^2 \mathrm{d}\omega
        \overset{\eqref{xRewrittenLocal}}{=} \int _{b+2ck}^{L_k}\n{P_{L_k}(g_k)(\omega-ck)}^2 \mathrm{d}\omega\\
        &= \int _{b+ck}^{b+2ck}\n{P_{L_k}(g_k)(y)}^2 \mathrm{d}y = \n{g_k}_{L^2([0,L_k])}^2.
    \end{flalign*}
    Therefore we have
        $\forall \, k\in \N: \n{\mathcal{T}_{\Omega_{L_k}}^\varphi(k)}\geq 1.$
\end{proof}

\noindent\Cref{Prop: NegativeResult} shows that a domain-uniform stabilization of the transport equation is not possible using a single bounded interval as control domain. In the following~\Cref{Ex: LocalDamping} we will further discuss this situation for the case of a constant state feedback.

\begin{example}[Transport equation with local damping via state feedback]
\label{Ex: LocalDamping}
    Consider a damped transport equation
    \begin{equation}
        \label{DampedTransportEquation}
        \forall (\omega,t) \in \Omega _L \times [0,T]: \dot{x}(\omega,t)=-k\chi _{[a,b]} x(\omega,t)-cx'(\omega,t)
    \end{equation}
    with damping constant $k>0$, boundary condition~\eqref{BoundaryCondition} and initial condition~\eqref{InitialCondition} and $0<a<b<\infty$. 
    It is easy to check that the solution is given by
    \begin{equation*}
        x(\omega,t)=e^{-\frac{k}{c}E_L(\omega,t)}P_{\Omega _L}(x^0)(\omega-ct) \quad
        \mathrm{with}
        \quad
        E_L(\omega,t):=\int _{\omega-ct}^\omega P_{\Omega _L}(\chi _{[a,b]})(y) \mathrm{d}y.
    \end{equation*}
    Further, we may rewrite~\eqref{DampedTransportEquation} via
    \begin{equation}
    \label{ClosedLoop}
        \dot{x}=(A_{\Omega _L}+B_{\Omega _L}K^B_{\Omega _L})x, \qquad x(0)=x_{\Omega _L}^0
    \end{equation}
    with feedback operator 
        $K_{\Omega _L} ^B\in L(L^2(\Omega _L), L^2(\Omega _L ^c))$ defined via $ K^B_{\Omega _L}x:= -k\chi _{[a,b]}x$.
    In view of the above solution formula, the closed-loop semigroup $\mathcal{T}^\varphi_{\Omega_L}$ generated by $A_{\Omega _L}+B_{\Omega _L}K^B_{\Omega _L}$ satisfies
    \begin{flalign*}
        &\n{\mathcal{T}^\varphi_{\Omega_L}(t)x^0}_{L^2(\Omega_L)}^2 = \int _0^L \n{e^{-\frac{k}{c}E_L(\omega,t)}P_{\Omega _L}(x^0)(\omega-ct)}^2 \mathrm{d}\omega\\
        \leq &e^{-\frac{2k}{c}\left \lfloor \frac{ct}{L}\right\rfloor(b-a)}\int _{-ct}^{L-ct} \n{P_{\Omega _L}(x^0)(y)}^2 \mathrm{d}y 
        = e^{-\frac{2k}{c}\left \lfloor \frac{ct}{L}\right\rfloor(b-a)} \n{x^0}_{L^2(\Omega _L)}^2
    \end{flalign*}
    such that we have the exponential decay
    \begin{equation*}
        \n{\mathcal{T}^\varphi_{\Omega_L}(t)}_{L(L^2([0,L]),L^2([0,L]))}^2 \leq e^{-\frac{2k}{c}\left \lfloor \frac{ct}{L}\right\rfloor(b-a)}.
    \end{equation*}
    The crucial observation now is that decay rate of this estimate depends on the domain size $L$ and in particular deteriorates for $L\to \infty$, hence preventing domain-uniform stabilization. 
    However, if we choose the control domain $\Omega _c = [0,\infty)$ and $\Omega_L^c=[0,L]$ for all $L>0$, the corresponding state feedback leads to an domain-uniformly exponentially stable semigroup.
    By straightforward calculation, one observes that the associated closed-loop semigroup satisfies $(\mathcal{T}^\varphi_{\Omega_L}(t)x^0)(\omega) := e^{-kt} P_{\Omega _L}(x^0)(\omega - ct)$ such that
    \begin{flalign*}
        \n{\mathcal{T}^\varphi_{\Omega_L}(t)x^0}_{L^2([0,L])}^2 &= \int _0^L \n{e^{-kt}P_{\Omega _L}(x^0)(\omega - ct)}^2 \mathrm{d}\omega = e^{-2kt}\int _{-ct}^{L-ct} \n{P_{\Omega _L}(x^0)(y)}^2 \mathrm{d}y\\ 
        &=e^{-2kt}\left(\int _{L-L_0}^{L} \n{x^0(y)}^2 \mathrm{d}y + \int _{0}^{L-L_0} \n{x^0(y)}^2 \mathrm{d}y\right)\\ 
        &= e^{-2kt}\int _{0}^{L} \n{x^0(y)}^2 \mathrm{d}y = e^{-2kt} \n{x^0}_{L^2([0,L])}^2,
    \end{flalign*}
    where we used Euclidean division to find a representation $ct=mL+L_0, m\in \N, L_0 \in [0,L)$. This implies
    \begin{equation*}
        \n{\mathcal{T}^\varphi_{\Omega_L}(t)}_{L(L^2([0,L]),L^2([0,L]))}^2 = e^{-kt}.
    \end{equation*}
\end{example}

\noindent In the first case of previous~\Cref{Ex: LocalDamping} stabilization of the transport equation via state feedback failed because of the direct dependency of the decay rate on the domain size $L$. In the following example we eliminate this dependency. This is achieved by ensuring that the maximum distance between an arbitrary point in the spatial domain and the control domain is uniformly bounded by some constant $L_0 > 0$.

\begin{example}[Transport equation with damping on equidistantly distributed intervals]
\label{Ex: Equidistant}
    We will show that if we replace the single control interval by a sequence of equidistant control intervals, then the family $(A_{\Omega _L},B_{\Omega _L})_{\Omega_L \in \mathcal{O}}$ is domain-uniformly exponentially stabilizable.
    More precisely, let $L_0>b>a>0$ be given and for control domain $\Omega_c := \overset{\infty}{\underset{k=0}{\cup}}\left[a+kL_0,b+kL_0\right]$ and $L\geq b$ define the feedback operator $K^B_{\Omega _L}\in  L(L^2(\Omega _L), L^2(\Omega _L^c))$ via
    \begin{equation*}
          K^B_{\Omega _L}x:= -k\chi _{\Omega _L^c}x \quad
        \mathrm{with} \quad
        \Omega _L^c:= \Omega_c  \cap [0,L].
    \end{equation*}
    Similar deliberations as in~\Cref{Ex: LocalDamping} yield, that the operator $A_{\Omega _L}+B_{\Omega _L}K^B_{\Omega _L}$ generates a strongly continuous semigroup $\left(\mathcal{T}^\varphi_{\Omega_L}(t)\right)_{t\geq 0}$ which fulfils the estimate
    \begin{equation*}
        \n{\mathcal{T}^\varphi_{\Omega_L}(t)}_{L(L^2([0,L]),L^2([0,L]))}^2 \leq e^{-\frac{k}{c}\left \lfloor \frac{ct}{L_0}\right\rfloor(b-a)}.
    \end{equation*}
    Here, the decay rate of this estimate only depends on the fixed constant $L_0$ and not on the domain size $L$ such that the closed-loop semigroup is domain-uniformly exponentially stable.
\end{example}
\noindent Examples~\ref{Ex: LocalDamping}-\ref{Ex: Equidistant} show that the domain-uniform stabilizability of the transport equation is strongly affected by the control domain. In the next section we present a general characterization of control domains which allow for this property.

\subsection{Characterization of 
control domains for domain-uniform stabilizability}
\label{Sec: CharacterisationConstant}
In this section we present our main result characterizing the domain-uniform stabilizability for the transport equation with constant velocity $c>0$. A similar result for space-dependent velocities is presented in~\eqref{Sec: TransportNonConstant}. We consider linear and bounded state feedbacks $K_{\Omega _L}^B\in L(L^2([0,L]), L^2(\Omega _L^c))$ of the form
\begin{equation}
\label{GeneralStateFeedback}
    K_{\Omega _L}^Bv:=-k_Bv, \quad k_B \in L^\infty (\R _{>0})
\end{equation}
on a control domain $\Omega _c \subset \R_ {\geq 0}$ which is given as a union of countably many intervals.
For such state feedbacks we find the following necessary and sufficient condition for domain-uniform stabilizability.
\begin{thm}[Domain-uniform stabilizability of the transport equation]
\label{Thm: Stabilizability}
The following three statements are equiva\-lent:
\begin{itemize}
    \item [(i)] The controlled transport equation~\eqref{eq: Transport} with constant transport velocity $c>0$ can be domain-uniformly stabilized via a state feedback of the form~\eqref{GeneralStateFeedback}.
    \item [(ii)] There exist constants $K>k>0$ such that for all $n,m \in \N$ with $n\geq m$
    \begin{equation}
    \label{IntervalConditionTheorem}
        k(a_n-b_m) - \sum _{j=m+1}^{n-1} K (b_j-a_j) \leq 1, 
    \end{equation}
    where $\Omega _c:=\underset{j \in \N}{\cup} [a_j,b_j]$, $\left(a_j\right) _{j\in \N}$ is an unbounded sequence with $a_0 = 0$, $a_i < b_i\leq a_{i+1}$ for $i \in \mathbb{N}$.
    \item [(iii)] There exist constants $c_0, c_1>0$ such that for all intervals $I\subset \R _{\geq 0}$ the inequality
    \begin{equation*}
        |\Omega _c \cap I| \geq c_1 |I| - c_0
    \end{equation*}
    is fulfilled and $\forall \, L>0: |\Omega _L^c| = |\Omega _L \cap \Omega _c| >0$.
\end{itemize}
\end{thm}
\noindent \Cref{Thm: Stabilizability} (ii) provides a simple algebraic characterization of control domains which ensure domain-uniform stabilizability of~\eqref{eq: Transport}.
This characterization has two main advantages: 
First, it is a useful tool in order to check domain-uniform stabilizability for a given control domain. 
Second, it even allows for the iterative construction of such a control domain. For example, this can be done by predefining constants $K>k>0$ and - starting from $a_1=0$ - successively choosing $a_n$ and $b_n$ such that condition (ii) is fulfilled for all $m\leq n$.
In (iii) this characterization is rephrased to allow for an interpretation in terms of measures. This characterization shows that the size of the control domain $|\Omega _c|$  has to grow at least linearly with the domain size $|\Omega|$ in order to guarantee domain-uniform stabilizability. In~\cite[Eq.(3.18)]{Göttlich2024} a similar condition is used for this purpose. However in this work the authors assume that the uncontrollable domain, i.e. $\Omega _L \setminus \Omega_L^c$, has $L$-uniformly bounded Lebesgue measure. Here, we stress that condition (iii) in~\Cref{Thm: Stabilizability} actually allows for the set $\R _{\geq} \setminus \Omega _c$ to be unbounded as long as the Lebesgue measure $|\Omega _L \setminus \Omega_L^c|$ does not grow more than linearly with the domain size $L$. In this sense~\Cref{Thm: Stabilizability}(iii) is more general than the previous result of \cite{Göttlich2024}.

\begin{remark}
    Note that it is strictly necessary for domain-uniform stabilizability/detectability of the transport equation, that the Lebesgue measure $|\Omega _L^c|$ of the control domain increases with the domain size $L$. To see this, assume that there is $c_\mathrm{Dom}>0$ such that for all $L>0$ we have $\mu (\Omega _L^c)<c_\mathrm{Dom}$. However, condition (ii) from~\Cref{Thm: Stabilizability} reads
\begin{equation*}
        1 \geq k(a_n-b_m) - K\sum _{j=m+1}^{n-1} (b_j-a_j) > k(a_n-b_m) - K c_\mathrm{Dom}.
\end{equation*}
Since the first term on the right-hand side diverges for fixed $m \in \N$ as $n\to \infty$ while the second one is bounded, this yields a contradiction. 
\Benedikt{Generally speaking the reasons for this requirement are the finite propagation velocity of the transport equation and the uniform boundedness of the state feedback operators $K_\Omega ^B$ in~\Cref{Def: DomUniStabDet}. 
The finite propagation velocity implies, that for any point $\omega \in \Omega$ the distance to the closest point of the control domain must be bounded from above uniformly in $T$ and $|\Omega|$. 
The uniform boundedness of the feedback operators mean, that the individual control intervals $(a_j,b_j)$ can not be to small since the impact of a bounded feedback on the solution's exponential decay depends on the size of the control domain it is acting on. 
The finite propagation velocity is a core property of hyperbolic equations which is independent of the specific choice of the boundary conditions, the domain set $\mathcal{O}$ and the families $(A_\Omega, B_\Omega)_{\Omega \in \mathcal{O}}$. Although we acknowledge, that scaling the control domain with the spatial domain may prove challenging in applications. We therefore have to emphasize, that this condition is indispensable in our setting}.
\end{remark}

 The remaining part of this section is dedicated to proving~\Cref{Thm: Stabilizability}.
To this end we will first derive a solution formula for a transport equation with a state feedback of the form~\eqref{GeneralStateFeedback}, i.e. for the equation
\begin{equation}
\label{TransportEquationwithStateFeedback}
    \frac{\partial}{\partial t}(\omega,t)=-\chi _{\Omega _c^L} (\omega) k_B(\omega)x(\omega,t)-c\frac{\partial}{\partial \omega}x(\omega,t).
\end{equation}
This solution formula is then leveraged to derive an algebraic condition on a piecewise constant state feedback $k_B$, such that the semigroup corresponding with~\eqref{TransportEquationwithStateFeedback} is domain-uniformly exponentially stable.
\begin{lemma}
    \label{Lem: SolutionFormulaStateFeedback1}
    The mild solution of~\eqref{TransportEquationwithStateFeedback} with boundary condition~\eqref{BoundaryCondition} and initial condition~\eqref{InitialCondition} is given by
    \begin{equation*}
        x(\cdot,t) = e^{-\frac{1}{c}\int_{\cdot-ct}^{\cdot}P_{\Omega _L}(k_B|_{\Omega _L^c})(y)\mathrm{d}y} P_{\Omega _L}(x^0)(\cdot-ct).
    \end{equation*}
\end{lemma}
\begin{proof}
    See~\Cref{App: Lemmata}.
\end{proof}

\noindent Due to~\Cref{Lem: SolutionFormulaStateFeedback1} the family $(\mathcal{T}^\varphi_{\Omega_L})_{\Omega _L \in \mathcal{O}}$ of strong\-ly continuous semigroups generated by the operator family $(A_{\Omega _L}+B_{\Omega _L}K^B_{\Omega _L})_{\Omega _L \in \mathcal{O}}$ corresponding to~\eqref{GeneralStateFeedback} is given by
\begin{equation}
\label{SemigroupGeneralDef}
    \mathcal{T}^\varphi_{\Omega_L}(t): L^2(\Omega _L) \rightarrow L^2(\Omega _L), \platz \mathcal{T}^\varphi_{\Omega_L}(t)x^0:= e^{-\frac{1}{c}\int _{\cdot-ct}^\cdot P_{\Omega _L}(k_B|_{\Omega _c^L})(y)\mathrm{d}y}P_{\Omega _L}(x^0)(\cdot-ct).
\end{equation}
In the second step of proving~\Cref{Thm: Stabilizability} we derive some necessary conditions on the control domain. We only consider control domains, which are countable unions of intervals in this section. Note that controls which only act on a nullset $N \subset \R _{\geq 0}$ do not have any influence on the solution of the transport equation since the input operator $B_{\Omega _L}$ is bounded. In particular they are irrelevant for the domain-uniform stabilizability of this equation. Therefore w.l.o.g.\ we only consider countable unions of \textit{closed} intervals.~\Cref{Thm: Stabilizability} can be extended to measurable control domains without major changes in conditions (ii) and (iii).
In view of applications and for the sake of readability of our results, we decided to limit ourselves to countable unions of intervals.

\begin{lemma}
\label{Lem: Necessary}
    Let $\Omega _c = \underset{j \in \N}{\cup} [a_j,b_j] \subset \R _{\geq 0}$ with $0\leq a_1 \leq b_1 \leq a_2 \leq \ldots$. If~\eqref{eq: Transport} with control domain $\Omega _c$ can be domain-uniformly stabilized via a state feedback of the form~\eqref{GeneralStateFeedback}, then the following are true:
    \begin{itemize}
        \item [(i)] $a_1 = 0$
        \item [(ii)] $|\Omega _c| = \infty$
    \end{itemize}
\end{lemma}
\begin{proof}
    See~\Cref{App: Lemmata}.
\end{proof}

\noindent In view of \Cref{Lem: Necessary}, we may thus w.l.o.g.\ assume that the control domains are of the form $\Omega _c:=\underset{j \in \N}{\cup} [a_j,b_j]$ where $\left(a_j\right) _{j\in \N}$ is an unbounded sequence and $0 = a_1 < b_1 \leq a_2 < b_2 \leq \cdots$.  We write $I_j:=[a_j,b_j]$. In the next step of proving~\Cref{Thm: Stabilizability} we will only consider piecewise constant feedbacks
\begin{align}
\label{PiecewiseConstantStateFeedback}
    k_B: \R _{>0} \rightarrow \R , \platz k_B(\omega):= \sum _{j=1}^\infty k_j \chi _{[a_j,b_j]}(\omega),
\end{align}
where $(k_j)_{j\in \N} \in \R _{>0}$ is a bounded sequence with $\hat{k}:=\sup_{j\in \N}\, k_j$. The state feedback~\eqref{PiecewiseConstantStateFeedback} domain-uniformly stabilizes the controlled transport equation~\eqref{TransportEquationControlled} if and only if the corresponding family of closed-loop semigroups is domain-uniformly exponentially stable. Due to~\eqref{SemigroupGeneralDef} this is equivalent to
\begin{equation*}
    \exists M,k>0 \, \forall  L>0 \, \forall  t\geq 0: \n{e^{-\frac{1}{c}\int_{\cdot-ct}^{\cdot}P_{\Omega _L}(k_B|_{\Omega _L})(y)\mathrm{d}y} P_{\Omega _L}(x^0)(\cdot-ct)}_{L^2([0,L])} \leq M e^{-kt} \n{x^0}_{L^2(\Omega _L)}.
\end{equation*}
Since the periodization operator preserves the $L^2(\Omega_L)$-norm (see~\eqref{Eq: NormPreservation}) it suffices to show an estimate of the form
\begin{equation}
    \label{CentralExpEstimate}
    \exists \,M,k>0 \, \forall \, L>0 \, \forall \,\omega \in [0,L] \,\forall \, t\geq 0: e^{-\frac{1}{c}\int_{\omega-ct}^{\omega}P_{\Omega _L}(k_B|_{\Omega _L})(y)\mathrm{d}y} \leq M e^{-kt}
\end{equation}
which only includes the first (exponential) term from the solution formula in~\Cref{Lem: SolutionFormulaStateFeedback1}. In the following~\Cref{Lem: IntervalCondition}  we are able to show such an estimate in three steps by finding algebraic conditions for estimates of the form
\begin{equation}
\label{ExponentialInequality}
    \forall t \geq 0: e^{-\int _\omega^{\omega + t} e(\tau )\mathrm{d}\tau} \leq M e^{-kt},
\end{equation}
where $e: \R_{\geq 0} \rightarrow \R_{\geq 0}$.

\begin{lemma}
\label{Lem: IntervalCondition}
Let $k_B \in L^\infty (\R _{>0},\R_{\geq 0})$ be a state feedback as in~\eqref{PiecewiseConstantStateFeedback}. Then the following statements are true:
\begin{itemize}
    \item [(i)] For a given $k>0$ there exists a constant $M>0$ such that~\eqref{ExponentialInequality} holds for $e=k_B$, $\omega = 0$, if and only if
    \begin{equation}
        \label{IntervalCondition1}
        \exists \, \tilde{M} > 0 \,\forall \, n \in \N: ka_n - \sum _{j=1}^{n-1} k_j (b_j-a_j) = ka_n - \sum _{j=1}^{n-1} k_j B_j \leq \tilde{M}.
    \end{equation}
    \item [(ii)] For a given $k>0$ there exists a constant $M>0$ such that~\eqref{ExponentialInequality} holds for $e=k_B$, $\omega \in \R_{\geq 0}$, if and only if there exists a constant $\tilde{M}>0$ such that for all $n,m \in \N$ with $n\geq m$
    \begin{equation}
    \label{IntervalCondition}
        k(a_n-b_m) - \sum _{j=m+1}^{n-1} k_j (b_j-a_j) = k(a_n-b_m) - \sum _{j=m+1}^{n-1} k_j B_j \leq \tilde{M}.
    \end{equation}
    \item [(iii)]For a given $k>0$ and $e_L:=k_B|_{\Omega _L}$ there exists a constant $M>0$ such that
    \begin{equation}
    \label{ExpInequalityPeriodified}
        \forall L>0 \,\forall  \omega \in [0,L] \, \forall t \geq 0: e^{-\frac{1}{c}\int _{\omega}^{\omega+ct}P_{\Omega _L}(e_L)(y )\mathrm{d}y} \leq Me^{-kt}.
    \end{equation}
    if and only if~\eqref{IntervalCondition} is fulfilled.
\end{itemize}
\end{lemma}
\begin{proof}
    See~\Cref{App: Lemmata}.
\end{proof}

\noindent Using the estimate from~\Cref{Lem: IntervalCondition}(iii) we can now derive a necessary and sufficient condition on the piecewise constant state feedback $k_B$ such that the corresponding family of semigroups is domain-uniformly exponentially stable.
\begin{thm}
\label{Thm: CharExpStab}
    The following two statements are equivalent:
    \begin{itemize}
        \item [(i)] The family of of closed-loop semigroups in~\eqref{SemigroupGeneralDef} with state feedback $k_B \in L^\infty (\R _{>0},\R_{\geq 0})$ as in~\eqref{PiecewiseConstantStateFeedback} is domain-uniformly exponentially stable, i.e.
    \begin{equation*}
        \exists \tilde{M},\tilde{k}>0\, \forall L>0 \, \forall t\geq 0: \n{\mathcal{T}_{\Omega _L}(t)}_{L(L^2([0,L]),L^2([0,L]))} \leq \tilde{M}e^{-\tilde{k}t}.
    \end{equation*}
    \item[(ii)] There exist constants $M,k>0$ such that for all $n,m \in \N$ with $n\geq m$,
    \begin{equation*}
        k(a_n-b_m) - \sum _{j=m+1}^{n-1} k_j (b_j-a_j) = k(a_n-b_m) - \sum _{j=m+1}^{n-1} k_j B_j \leq M.
    \end{equation*}
    \end{itemize}
\end{thm}
\begin{proof}
    \textbf{(i) $\implies$ (ii):} We show this implication by contraposition. Assume that (ii) is not fulfilled. Taking~\Cref{Lem: IntervalCondition} (iii) into account we know, that for all $\tilde{M},\tilde{k}>0$ there exist $L_0>0$ and $\omega_0 \in [0,L],t_0 \geq 0$ such that 
    \begin{equation*}
        e^{-\frac{1}{c}\int _{\omega_0}^{\omega_0+ct_0}P_{L_0}(e_{L_0})(y )\mathrm{d}y} > \tilde{M}e^{-\tilde{k}t_0}.
    \end{equation*}
    Since $P_{L_0}(e_{L_0})$ is a piecewise constant function the expression $e^{-\frac{1}{c}\int _{\omega_0}^{\omega_0+ct_0}P_{L_0}(e_{L_0})(y )\mathrm{d}y}$ is continuous in $\omega_0$. Therefore there exists an interval $[\varepsilon _1, \varepsilon _2] \subset [0,L]$ such that
    \begin{equation*}
        \forall \omega \in [\varepsilon _1, \varepsilon _2]: e^{-\frac{1}{c}\int _{\omega_0}^{\omega_0+ct_0}P_{L_0}(e_{L_0})(y )\mathrm{d}y} > \tilde{M}e^{-\tilde{k}t_0}.
    \end{equation*}
    and $\varepsilon _2>\varepsilon _1$.
    Define $\varepsilon:= \varepsilon _2 - \varepsilon _1 >0$. Using Euclidean division we know, that there exist $k \in \N$ and $t_1 \in \left[0,\frac{L_0}{c}\right]$ such that $t_0=k\frac{L_0}{c}+t_1$. Choose $g_0 \in L^2([0,L_0])$ such that
    \begin{equation*}
         g_0(\omega):=\left\{ \begin{array}{cc}
            \frac{1}{\sqrt{\varepsilon}}  & \omega \in [\max \{0,\varepsilon _1 -ct_1\}, \max \{0,\varepsilon _2 -ct_1\}]\\
             \frac{1}{\sqrt{\varepsilon}}  & \omega \in [\min \{L,\varepsilon _1 -ct_1+L\}, \min \{L,\varepsilon _2 -ct_1+L\}]\\
             0 & \mathrm{else}
         \end{array}\right..
    \end{equation*}
    We find
    \begin{flalign*}
        &\n{(\mathcal{T}_{L_0}(t_0)g_0)}_{L^2([0,L])}^2= \int _0^L \n{e^{-\frac{1}{c}\int _{\omega-ct_0}^\omega P_{\Omega _L}(e_L)(y)\mathrm{d}y}P_{\Omega _L}(x^0)(\omega-ct_0)}^2 \mathrm{d}\omega\\
        &= \int _0^L \n{e^{-\frac{1}{c}\int _{\omega-ct_0}^\omega P_{\Omega _L}(e_L)(y)\mathrm{d}y}P_{\Omega _L}(x^0)(\omega-ct_1)}^2 \mathrm{d}\omega= \int _{\varepsilon_1}^{\varepsilon_2} e^{-2\frac{1}{c}\int _{\omega-ct_0}^\omega P_{\Omega _L}(e_L)(y)\mathrm{d}y}\frac{1}{\varepsilon } \mathrm{d}\omega > \tilde{M}e^{2-\tilde{k}t_0}.
    \end{flalign*}
    This shows
    \begin{equation*}
        \n{\mathcal{T}_{L_0}(t_0)}_{L(L^2([0,L]),L^2([0,L]))}\geq \tilde{M}e^{2-\tilde{k}t_0}
    \end{equation*}
    since $\n{g_0}_{L^2([0,L])}=1$. Therefore (i) is not fulfilled. \\
    \textbf{(ii) $\implies$ (i):} Let $x^0 \in L^2([0,L])$ be arbitrary. Using~\Cref{Lem: IntervalCondition} (iii) we find $\tilde{M},\tilde{k}>0$ such that the equalities
    \begin{flalign*}
        \n{(\mathcal{T}_{\Omega _L}(t)x^0)}_{L^2([0,L])}^2&= \int _0^L \n{e^{-\frac{1}{c}\int _{\omega-ct}^\omega P_{\Omega _L}(e_L)(y)\mathrm{d}y}P_{\Omega _L}(x^0)(\omega-ct)}^2 \mathrm{d}\omega\\
        &= \int _0^L e^{-\frac{2}{c}\int _{\omega-ct}^\omega P_{\Omega _L}(e_L)(y)\mathrm{d}y} \n{P_{\Omega _L}(x^0)(\omega-ct)}^2 \mathrm{d}\omega\\
        &= \tilde{M}e^{-2\tilde{k}t}\int _0^L \n{P_{\Omega _L}(x^0)(\omega-ct)}^2 \mathrm{d}\omega\\
        &= \tilde{M}e^{-2\tilde{k}t}\n{x^0}_{L^2([0,L])}^2
    \end{flalign*}
    hold. This shows
        $\n{\mathcal{T}_{\Omega _L}(t)}_{L(L^2([0,L]),L^2([0,L]))} \leq \tilde{M}e^{-\tilde{k}t}.$
\end{proof}

\noindent Our main result of this section, i.e., \Cref{Thm: Stabilizability}, follows from~\Cref{Thm: CharExpStab}. This is shown as follows:

\begin{proof}[Proof of~\Cref{Thm: Stabilizability}]
    \textbf{(i) $\Lra$ (ii): } Note that (ii) is equivalent to
    \begin{equation*}
        \exists \, M,K,k>0 \, \forall \, n,m \in \N \textrm{ with }n\geq m: k(a_n-b_m) - \sum _{j=m+1}^{n-1} K (b_j-a_j) \leq M. 
    \end{equation*}
    We first show the claim for piecewise constant state feedbacks of the form~\eqref{PiecewiseConstantStateFeedback}. For necessity we assume, that there exists a piecewise constant stabilizing feedback of the form~\eqref{PiecewiseConstantStateFeedback}. Then~\Cref{Thm: CharExpStab} implies that (ii) is fulfilled for $K:=\sup k_j$. If on the other hand (ii) is fulfilled then we can choose the piecewise constant feedback with $k_j := K$ for all $j\in \N$. For this constant feedback the algebraic conditions from~\Cref{Thm: Stabilizability} and~\Cref{Thm: CharExpStab} are equivalent which means that it is domain-uniformly stabilizing. Finally note that if their is a general domain-uniformly stabilizing feedback $k_B$ of the form~\eqref{GeneralStateFeedback} with $\sup k_B(\omega)=K$ then there also exists a piecewise constant domain-uniformly stabilizing feedback. The piecewise constant feedback can be choosen via $\forall \omega \in \R_{>0}: k_B^{pc}(\omega) := K$.\\
    \textbf{(ii) $\Lra$ (iii): }
    Let $a:= \inf I$ and $b:= \sup I$. Again we assume w.l.o.g. that the control domain takes the form $\Omega_c:=\underset{j \in \N}{\cup} I_j$, $I_j:=[a_j,b_j]$, $\left(a_j\right) _{j\in \N}$ is an unbounded sequence and $0 =a_0 < b_1 \leq a_2 < b_2 \leq \cdots$.\\
    We first show (iii) $\implies$ (ii): 
    Note, that by choosing $n=m+1$, \eqref{IntervalConditionTheorem} implies
    \begin{equation}
        \label{Eq: MaxDistance}
        \forall \, m \in \N _0: k(a_{m+1}-b_m)\leq 1 \quad \Leftrightarrow \quad a_{m+1}-b_m\leq \frac{1}{k}.
    \end{equation}
    Define $m_a = \min \{m\in \N _0: b_{m} \geq a\}$ and $m_b = \max \{m\in \N _0: a_{m} \leq b\}$. In the case $a_{m_a}<a$ we have the inequality
    \begin{equation}
    \label{Eq: Case Distinction}
        b_{m_a} - \max \{a_{m_a},a\} = b_{m_a} - a = b_{m_a} - a_{m_a} + a_{m_a} - a \overset{\eqref{Eq: MaxDistance}}{\geq} b_{m_a} - a_{m_a} - \frac{1}{k}.
    \end{equation}
    This implies $b_{m_a} - \max \{a_{m_a},a\}\geq b_{m_a} - a_{m_a} - \frac{1}{k}$ for arbitrary $a_{m_a} \in \R_{\geq 0}$. Analogously we find the estimate $\max \{b_{m_b},b\} - a_{m_b}\geq b_{m_b} - a_{m_b} - \frac{1}{k}$. Overall this leads to
    \begin{flalign*}
        |\Omega _c \cap I| &= b_{m_a}-\max \{a, a_{m_a}\} + \sum _{j=m_a+1}^{m_b -1} (b_j - a_j) + \min \{b, b_{m_b}\} - a_{m_b} \overset{\eqref{Eq: MaxDistance}}{\geq}
        \sum _{j=m_a}^{m_b} (b_j - a_j) - \frac{2}{k}\\
        &\overset{\eqref{IntervalConditionTheorem}}{\geq} \frac{k}{K} (a_{m_b + 1} - b_{m_a}-1) - \frac{1}{K} - \frac{2}{k} \geq \frac{k}{K} |I| - \frac{1}{K} - \frac{2}{k}.
    \end{flalign*}
    This shows (ii) with $c_1 = \frac{k}{K}$ and $c_0 = \frac{1}{K}-\frac{2}{k}$.\\
    To show (ii) $\implies$ (iii) we
    choose $I = [b_m, a_n]$ for arbitrary $n,m \in \N _0$. Thus, (ii) implies
    \begin{equation*}
       |\Omega _c \cap I| = \sum _{j=m+1}^{n-1} (b_j-a_j) \geq c_1 |I| - c_0 =  c_1(a_n-b_m) - c_0,
    \end{equation*}
    i.e., (iii) with $k=\frac{c_1}{c_0}$ an $K=\frac{1}{c_0}$.
\end{proof}

\noindent We now briefly state the resulting exponential sensitivity result for the optimal control problem 
\begin{flalign}
\tag{$\mathrm{OCP}_{L} ^T$}
\label{TransportOCPAbstract}
    \begin{split}
    \underset{(x,u)}{\min }\,\,\, \frac{1}{2} \int _0^T \n{C_{\Omega _L} (x(t)-x^\mathrm{ref}_{\Omega _L})}_{L^2([0,L])}^2 \, + \, & \n{R_{\Omega _L}(u(t)-u^\mathrm{ref} _{\Omega _L})}_{L^2(\Omega_L^c)}^2 \, \mathrm{d}t\\
    \textrm{s.t.}: \platz \dot{x}= A_{\Omega _L} x + B_{\Omega _L} u &= 0, \,\, x(0)=x_{\Omega _L}^0
    \end{split}
\end{flalign}
which is constrained by the transport equation with periodic boundary conditions, i.e. $\mathcal{O} := \{[0,L]:L>0\}$, $A_{\Omega _L}$ and $B_{\Omega _L}$ are given by~\eqref{Eq: TransportOperator} and~\eqref{eq: ControlOperator} and $\Omega_L^c:=[0,L]\cap \Omega _c$ where $\Omega _c \subset \R_{\geq 0}$ is a countable union of closed intervals. By combining~\Cref{Thm: SensitivityResult},~\Cref{Thm: Boundedness} and~\Cref{Thm: Stabilizability} we find the following Corollary.

\begin{corollary}
\label{Cor: TransportDecay}
    Assume that $\Omega _c$ and $\Omega _o$ fulfill one of the conditions (ii) and (iii) in~\Cref{Thm: Stabilizability}.
    Consider a disturbance $\varepsilon \!\in \!(L^1 (\R _{\geq 0};L^2(\R_{\geq 0}))\times L^2(\R_{\geq 0}))^2$ 
    for which the family $\left(\varepsilon _{\Omega_L} ^T\right)_{\Omega_L \in \mathcal{O}, T>0}\! \in \! W_\mathcal{O}^1$ is exponentially localized in
    $F_\mathcal{O} = W_\mathcal{O}^1$ or $F_\mathcal{O} = W_\mathcal{O}^2$ with $\n{e^{\mu \n{P-\cdot}}\varepsilon _\Omega  ^T}_{F_{\Omega _L} ^T}<C _\varepsilon < \infty$. Let $\delta x_{\Omega _L}^T$, $\delta \lambda_{\Omega _L}^T$ and $\delta u_{\Omega _L}^T$ be the solution of the corresponding error system~\eqref{ErrorOptimalitySystem}.
    Then there exists constants $\mu,K>0$ such that for all $T>0$, $L>0$
    \begin{equation*}
        \n{e^{\mu \n{P-\cdot}_1} \delta x_{\Omega _L}^T}_{2\land \infty} + 
        \n{e^{\mu \n{P-\cdot}_1} \delta \lambda_{\Omega _L}^T}_{2\land \infty} + \n{e^{\mu \n{P-\cdot}_1} \delta u_{\Omega _L}^T}_{2\land \infty} \leq K \n{e^{\mu \n{P-\cdot}_1} \varepsilon_{\textcolor{black}{_{\Omega _L}}}^T}_{1\lor 2} \leq K\cdot C_\varepsilon.
    \end{equation*}
\end{corollary}
\begin{proof}
    From~\Cref{Thm: Stabilizability} we know that the family $\left(A_{\Omega _L}, B_{\Omega _L} \right)_{L>0, T>0}$ is domain-uniformly stabilizable. To show domain-uniform detectability of $\left(A_{\Omega _L}, C_{\Omega _L}\right)_{L>0, T>0}$ we first compute $A_{\Omega _L}^*$.  By definition of $A_{\Omega _L}$ and using partial integration we find
\begin{equation*}
    \forall \, x,v\in \mathrm{dom}(A_{\Omega _L}): \left \langle A_{\Omega _L}x, v\right\rangle_{X_{\Omega _L}} = \left \langle -cx', v\right\rangle_{X_{\Omega _L}} = \left \langle x, cv'\right\rangle_{X_{\Omega _L}}.
\end{equation*}
Therefore the operator $-A_{\Omega _L}: \mathrm{dom}(A_{\Omega _L}) \rightarrow L^2([0,L])$ is a formal adjoint operator of $A_{\Omega _L}$. Since $\mathrm{dom}(A_{\Omega _L})$ is a dense subset of $L^2([0,L])$ it can be extended to a maximal formal adjoint operator which is \textit{the} adjoint operator~\cite[Section 2.8]{Weiss2009}. This shows $\mathrm{dom}(A_{\Omega _L}) \subset \mathrm{dom}(A_{\Omega _L}^*)\subset H_1([0,L])$. Let $v \in H_1([0,L])$ such that for all $x \in \mathrm{dom}(A_{\Omega _L})$
\begin{equation*}
    \left \langle A_{\Omega _L}x,v \right\rangle _{X_{\Omega _L}} \!= \left \langle x,A_{\Omega _L}^*v \right\rangle _{X_{\Omega _L}} \!= \left \langle x,-A_{\Omega _L}v \right\rangle _{X_{\Omega _L}} \!\!\!\!\!\implies \!\! x(L)v(L)-x(0)v(0) = x(0)(v(L)-v(0))=0.
\end{equation*}
This implies $v\in \mathrm{dom}(A_{\Omega _L})$. Therefore $\mathrm{dom}(A_{\Omega _L})=\mathrm{dom}(A_{\Omega _L}^*)$, i.e. $A_{\Omega _L}$ is skew-adjoint. Therefore the domain-uniform stabilizability of $\left(A_{\Omega _L}, B_{\Omega _L}\right)_{L>0, T>0}$ is equivalent to the domain-uniform detectability of $\left(A_{\Omega _L}, C_{\Omega _L}\right)_{L>0, T>0}$. The opposite direction of transport does not matter for stabilizability/detectability. Using domain-uniform stabilizability and detectability~\Cref{Thm: Boundedness} implies boundedness of the solution operator $\M^{-1}$ from~\Cref{Sec: ProblemStatement}. Using~\Cref{Thm: SensitivityResult} concludes the proof.
\end{proof}

\section{Characterization of the domain-uniform stabilizability of the transport equation for space-dependent transport velocities}
\label{Sec: TransportNonConstant}
In this section we consider the controlled transport equation~\eqref{eq: Transport} with space-dependent transport velocity. Let $c: \R_{\geq 0} \rightarrow \R _{>0}$ be a piecewise Lipschitz-continuous function which fulfills the condition
\begin{equation}
\label{eq: VelBoundedness}
    \exists \, c_\mathrm{min}, c_\mathrm{max} >0 \, \forall \omega \in \R _{\geq 0}: c_\mathrm{min}\leq c(\omega) \leq c_\mathrm{max}.
\end{equation}
In the following we write $c_L$ for $c|_{\Omega_L}$. In this case domain-uniform stabilizability can be characterized in the same way as for the case of constant transport velocity (see~\Cref{Thm: Stabilizability}). 
\begin{thm}
\label{Thm: StabilizabilityNonConstantCase}
The following two statements are equivalent:
\begin{itemize}
    \item [(i)] The controlled transport equation~\eqref{eq: Transport} can be domain-uniformly stabilized via a state feedback of the form~\eqref{GeneralStateFeedback}.
    \item [(ii)] There exist constants $K>k>0$ such that for all $n,m \in \N$ with $n\geq m$
    \begin{equation}
        k(a_n-b_m) - \sum _{j=m+1}^{n-1} K (b_j-a_j) \leq 1, 
    \end{equation}
    where $\Omega _c:=\underset{j \in \N}{\cup} [a_j,b_j]$, $\left(a_j\right) _{j\in \N}$ is an unbounded sequence with $a_0 = 0$, $a_i < b_i\leq a_{i+1}$ for $i \in \mathbb{N}$.
    \item [(iii)] There exist constants $c_0, c_1>0$ such that for all intervals $I\subset \R _{\geq 0}$ the inequality
    \begin{equation*}
        |\Omega _c \cap I| \geq c_1 |I| - c_0
    \end{equation*}
    is fulfilled and $\forall \, L>0: |\Omega _L^c| = |\Omega _L \cap \Omega _c| >0$.
\end{itemize}
\end{thm}

\noindent The remaining part of this section is devoted to proving~\Cref{Thm: StabilizabilityNonConstantCase}. For this purpose we will again consider feedback operators of the form~\eqref{GeneralStateFeedback}. The corresponding equation with feedback is given by
\begin{equation}
\label{GeneralTransportEquationwithStateFeedback}
    \forall t \in [0,T]\, \forall \omega \in [0,L]: \frac{\partial}{\partial t}{x}(\omega,t)=-k(\omega)x(\omega,t)-c_L(\omega)\frac{\partial}{\partial \omega}x(\omega,t).
\end{equation}
Again, ignoring the feedback term for now, we can rewrite the controlled transport equation as an inhomogeneous Cauchy problem~\eqref{eq: CauchyTransport} with generator
\begin{equation}
\label{Eq: TransportOperatorNonConstant}
   A_{\Omega _L}: D(A_{\Omega _L}):= \{x \in H^1(\Omega _L): x(0)=x(L)\}\subset L^2(\Omega_L)  \rightarrow L^2(\Omega _L), \, A_{\Omega _L}x= -c_L \frac{\partial}{\partial \omega}x
\end{equation}
and input operator~\eqref{eq: ControlOperator}.

In the following, we utilize ordinary differential equations to provide solution formulas for the transport equation, closely related to the method of characteristics. To this end, consider a small particle located at position $p_0 \in [0,L]$ at time $t_0=0$. Its movement can be described by the ordinary differential equation
\begin{equation}
    \label{ODETransportForward}
    \dot{p}(t)=P_{\Omega _L}(c_L)(p(t)), \platz p(0)=p_0,
\end{equation}
where $p(t)$ models the 
position of the particle at time $t$.
Conversely if we observe the particle at a position $\omega \in [0,L]$ at time $t$, then its previous position a time $t_0=0$ zero can be computed via the solution of
\begin{equation}
    \label{ODETransportBackward}
    \dot{q}(t)=-P_{\Omega _L}(c_L)(q(t)), \platz q(0)=q_0.
\end{equation}
We denote the solutions of~\eqref{ODETransportForward} and~\eqref{ODETransportBackward} at time $t\geq 0$ with initial values $p_0, q_0\in [0,L]$ by $p(t,p_0)$ and $q(t,q_0)$, respectively.

In the following auxiliary result, we derive a solution for the transport equation with space-dependent velocity in the uncontrolled and the state feedback case. This is done by replacing the term $\omega - ct$ in~\Cref{Lem: SolConstant} resp.~\Cref{Lem: SolutionFormulaStateFeedback1} by the solution $q(t,\omega)$ of~\eqref{ODETransportBackward}. 
\begin{lemma}[Backward motion]
    \label{Lem: SolNonConstantTransport}
    Consider the controlled transport equation with transport velocity $c: [0,L] \rightarrow \R_{\geq 0}$ and the corresponding differential equations~\eqref{ODETransportForward} and~\eqref{ODETransportBackward}. Then the following hold:
    \begin{itemize}
    \item [(i)] The differential equations~\eqref{ODETransportForward} and~\eqref{ODETransportBackward} each have a unique absolutely continuous (i.e., a.e.\ differentiable) solution. 
    \item [(ii)] Let $p_0 \in [0,L]$, $T>0$ and $q_0:=p(T,p_0)$. Then for all $t\in [0,T]$
    \begin{equation*}
        q(t,q_0)=p(T-t, p_0).
    \end{equation*}
    \item [(iii)] For all $q_0 \in [0,L]$ and all $t \geq 0$ the unique solution $q(t,q_0)$ of~\eqref{ODETransportBackward} fulfills the equality
    \begin{equation*}
        \frac{\partial}{\partial q_0}q(t,q_0)=-\frac{1}{P_{\Omega _L}(c_L)(q_0)} \frac{\partial}{\partial t}{q}(t,q_0) = -\frac{1}{P_{\Omega _L}(c_L)(q_0)} P_{\Omega _L}(c_L)(q(t,q_0)).
    \end{equation*}
    \item [(iv)] The mild solution of~\eqref{eq: Transport} with $u\equiv 0$ is given by
    \begin{equation*}
        \forall t\in [0,T]: x(\cdot,t) = P_{\Omega _L}(x^0)(q(t,\cdot )).
    \end{equation*}
    \item [(v)] The mild solution of~\eqref{GeneralTransportEquationwithStateFeedback} with boundary condition~\eqref{BoundaryCondition} and initial condition~\eqref{InitialCondition} is
    \begin{equation*}
        \forall t\in [0,T]: x(\cdot,t) = e^{-\int_{q(t,\cdot)}^{\cdot}\frac{P_{\Omega _L}(k|_{\Omega_L})(y)}{P_{\Omega _L}(c_L)(y)}\mathrm{d}y} P_{\Omega _L}(x^0)(q(t,\cdot)).
    \end{equation*}
    \end{itemize}
\end{lemma}
\begin{proof}
    See~\Cref{App: LemmataNonConstant}.
\end{proof}

If the control is chosen as a state feedback $u=K^B_{\Omega _L}x$ with a feedback operator as defined in~\eqref{GeneralStateFeedback} then the corresponding operator $A_{\Omega _L}+B_{\Omega _L}K^B_{\Omega _L}$ generates the semigroup
\begin{equation*}
   \mathcal{T}^\varphi_{\Omega_L}(t): L^2([0,L],\R) \rightarrow L^2([0,L],\R), \platz \mathcal{T}^\varphi_{\Omega_L}(t)x^0:=e^{-\int_{q(t,\cdot)}^{\cdot}\frac{P_{\Omega _L}(k|_{\Omega_L})(y)}{P_{\Omega _L}(c_L)(y)}\mathrm{d}y} P_{\Omega _L}(x^0)(q(t,\cdot)).
\end{equation*}
as a direct consequence of~\Cref{Lem: SolNonConstantTransport}. This explicit representation can be used to prove the main result of this part. 

\begin{proof}[Proof of~\Cref{Thm: StabilizabilityNonConstantCase}]
The equivalence (ii) $\Lra$ (iii) was already shown in~\Cref{Thm: Stabilizability}. Therefore we only need to show (i) $\Lra$ (ii) and we start with preliminary derivations. Let $K_\Omega ^B$ be an arbitrary state feedback as in~\eqref{GeneralStateFeedback} with $K:=\n{k_B}_\infty$. For any $L>0$ we find 
    \begin{flalign}
    \label{SemigroupRewritten}
    \begin{split}
        \forall t \geq 0: \n{\mathcal{T}^\varphi_{\Omega_L}(t)x^0}^2_{L^2([0,L],\R)}&=\n{e^{-\int_{q(t,\cdot)}^{\cdot}\frac{P_{\Omega _L}(k_B)(y)}{P_{\Omega _L}(c_L)(y)}P_{\Omega _L}(\chi _{\Omega _L^c})(y)\mathrm{d}y} P_{\Omega _L}(x^0)(q(t,\cdot))}_{L^2([0,L],\R)}^2\\
        &= \int _0^L e^{-2\int_{q(t,\omega)}^{\omega}\frac{P_{\Omega _L}(k_B)(y)}{P_{\Omega _L}(c_L)(y)}P_{\Omega _L}(\chi _{\Omega _L^c})(y)\mathrm{d}y}\n{P_{\Omega _L}(x^0)(q(t,\omega))}^2 \mathrm{d}\omega.
    \end{split}
    \end{flalign}
For each $t\geq 0$ we define a transformation $\mathcal{T}_t: [0,L] \rightarrow [q(t,0), q(t,L)], \, \mathcal{T}_t(\omega ):= q (t,\omega)$
where $q(\cdot, \omega)$ is the solution of the initial value problem~\eqref{ODETransportBackward} with initial value $q _0 = \omega$.
\Cref{Lem: SolNonConstantTransport}(ii) shows
$\mathcal{T}_t(p(t,\omega))=q(t,p(t,\omega))=\omega$ 
such that the inverse transformation is given by
$\mathcal{T}_t^{-1}: [q(t,0), q(t,L)] \rightarrow [0,L], \, \mathcal{T}_t^{-1}(\omega )= p (t,\omega)$.~\Cref{Lem: SolNonConstantTransport}(iii) implies
\begin{equation*}
    \frac{\mathrm{d}}{\mathrm{d}\omega}\mathcal{T}_t(\omega) = -\frac{\frac{\partial}{\partial t}{q}(t,\omega)}{c_L(\omega)}=\frac{c_L(q(t,\omega))}{c_L(\omega)}.
\end{equation*}
Therefore, we have
\begin{flalign*}
    \int _0^L \n{P_{\Omega _L}(x^0)(z)}^2 \mathrm{d}z &= \int _{q(t,0)}^{q(t,L)} \n{P_{\Omega _L}(x^0)(q(t,\omega))}^2\frac{c_L(q(t,\omega))}{c_L(\omega)} \mathrm{d}\omega\\
    &= \int _0^L \n{P_{\Omega _L}(x^0)(q(t,\omega))}^2\frac{c_L(q(t,\omega))}{c_L(\omega)} \mathrm{d}\omega
\end{flalign*}
which implies
\begin{flalign}
\label{EstimatePLG}
\begin{split}
    &\frac{c_\mathrm{min}}{c_\mathrm{max}}\n{x^0}_{L^2([0,L])}^2 \leq \int _0^L \n{P_{\Omega _L}(x^0)(q(t,\omega))}^2\\
    =&\int _0^L \n{P_{\Omega _L}(x^0)(z)}^2 \frac{c_L(p(t,z))}{c_L(z)} \mathrm{d}z \leq \frac{c_\mathrm{max}}{c_\mathrm{min}}\n{x^0}_{L^2([0,L])}^2
\end{split}
\end{flalign}
since
\begin{equation*}
    \int _0^L \n{P_{\Omega _L}(x^0)(z)}^2 \mathrm{d}z= \int _0^L \n{x^0(z)}^2 \mathrm{d}z = \n{x^0}_{L^2([0,L])}^2.
\end{equation*}
Furthermore, we have 
\begin{equation}
    \label{ExpTermEstimate}
    e^{-2\int_{q(t,\omega)}^{\omega}\frac{P_{\Omega _L}(k_B)(y)}{P_{\Omega _L}(c_L)(y)}P_{\Omega _L}(\chi _{\Omega _L^c})(y)\mathrm{d}y}\geq e^{-2\frac{K}{c_\mathrm{min}}\int_{\omega - c_\mathrm{max}t}^{\omega}P_{\Omega _L}(\chi _{\Omega _L^c})(y)\mathrm{d}y}.
\end{equation}

\noindent \textbf{(i) $\Rightarrow$ (ii):} We prove this claim by contraposition. Assume that (ii) is not fulfilled. Let $K_\Omega ^B$ be an arbitrary feedback operator as in~\eqref{GeneralStateFeedback} with $\n{k_B}_\infty =: K>0$. Define $\forall j \in \N: k_j=K$.
Since (ii) is not fulfilled we find that condition~\eqref{IntervalCondition} from~\Cref{Lem: IntervalCondition}(ii) is also not fulfilled which implies for arbitrary $k>0$ that
\begin{equation}
\label{ExpTermEstimate2}
    \forall M>0 \,\exists L_M>0 \,\exists \,  \omega _M \in [0,L_M] \, \exists \, t_M \geq 0: e^{-\frac{1}{c_\mathrm{max}}\int _{\omega-c_\mathrm{max}t_M}^{\omega}K \, \chi _{\Omega_{L_M} ^c}(y)\mathrm{d}y} \geq Me^{-kt}.
\end{equation}
Inserting~\eqref{EstimatePLG},~\eqref{ExpTermEstimate} and~\eqref{ExpTermEstimate2}into~\eqref{SemigroupRewritten} (and leaving the index $M$ out for readability) yields
\begin{flalign}
    \label{SemigroupEstimate}
    \begin{split}
    \n{\mathcal{T}^\varphi_{\Omega_L}(t)x^0}^2_{L^2([0,L],\R)} &\geq \int _0^L e^{-2\frac{K}{c_\mathrm{min}}\int_{\omega - c_\mathrm{max}t}^{\omega}P_{\Omega _L}(\chi _{\Omega _L^c})(y)\mathrm{d}y}\n{P_{\Omega _L}(x^0)(q(t,\omega))}^2 \mathrm{d}\omega \\
    &\geq \int _0^L \left(M e^{-kt}\right)^{2\frac{c_\mathrm{min}}{c_\mathrm{max}}}\n{P_{\Omega _L}(x^0)(q(t,\omega))}^2 \mathrm{d}\omega\\
    &\geq \frac{c_\mathrm{min}}{c_\mathrm{max}}\left(M e^{-kt}\right)^{2\frac{c_\mathrm{min}}{c_\mathrm{max}}}\n{x^0}_{L^2(\Omega_L)}^2
    \end{split}
\end{flalign}
and therefore (i) is not fulfilled since
\begin{equation*}
    \n{\mathcal{T}^\varphi_{\Omega_L}(t)}\geq \sqrt{\frac{c_\mathrm{min}}{c_\mathrm{max}}}\left(M e^{-kt}\right)^{\frac{c_\mathrm{min}}{c_\mathrm{max}}}.
\end{equation*}
\textbf{(ii) $\Rightarrow$ (i):} Let $M,k,K$ be such that (ii) is fulfilled. Define $K_\Omega ^B$ as in~\eqref{GeneralStateFeedback} with $k_B\equiv K$. Then we have
\begin{flalign*}
    \n{\mathcal{T}^\varphi_{\Omega_L}(t)x^0}^2_{L^2([0,L],\R)}&\overset{\eqref{SemigroupRewritten}}{=}\int _0^L e^{-2\int_{q(t,\omega)}^{\omega}\frac{K}{P_{\Omega _L}(c_L)(y)}P_{\Omega _L}(\chi _{\Omega _L^c})(y)\mathrm{d}y}\n{P_{\Omega _L}(x^0)(q(t,\omega))}^2 \mathrm{d}\omega\\
    &\overset{}{\leq}\int _0^L e^{-2\int_{\omega-c_\mathrm{min}t}^{\omega}\frac{K}{c_\mathrm{max}}P_{\Omega _L}(\chi _{\Omega _L^c})(y)\mathrm{d}y}\n{P_{\Omega _L}(x^0)(q(t,\omega))}^2 \mathrm{d}\omega\\
    &\leq \int _0^L (Me^{-2kt})^{\frac{c_\mathrm{max}}{c_\mathrm{min}}}\n{P_{\Omega _L}(x^0)(q(t,\omega))}^2 \mathrm{d}\omega\\
    &\leq \frac{c_\mathrm{max}}{c_\mathrm{min}}M^{2\frac{c_\mathrm{max}}{c_\mathrm{mn}}} e^{-2\frac{c_\mathrm{max}}{c_\mathrm{min}}kt}\n{x^0}_{L^2([0,L])}^2.
\end{flalign*}
This shows domain-uniform exponential stability of the semigroup and therefore (i).
\end{proof}

In the case of space-dependent transport velocity the differential operator $A_{\Omega _L}$ is no longer skew-adjoint. Therefore domain-uniform detectability has to be discussed separately from domain-uniform stabilizability in this section. For this purpose we first compute the adjoint operator in the following Lemma.
\begin{lemma}[Adjoint of $A_{\Omega _L}$ in the case of space-dependent transport velocity]
    Let $c \in H^1(\R_{\geq 0})$. Then for all $L>0$ the adjoint of the operator $A_{\Omega _L}$ defined in~\eqref{Eq: TransportOperatorNonConstant} is given by
    \begin{equation*}
        A_{\Omega _L}^*: D(A_{\Omega _L}^*)\subset L^2(\Omega _L) \rightarrow L^2(\Omega _L), \platz (A_{\Omega _L}^*x)(\omega):= c'_L(\omega) x(\omega)+c_L(\omega) x'(\omega),
    \end{equation*}
    where $D(A_{\Omega _L}^*):=\{v \in H^1(\Omega _L): c(0)v(0)=c(L)v(L)\}$.
\end{lemma}
\begin{proof}
    By definition of $A_{\Omega _L}$ and using integration by parts we find
    \begin{equation*}
        \forall \, x\in \mathrm{dom}(A_{\Omega _L}), v\in H^1(\Omega_L): \left \langle A_{\Omega _L}x, v\right\rangle_{X_{\Omega _L}} = \left \langle -c_Lx', v\right\rangle_{X_{\Omega _L}} = \left \langle x, c_L'v+c_Lv'\right\rangle_{X_{\Omega _L}} - \left[c_Lxv\right]_0^L.
    \end{equation*}
    Since $x(0)=x(L)$ the right-hand boundary value term only vanishes if $c(0)v(0)=c(L)v(L)$. An analogous argumentation as in the proof of~\Cref{Cor: TransportDecay} yields the result.
\end{proof}
\noindent By definition, domain-uniform detectability of $(A_{\Omega _L},C_{\Omega _L})$ is equivalent to domain-uniform stabilizability of $(A_{\Omega _L}^*,C_{\Omega _L}^*)$. It is easy to see, that the adjoint output operator $C_{\Omega _L}^*$ corresponds to the input operator $B_{\Omega _L}$ defined in~\eqref{eq: ControlOperator}. Therefore domain-uniform detectability of $(A_{\Omega _L},C_{\Omega _L})$ leads to the domain-uniform stabilizability of the continuity equation
\begin{subequations}
\label{eq: Continuity}
\begin{align}
\label{ContinuityEquationControlled}
    \forall (\omega,t) \in \Omega_L \times [0,T]: \frac{\partial}{\partial t}{x}(\omega,t)&=-\frac{\partial}{\partial \omega}\left(c_L(\omega)x(\omega,t)\right) + \chi _{\Omega_L^o}(\omega) u(\omega,t)\\
    \label{BoundaryConditionContinuity}
    \forall t \in [0,T]: c(0)x(0,t) &= c(L)x(L,t) \\
    \label{InitialConditionContinuity}
    \forall \omega \in \Omega_L: x(\omega,0) &= x^0(\omega) = x_{\Omega _L}^0.
\end{align}
\end{subequations}
If the control is chosen as a feedback of the form~\eqref{GeneralStateFeedback} then an analogous argumentation as in~\Cref{Lem: SolNonConstantTransport} yields that the solution of~\eqref{eq: Continuity} is given by
\begin{equation}
\label{eq: ContSol}
    x(\omega, t)=\left( \frac{c(0)}{c(L)}\right)^{N_L(\omega,t)} e^{\int_{\omega}^{p(t,\omega)}\frac{P_{\Omega _L}\left(k_L^C\right)(y)}{P_{\Omega _L}(c_L)(y)}\mathrm{d}y} P_{\Omega _L}(x^0)(p(t,\omega)),
\end{equation}
where for all $L>0$
\begin{equation*}
    k_L^C: \Omega _L \rightarrow \R, \quad k_L^C(\omega):=c_L'(\omega ) + \chi _{\Omega _L^o}(\omega) k_B(\omega )
\end{equation*}
and
\begin{equation*}
    N_L: [0,L]\times \R_{\geq 0} \rightarrow \N _0, \quad N_L(\omega,t):=\#\{s \in [\omega, p(t,\omega) )|\, \exists \, m_s\in \Z: s=m_s L\} =\left \lfloor \frac{p(t,\omega)}{L} \right \rfloor.
\end{equation*}
The corresponding semigroup is given by
\begin{equation}
\label{eq: ContSemigroup}
    T_{\Omega _L}^\mathrm{Cont}(t): L^2(\Omega _L) \rightarrow L^2(\Omega _L), \, \left(T_{\Omega _L}^\mathrm{Cont}(t)x^0\right) = \left( \frac{c(0)}{c(L)}\right)^{N_L(\omega,t)} e^{\int_{\omega}^{p(t,\omega)}\frac{P_{\Omega _L}\left(k_L^C\right)(y)}{P_{\Omega _L}(c_L)(y)}\mathrm{d}y} P_{\Omega _L}(x^0)(p(t,\omega)).
\end{equation}
\begin{thm}[Domain-uniform stabilizability of the continuity equation]
\label{Thm: StabilizabilityCont}
    Let $c \in H^1(\R_{\geq 0})$ such that~\eqref{eq: VelBoundedness} is fulfilled and $c' \in L^2(\R_{\geq 0}) \cap L^\infty(\R_{\geq 0})$. Then for arbitrary $\gamma >0$ the controlled continuity equation equation~\eqref{eq: Continuity} can be domain-uniformly stabilized with regard to the domain set $\mathcal{O}_\gamma := \{[0,L]: L\geq \gamma\}$ via a state feedback of the form~\eqref{GeneralStateFeedback} if one of the conditions (ii) and (iii) in~\Cref{Thm: StabilizabilityNonConstantCase} is fulfilled by $\Omega _o$.
\end{thm}
\begin{proof}
    In~\Cref{Thm: Stabilizability} equivalence of (ii) and (iii) was already shown. Therefore it suffices to consider condition (ii).
    Define $\alpha _L:= \frac{c(0)}{c(L)}$. We find the estimates
    \begin{equation*}
        \forall L \geq \gamma: \left( \alpha _L\right)^{N_L(\omega,t)} \leq 1
    \end{equation*}
    for $\alpha_L \leq 1$ and
    \begin{equation*}
        \forall L\geq \gamma: \left( \alpha_L\right)^{N_L(\omega,t)} = e^{\log \left(\alpha_L\right)\left \lfloor \frac{p(t,\omega)}{L} \right \rfloor} \leq e^{\log \left(\frac{c(0)}{c(L)}\right) \frac{p(t,\omega)}{L}} \leq e^{\log \left(\frac{c_\mathrm{max}}{c_\mathrm{min}}\right) \frac{\omega + c_\mathrm{max}t}{\gamma}}
    \end{equation*}
    for $\alpha > 1$. Therefore, there exist $M_\alpha, k_\alpha > 0$ such that
    \begin{equation*}
        \forall L\geq \gamma: \left( \alpha_L\right)^{N_L(\omega,t)} \leq M_\alpha e^{k_\alpha t}.
    \end{equation*}
    Assume, that condition (ii) is satisfied. In this case there exist constants $M_e,k_e,K_e>0$ such that
    \begin{equation*}
        \forall \, L \geq \gamma \, \forall \, \omega \in [0,L] \, \forall \, t\geq 0: e^{-\int_{\omega}^{\omega+c_\mathrm{min}t}K_eP_{\Omega _L}(\chi _{\Omega _L^o})(y)\mathrm{d}y} \leq M_e e^{-k_e t}.
    \end{equation*}
    For $k>0$ define a domain-uniformly stabilizing state feedback via
    \begin{equation*}
        k_B: \R_{\geq 0} \rightarrow \R_{<0},\quad k_B(\omega):= -K, \quad K:= c_\mathrm{max}\frac{k+k_\alpha + \frac{c_\mathrm{max}}{c_\mathrm{min}}\n{c'}_\infty}{k_e}K_e.
    \end{equation*}
    Overall we find the estimate
    \begin{flalign*}
        &\n{\mathcal{T}^\mathrm{Cont}_{\Omega_L}(t)x^0}^2_{L^2([0,L],\R)}\overset{\eqref{eq: ContSemigroup}}{=}\int _0^L \left( \frac{c(0)}{c(L)}\right)^{2N_L(\omega,t)} e^{2\int_{\omega}^{p(t,\omega)}\frac{P_{\Omega _L}\left(k_L^C\right)(y)}{P_{\Omega _L}(c_L)(y)}\mathrm{d}y} \n{P_{\Omega _L}(x^0)(p(t,\omega))}^2 \mathrm{d}\omega\\
        \overset{}{\leq}&M_\alpha^2 e^{2k_\alpha t}\int _0^L e^{2\int_{\omega}^{p(t,\omega)}\frac{P_{\Omega _L}\left(c_L'\right)(y)}{P_{\Omega _L}(c_L)(y)}\mathrm{d}y} e^{-2\int_{\omega}^{p(t,\omega)}\frac{KP_{\Omega_L}\left(\chi _{\Omega_L^o}\right)(y)}{P_{\Omega _L}(c_L)(y)}\mathrm{d}y}\n{P_{\Omega _L}(x^0)(p(t,\omega))}^2 \mathrm{d}\omega\\
        \overset{}{\leq}&M_\alpha^2 e^{2k_\alpha t}\int _0^L e^{2\frac{c_\mathrm{max}}{c_\mathrm{min}}\n{c'}_\infty t} e^{-2\int_{\omega}^{\omega + c_\mathrm{min}t}\frac{KP_{\Omega_L}\left(\chi _{\Omega_L^o}\right)(y)}{c_\mathrm{max}}\mathrm{d}y}\n{P_{\Omega _L}(x^0)(p(t,\omega))}^2 \mathrm{d}\omega\\
        \leq & M_\alpha^2 M_e^{2\frac{K}{K_e}}e^{-2kt}\int _0^L\n{P_{\Omega _L}(x^0)(p(t,\omega))}^2 \mathrm{d}\omega\\
        \leq &\frac{c_\mathrm{max}}{c_\mathrm{min}}M_\alpha^2 M_e^{2\frac{K}{K_e}}e^{-2kt}\n{x^0}_{L^2([0,L])}^2
\end{flalign*}
for the corresponding closed-loop semigroup where the last inequality can be obtained analogously to the proof of~\Cref{Thm: StabilizabilityNonConstantCase}.
This shows domain-uniform exponential stability of the semigroup.
\end{proof}

\noindent Again consider the optimal control problem ~\eqref{TransportOCPAbstract} however with differential operator $A_{\Omega _L}$ as defined in~\eqref{Eq: TransportOperatorNonConstant}. By combining~\Cref{Thm: SensitivityResult},~\Cref{Thm: Boundedness},~\Cref{Thm: Stabilizability} and~\Cref{Thm: StabilizabilityCont} we find a Corollary similar to~\Cref{Cor: TransportDecay}.
\begin{corollary}
\label{Cor: TransportDecayNonConstant}
    Assume that $\Omega _c$ and $\Omega_o$ fulfill one of the conditions (ii) and (iii) in~\Cref{Thm: Stabilizability}.
    Let $c \in H^1(\R_{\geq 0})$ such that~\eqref{eq: VelBoundedness} is fulfilled and $c' \in L^2(\R_{\geq 0}) \cap L^\infty(\R_{\geq 0})$. Let $\gamma >0$ be an arbitrary constant and $\mathcal{O}:=\{[0,L]: L\geq \gamma\}$ be the correponding set of spatial domains.
    Consider a disturbance $\varepsilon \!\in \!(L^1 (\R _{\geq 0};L^2(\R_{\geq 0}))\times L^2(\R_{\geq 0}))^2$ 
    for which the family $\left(\varepsilon _{\Omega_L} ^T\right)_{\Omega_L \in \mathcal{O}, T>0}\! \in \! W_\mathcal{O}^1$ is exponentially localized in
    $F_\mathcal{O} = W_\mathcal{O}^1$ or $F_\mathcal{O} = W_\mathcal{O}^2$ with $\n{e^{\mu \n{P-\cdot}}\varepsilon _\Omega  ^T}_{F_{\Omega _L} ^T}<C _\varepsilon < \infty$. Let $\delta x_{\Omega _L}^T$, $\delta \lambda_{\Omega _L}^T$ and $\delta u_{\Omega _L}^T$ be the solution of the corresponding error system~\eqref{ErrorOptimalitySystem}.
    Then there exist $\mu,K>0$ such that
    \begin{flalign*}
        \forall \, T>0 \, L\geq \gamma: \n{e^{\mu \n{P-\cdot}_1} \delta x_{\Omega _L}^T}_{2\land \infty} + 
        \n{e^{\mu \n{P-\cdot}_1} \delta \lambda_{\Omega _L}^T}_{2\land \infty} + \n{e^{\mu \n{P-\cdot}_1} \delta u_{\Omega _L}^T}_{2\land \infty} 
        &\leq K\cdot C_\varepsilon.
    \end{flalign*}
\end{corollary}
\begin{proof}
From~\Cref{Thm: Stabilizability} we know that the family $\left(A_{\Omega _L}, B_{\Omega _L} \right)_{L>0, T>0}$ is domain-uniformly stabilizable. From~\Cref{Thm: StabilizabilityCont} we know that the family $\left(A_{\Omega _L}^*, C_{\Omega _L}^* \right)_{L\geq\gamma, T>0}$ is domain-uniformly stabilizable.
Therefore the family $\left(A_{\Omega _L}, C_{\Omega _L}\right)_{L\geq \gamma, T>0}$ is domain-uniformly detectable. Using domain-uniform stabilizability and detectability~\Cref{Thm: Boundedness} implies boundedness of the solution operator $\M^{-1}$ from~\Cref{Sec: ProblemStatement}. Using~\Cref{Thm: SensitivityResult} concludes the proof.
\end{proof}

\section{Domain-uniform stabilizability of the wave equation on a one-dimensional domain}
\label{Sec: WaveEq}
In the previous sections~\ref{Sec: TransportConstant} and~\ref{Sec: TransportNonConstant}, we derived necessary and sufficient conditions for domain-uniform stabilizability and detectability of the transport equation with distributed control and periodic boundary condition. In this section, we extend these results to the wave equation. This equation plays a key role in the modeling of vibrations and oscillations with a large variety of applications in electrical engineering, mechanics, hydraulics and even acoustics/music theory. In order to enable the further use of our previously presented results on domain-uniform stabilizability, we will consider a suitable non-local state feedback law, which leads to a coupled system of two damped transport equations.

\subsection{System formulation and unitary group}

\noindent As before, we consider the family of domains $\mathcal{O}:=\{\Omega _L:=[0,L]: L>0\}$ where the corresponding control domains are given by $\Omega _L^c := \Omega_c \cap \Omega _L$ and $\Omega_c \subset \R_{\geq 0}$ denotes some global control domain of positive Lebesgue-measures. On these domains the wave equation with distributed control and Dirichlet boundary conditions is given by
\begin{subequations}
\label{eq: Wave}
\begin{align}
\label{WaveEquationControlled}
    \forall \, (\omega,t) \in \Omega_L \times [0,T]: \frac{\partial^2}{\partial t^2}{x}(\omega,t)&=c^2\frac{\partial ^2}{\partial \omega ^2}x(\omega,t) + \chi _{\Omega_L^c}(\omega) u(\omega,t)\\
    \label{WaveBoundaryCondition}
    \forall \, t \in [0,T]: x(0,t) &= x(L,t)=0\\
    \label{WaveInitialCondition}
    \forall \, \omega \in \Omega_L: x(\omega,0) &= x_{\Omega _L}^0 \textrm{ and } \frac{\partial}{\partial t}x(\omega,0) = x_{\Omega _L}^1
\end{align}
\end{subequations}
with time horizon $T>0$, control $u \in L^2(\Omega _L^c\times [0,T])=:U_{\Omega _L}$, initial distributions  $x_{\Omega _L}^0\in H^2(\Omega_L) \cap H_0^1(\Omega_L)$ and $x_{\Omega _L}^1\in H^1(\Omega_L)$ and velocity $c>0$. 
We first consider the uncontrolled equation ($u\equiv 0$). In this case, the wave equation can be rewritten as a first-order PDE of the form $\dot{z}=A _{\Omega _L}z$ where $z = \begin{pmatrix}
    x & \frac{\partial}{\partial t}x
\end{pmatrix}^\top$. For this purpose, we define the Dirichlet Laplacian via
\begin{equation*}
    A_{\Omega_L}^0: D(A_{\Omega_L}^0):= \left\{v \in H_0^1(\Omega _L): \frac{\partial^2}{\partial \omega ^2}v \in L^2(\Omega _L)\right\} \rightarrow L^2(\Omega _L), \qquad A_{\Omega_L}^0x:=-c^2\frac{\partial^2}{\partial \omega ^2}x.
\end{equation*}
From~\cite[Proposition 3.6.1]{Weiss2009} we know that $A_{\Omega_L}^0$ is strictly positive, i.e. $A_{\Omega_L}^0=(A_{\Omega_L}^0)^*$ and
\begin{equation*}
    \exists m >0 \, \forall x \in D(A_{\Omega_L}^0): \scp{x}{A_{\Omega_L}^0x}\geq m\n{x}_{L^2(\Omega _L)}^2.
\end{equation*}
Furthermore, this proposition states for the square-root $(A_{\Omega_L}^0)^{\frac{1}{2}}$ of $A_{\Omega_L}^0$, i.e., the unique, strictly positive operator which fulfills $\left((A_{\Omega_L}^0)^{\frac{1}{2}}\right)^2=A_{\Omega_L}^0$, that $D((A_{\Omega_L}^0)^{\frac{1}{2}}=H_0^1(\Omega _L)$ and $\forall x \in D((A_{\Omega_L}^0)^{\frac{1}{2}}): \n{x}_{\frac{1}{2}}^2:=\scp{(A_{\Omega_L}^0)^\frac{1}{2}x}{(A_{\Omega_L}^0)^\frac{1}{2}x}=\n{x}_{H^1}^2$. We now define the block operator $A_{\Omega_L}: D(A_{\Omega_L}):=D(A_{\Omega_L}^0)\times D(A_{\Omega_L}^0)^{\frac{1}{2}} = H^2(\Omega _L) \cap H_0^1(\Omega _L) \times H_0^1(\Omega _L)\rightarrow H_0^1(\Omega _L)\times L^2(\Omega _L)$ via
\begin{equation*}
    A_{\Omega_L}z:= \begin{pmatrix}
        0 & I\\
        -A_{\Omega_L}^0 & 0
    \end{pmatrix} \begin{pmatrix}
        z_1\\z_2
    \end{pmatrix} = \begin{pmatrix}
        z_2\\c^2 \frac{\partial ^2}{\partial \omega ^2}z_1
    \end{pmatrix}.
\end{equation*}
\cite[Proposition 3.7.6]{Weiss2009} tells us, that $A$ is skew-adjoint. Therefore, we can apply the Theorem of Stone~\cite[Prop. 3.8.6]{Weiss2009} to show, that $A$ generates a unitary group $\left(T_L(t)\right)_{t\in \R}$ on $H_0^1(\Omega _L)\times L^2(\Omega _L)$. It is easy to see that $x \in D(A_{\Omega_L}^0)$ is a classical solution of the uncontrolled wave equation if and only if $z = \begin{pmatrix}
    x & \frac{\partial}{\partial t}x
\end{pmatrix}^\top \in D(A)$ is a classical solution of
\begin{equation}
\label{WaveEquation0}
    \dot{z} = A_{\Omega_L}z, \qquad z(0)=\begin{pmatrix}
        x_{\Omega _L}^0\\ x_{\Omega _L}^1
    \end{pmatrix}.
\end{equation}
Since in this case $z(t)=T_L(t)z_0$ (see~\cite[Proposition 6.2]{Engel2000}) and $\left(T_L(t)\right)_{t\in \R}$ is a unitary group we find the equality (see also \cite[Proposition 3.8.7]{Weiss2009})
\begin{flalign*}
    \forall t \geq 0: \n{z(t)}^2&=\n{z_1(t)}_{H_0^1(\Omega _L)}^2+\n{z_2(t)}_{L^2(\Omega _L)}^2 = \n{\frac{\partial}{\partial \omega}x(t)}_{L^2(\Omega _L)}^2+\n{\frac{\partial}{\partial t}x(t)}_{L^2(\Omega _L)}^2 \\
    &= \n{z(0)}^2 = \n{\frac{\partial}{\partial \omega}x_{\Omega _L}^0}_{L^2(\Omega _L)}^2+\n{x_{\Omega _L}^1}_{L^2(\Omega _L)}^2.
\end{flalign*}
Using the formula of d'Alembert we find the solution formula
\begin{flalign*}
    x(\omega ,t) = \frac{1}{2}(\tilde{x}_{\Omega _L}^0(\omega + ct) + \tilde{x}_{\Omega _L}^0(\omega -ct))+\frac{1}{2c} \int _{\omega -ct}^{\omega +ct} \tilde{x}_{\Omega _L}^1(s)\mathrm{d}s
\end{flalign*}
for the uncontrolled wave equation where $\tilde{x}_{\Omega _L}^0$ and $\tilde{x}_{\Omega _L}^1$ are the unique odd and $2L$-periodic extensions of $x_{\Omega _L}^0$ and $x_{\Omega _L}^1$, i.e. for $i \in \{0,1\}$ we have
\begin{itemize}
    \item [(i)] $\forall k \in \Z\, \forall \omega \in (2kL,(2k+1)L): \tilde{x}_{\Omega _L}^i(\omega)=x_{\Omega _L}^i(\omega-2kL)$
    \item [(ii)] $\forall k \in \Z\, \forall \omega \in ((2k+1)L,(2k+2)L): \tilde{x}_{\Omega _L}^i(\omega)=-x_{\Omega _L}^i((2k+2)L-\omega))$.
\end{itemize}
All of the above equations can be extended to mild solutions of~\eqref{eq: Wave} using standard density arguments.

\subsection{Domain-uniformly stabilizing control}
In this section, we use a state feedback of the form
\begin{equation}
\label{WaveFeedback}
    u(\omega, t)=-2k\frac{\partial}{\partial t}x(\omega,t)-k^2x(\omega,t)
\end{equation}
to stabilize the system~\eqref{eq: Wave} domain-uniformly. In order to prove that this approach fulfills its purpose we will transform the wave equation into a pair of damped transport equations. After that we can use~\Cref{Thm: Stabilizability} to find suitable control domains.

Inserting the state feedback~\eqref{WaveFeedback} into~\eqref{eq: Wave} leads to the closed-loop system
\begin{flalign}
\label{eq: ClosedLoop}
\begin{split}
    \forall \, (\omega,t) \in \Omega_L \times [0,T]: \frac{\partial^2}{\partial t^2}{x}(\omega,t)&=c^2\frac{\partial ^2}{\partial \omega ^2}x(\omega,t) - \chi _{\Omega_L^c}(\omega) \left(2k \frac{\partial}{\partial t}x(\omega,t) +k^2x(\omega,t)\right)\\
    \forall \, t \in [0,T]: x(0,t) &= x(L,t)=0\\
    \forall \, \omega \in \Omega_L: x(\omega,0) &= x_{\Omega _L}^0 \textrm{ and } \frac{\partial}{\partial t}x(\omega,0) = x_{\Omega _L}^1.
\end{split}
\end{flalign}
This corresponds to the closed-loop Cauchy problem
\begin{equation}
\label{WaveClCauchy}
    \dot{z}(t) = (A_{\Omega _L}+B_{\Omega _L}K^B_{\Omega _L})z(t), z(0)=\begin{pmatrix}
        x_{\Omega _L}^0\\
        x_{\Omega _L}^1
    \end{pmatrix},
\end{equation}
where $A_{\Omega _L}$ is defined as in~\eqref{WaveEquation0}, the input operator is specified by
\begin{equation*}
    B_{\Omega _L}: L^2(\Omega _L^c) \rightarrow H_0^1(\Omega _L)\times L^2(\Omega _L), \qquad B_{\Omega _L}u := \begin{pmatrix}
        0\\ \chi _{\Omega _L^c} u
    \end{pmatrix}
\end{equation*}
and we have the bounded feedback operator
\begin{equation*}
    K^B_{\Omega _L}: H_0^1(\Omega _L)\times L^2(\Omega _L) \rightarrow L^2(\Omega _L^c), \qquad K^B_{\Omega _L} \begin{pmatrix}
        x\\y
    \end{pmatrix} = -2kx-k^2 y.
\end{equation*}
Note that boundedness of $K^B_{\Omega _L}$ follows from the Poincaré inequality. We already know that $A_{\Omega _L}$ is the generator of a unitary group on $H_0^1(\Omega _L)\times L^2(\Omega _L)$. Therefore, we can apply the bounded perturbation theorem for semigroups~\cite[Theorem 1.3]{Engel2000} to conclude that $A_{\Omega _L}+B_{\Omega _L}K^B_{\Omega _L}: D(A_{\Omega _L})\rightarrow H_0^1(\Omega _L)\times L^2(\Omega _L)$ generates a strongly continuous semigroup $(T_L^\mathrm{cl})_{t\geq 0}$, where $\n{T_L^\mathrm{cl}(t)}\leq e^{\n{B_{\Omega _L}K^B_{\Omega _L}}t}$.
The following~\Cref{Prop: Equivalence} gives an equivalent system representation of the closed-loop~\eqref{eq: ClosedLoop} which simplifies the analysis.

\begin{proposition}[Equivalent Representation of the damped wave equation]
\label{Prop: Equivalence}
    $x\in C^2([0,T],H^2(\Omega _L)\cap H_0^1(\Omega _L))$ is a classical solution of~\eqref{eq: ClosedLoop} if and only if there exists $v\in C^1([0,T],H^1(\Omega _L))$ such that $(\xi _1, \xi _2) = (x,v)$ is a classical solution of
    \begin{flalign}
\label{WaveEquation1}
\begin{split}
    \forall \, (\omega,t) \in \Omega_L \times [0,T]:\frac{\partial}{\partial t} \xi(\omega,t) &= \frac{\partial}{\partial t} \begin{pmatrix}
        \xi _1(\omega,t) \\ \xi _2(\omega,t)
    \end{pmatrix} = -\underset{=:D}{\underbrace{\begin{pmatrix}
        0 & 1\\
        c^2 & 0
    \end{pmatrix}}}
    \frac{\partial}{\partial \omega} \xi(\omega ,t) - \chi _{\Omega_L^c}(\omega)\underset{=:F}{\underbrace{\begin{pmatrix}
        k & 0\\ 0 &k
    \end{pmatrix}}}
    \xi(\omega,t)\\
    \forall \, t \in [0,T]: \xi_1(0,t) &= \xi _1(L,t)=0\\
    \forall \, \omega \in \Omega_L: \xi _1(\omega,0) &= x_{\Omega _L}^0 \textrm{ and } \frac{\partial}{\partial t}\xi _1(\omega,0) = x_{\Omega _L}^1.
\end{split}
\end{flalign}
\end{proposition}
\begin{proof}
    See~\Cref{App: WaveEq}.
\end{proof}

\noindent \Cref{Prop: Equivalence} yields a system representation, which is very similar to the strain variable formulation in~\Cref{ex:vs}. The important difference can be found in the differential equation~\eqref{WaveODE} which yields a damping term in both equations of~\eqref{WaveEquation1}. Using strain variables the transformed closed-loop system would only exhibit such a term in the second equation. This property allows us to show the following theorem.

\begin{thm}[Domain-uniform stabilizability of the wave equation]
\label{Thm: StabilizabilityWave}
The following three statements are equiva\-lent:
\begin{itemize}
    \item [(i)] The controlled wave equation~\eqref{eq: Wave} with constant transport velocity $c>0$ can be domain-uniformly stabilized via a state feedback of the form~\eqref{WaveFeedback}.
    \item [(ii)] There exist constants $K>k>0$ such that for all $n,m \in \N$ with $n\geq m$
    \begin{equation}
        k(a_n-b_m) - \sum _{j=m+1}^{n-1} K (b_j-a_j) \leq 1, 
    \end{equation}
    where $\Omega _c:=\underset{j \in \N}{\cup} [a_j,b_j]$, $\left(a_j\right) _{j\in \N}$ is an unbounded sequence with $a_0 = 0$, $a_i < b_i\leq a_{i+1}$ for $i \in \mathbb{N}$.
    \item [(iii)] There exist constants $c_0, c_1>0$ such that for all intervals $I\subset \R _{\geq 0}$ the inequality
    \begin{equation*}
        |\Omega _c \cap I| \geq c_1 |I| - c_2
    \end{equation*}
    is fulfilled and $\forall \, L>0: |\Omega _L^c| = |\Omega _L \cap \Omega _c| >0$.
\end{itemize}
\end{thm}
\begin{proof}
    \noindent Using the simple transformation
\begin{flalign*}
    \begin{pmatrix}
        \xi_1\\
        \xi_2
    \end{pmatrix} = \begin{pmatrix}
        \frac{1}{c}\left(\zeta_1 + \zeta_2\right)\\
        \zeta_1 - \zeta_2
    \end{pmatrix} = \begin{pmatrix}
        \frac{1}{c} & \frac{1}{c}\\
        1 & -1
    \end{pmatrix}\begin{pmatrix}
        \zeta_1\\\zeta_2
    \end{pmatrix} =:T\zeta \implies 
    \begin{pmatrix}
        \zeta_1\\
        \zeta_2
    \end{pmatrix} 
    = \begin{pmatrix}
        \frac{c}{2} & \frac{1}{2}\\
        \frac{c}{2} & -\frac{1}{2}
    \end{pmatrix}\begin{pmatrix}
        \xi_1\\\xi_2
    \end{pmatrix} = T^{-1} \xi ,
\end{flalign*}
we can convert~\eqref{WaveEquation1} into a system of two coupled damped transport equations which is given by
\begin{flalign}
\label{WaveEquation2}
\begin{split}
    \forall \, (\omega,t) \in \Omega_L \times [0,T]:\frac{\partial}{\partial t} \zeta(\omega,t) &= \underset{=T^{-1}DT}{\underbrace{\begin{pmatrix}
        -c & 0\\
        0 & c
    \end{pmatrix}}}
    \frac{\partial}{\partial \omega} \zeta(\omega ,t) - \underset{=T^{-1}FT}{\underbrace{\begin{pmatrix}
        k & 0\\
        0 & k
    \end{pmatrix}}}\chi _{\Omega_L^c}(\omega) \zeta(\omega,t)\\
    \forall \, t \in [0,T]: \zeta _1(0,t) &= \zeta _2(0,t) \quad \textrm{and}\quad \zeta _1(L,t) = \zeta _2(L,t)\\
    \forall \, \omega \in \Omega_L: \zeta _1(\omega,0) &= \frac{1}{2}\left(cx_{\Omega _L}^0(\omega)-\int _0^\omega x_{\Omega _L}^1(s)+k\,\chi _{\Omega_L^c}(s) x_{\Omega _L}^0(s)\mathrm{d}s\right)\\
    \textrm{ and } \zeta _2(\omega,0) &= \frac{1}{2}\left(cx_{\Omega _L}^0(\omega)+\int _0^\omega x_{\Omega _L}^1(s)+k\,\chi _{\Omega_L^c}(s) x_{\Omega _L}^0(s)\mathrm{d}s\right).
\end{split}
\end{flalign}
Let $\zeta \in H^1(\Omega _L)^2$ be a classical solution of~\eqref{WaveEquation2}. Define $v\in H^1(\Omega _{2L})$ via
\begin{equation*}
    v(\omega, t):=\left\{\begin{array}{cc}
        \zeta _1 (\omega, t),& \omega \in [0,L]\\
        \zeta _2 (2L-\omega, t),& \omega \in [L,2L].
    \end{array}\right.
\end{equation*}
Then, $v$ solves the damped transport equation with periodic boundary condition
\begin{flalign*}
\begin{split}
    \forall \, (\omega,t) \in \Omega_{2L} \times [0,T]:\frac{\partial}{\partial t} \zeta(\omega,t) &= c
    \frac{\partial}{\partial \omega} v(\omega ,t) -k\chi _{\Omega_L^c\cup (\Omega_L^c+L)}(\omega) v(\omega,t)\\
    \forall \, t \in [0,T]: v(0,t) &= v(2L,t)
\end{split}
\end{flalign*}
with initial condition
\begin{equation*}
\zeta _1(\omega,0) = \left\{\begin{array}{cc}
        \frac{1}{2}\left(cx_{\Omega _L}^0(\omega)-\int _0^\omega x_{\Omega _L}^1(s)+k\,\chi _{\Omega_L^c}(s) x_{\Omega _L}^0(s)\mathrm{d}s\right),& \omega \in [0,L]\\
        \frac{1}{2}\left(cx_{\Omega _L}^0(2L-\omega)+\int _0^{2L-\omega} x_{\Omega _L}^1(s)+k\,\chi _{\Omega_L^c}(s) x_{\Omega _L}^0(s)\mathrm{d}s\right),& \omega \in [L,2L]
    \end{array}\right.
\end{equation*}
for $\omega \in \Omega _{2L}$.
To this equation we can directly apply~\Cref{Thm: Stabilizability} to characterize domain-uniform stabilizability. Note that $T$ is a (domain-uniformly) bounded transformation, i.e., the damped transport equation is domain-uniformly exponentially stable if and only if the same holds for~\eqref{WaveEquation1} which contains the solution of~\eqref{eq: ClosedLoop} in the first state. This shows the claim.
\end{proof}

\section{Numerical examples}
\label{Sec: Numerics}
To visualize our findings from Sections~\ref{Sec: MainResults} and~\ref{Sec: TransportConstant} the results of several numerical simulations are discussed in this section. For this purpose  we solve the optimal control problem
\begin{flalign}
\label{TransportOCP}
\begin{split}
    \underset{(x,u)}{\min }\,\,\, \frac{1}{2} \int _0^T \n{x(\cdot,t)}_{L^2([0,L])}^2 \, + \, & \alpha ^2 \n{u(\cdot, t)}_{L^2(\Omega _c^L)}^2 \, \mathrm{d}t, \quad \alpha > 0\\
    \textrm{s.t.}: \platz \forall (\omega,t) \in [0,L] \times [0,T]: \frac{\partial}{\partial t}{x}(\omega,t)&=-c(\omega)\frac{\partial}{\partial \omega}x(\omega,t) + \chi _{\Omega_c}(\omega) u(\omega,t)\\
    \forall t \in [0,T]: x(0,t) &= x(L,t) \\
    \forall \omega \in [0,L]: x(\omega,0) &= x^0(\omega) = x_{\Omega _L}^0.
\end{split}
\end{flalign}
The transport velocity is assumed to be constant ($c=2$) in all simulations.
The problem is solved for two different types of control domains:
\begin{itemize}
    \item[(1)] The control domain consists of a single interval which is given by $\Omega_c := [a,b] := [0,0.2]$.
    \item[(2)] The control domain consists of a series of equidistantly distributed intervals. These intervals are parametrized via $\Omega_c := \cup_{j=1}^{\infty} \left[a_j,b_j\right]$ where $a_j := a + (j-1)L_0$, $b_j := b + (j-1)L_0$ and $L_0 := 1$.
\end{itemize}

\noindent To numerically solve~\eqref{TransportOCP} we discretize the optimality system~\eqref{DisturbedOptimalitySystem} via a finite difference method with symmetric difference quotient $(D_hx)(\omega):=\frac{x(\omega + h) - x(\omega -h)}{2h}$ in space and an implicit midpoint rule for discretization in time. The space discretization leads to a skew-symmetric matrix while the implicit midpoint rule is a symplectic integrator. Therefore this discretization method preserves unitarity of the semigroup on a discrete level and thus avoids a corruption of our results by numerical dissipation. 

\begin{figure}[htb]
\centering
\includegraphics[width=0.5\linewidth]{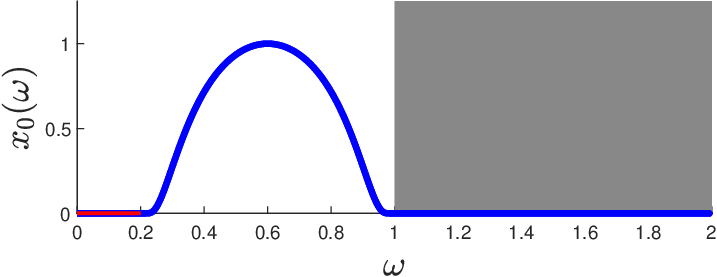}
  \caption{Initial condition function (blue) of \eqref{TransportOCP} with single control interval (red) and distance $| (\Omega \setminus (\Omega _c \cup \mathrm{supp}(\varepsilon))|$ between perturbation and control domain (grey)
  }
  \label{InitialValue}
\end{figure}

In the following we will consider the solution corresponding to a non-zero initial value as a perturbation of the zeros solution of \eqref{TransportOCP} for vanishing initial value. Thus, setting $\varepsilon _\mathrm{width} = 0.8$ and $\varepsilon _\mathrm{center} = 0.6$, we set the initial value (see also \Cref{InitialValue})
\begin{equation}
\label{Perturbation}
    \varepsilon (\omega):= \left\{ \begin{array}{cc}
    e^{\mu(\omega)}, & \omega 
    \!\in \!\left (
    \varepsilon_\mathrm{center}\!-\!\frac{\varepsilon_\mathrm{width}}{2}, \varepsilon_\mathrm{center}\!+\!\frac{\varepsilon_\mathrm{width}}{2}
    \right )
    \\
    0, & \mathrm{else}
    \end{array}\right.\!\mathrm{with}\, \mu(\omega) = 1+\frac{1}{\left(\frac{2}{\varepsilon_\mathrm{width}}\left(\omega\!-\!\varepsilon _\mathrm{center}\right)\right)^2\!-\!1}.
\end{equation}

Here, the choice of the center of perturbation $\varepsilon _\mathrm{center}$ corresponds to a worst case scenario in the sense that the distance $| (\Omega \setminus (\Omega _c \cup \mathrm{supp}(\varepsilon))|$ (highlighted in grey in~\Cref{InitialValue}) over which the perturbation has to be transported until it reaches any part of the control domain is maximal. Therefore the part of the spatial domain influenced by the undamped perturbation is as large as possible. 

\begin{figure}[htb]
\centering
\includegraphics[trim={1.5cm 0.6cm 1.5cm 0},clip,width=0.975\linewidth]{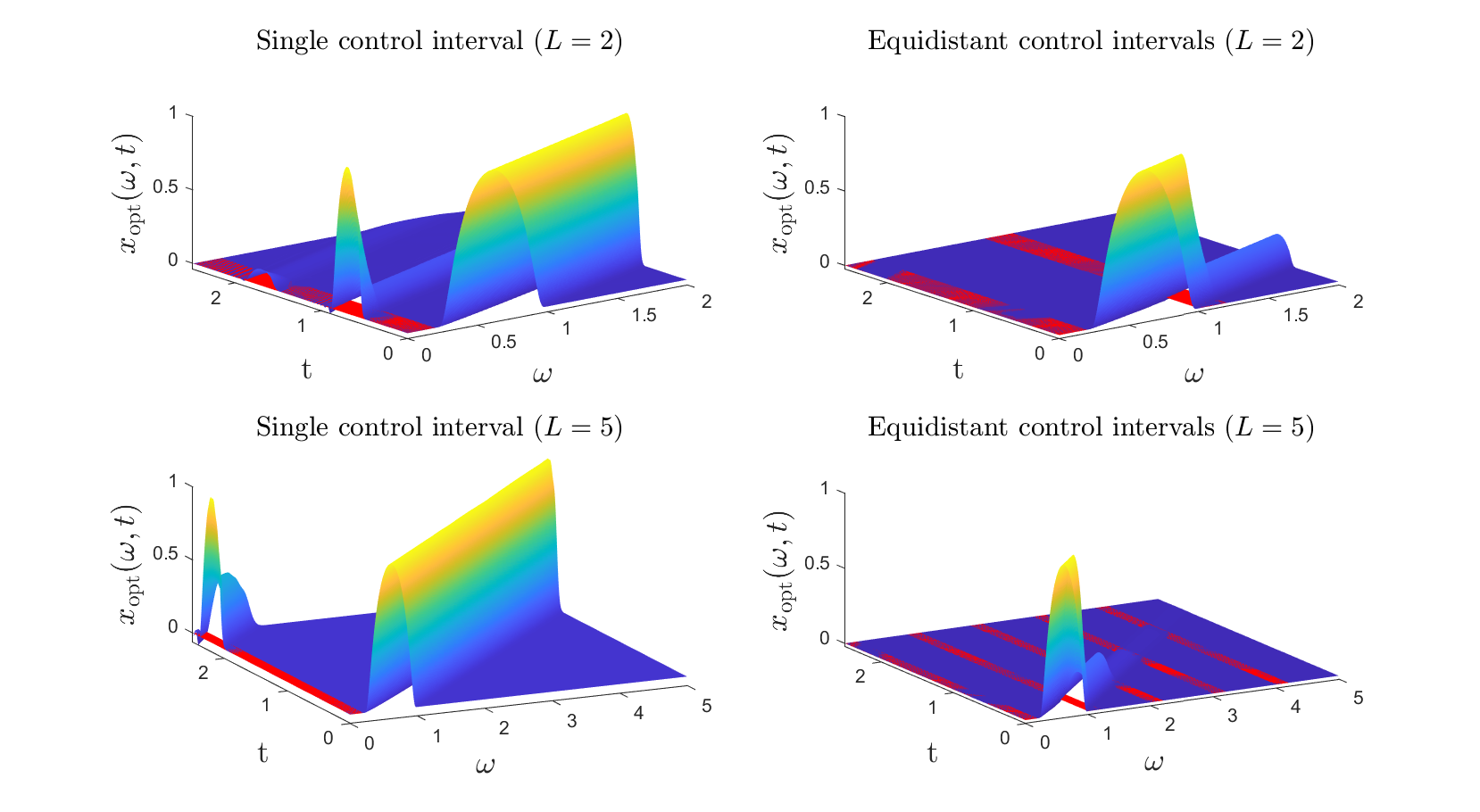}
  \caption{Optimal state (blue, space variable $\omega$) for two different control domains (red) and domain sizes ($T = 2.5$, $\alpha = 0.125$)}
  \label{3DPlotSpace}
\end{figure}

\Cref{3DPlotSpace} shows the state trajectory of an optimal controlled transport equation for two different domain sizes. Over time the perturbation is transported along the spatial domain until it is damped at the control domain. In case of a single control interval the spatial area on which the perturbation's influence on the solution is undiminished increases with the domain size (see the left two plots in~\Cref{3DPlotSpace}). If the control acts on a series of equidistant intervals then the speed of decay of the perturbed solution does not depend on the domain size anymore (see the right two plots in~\Cref{3DPlotSpace}).

\begin{figure}[htb]
\centering
\includegraphics[width=\linewidth]{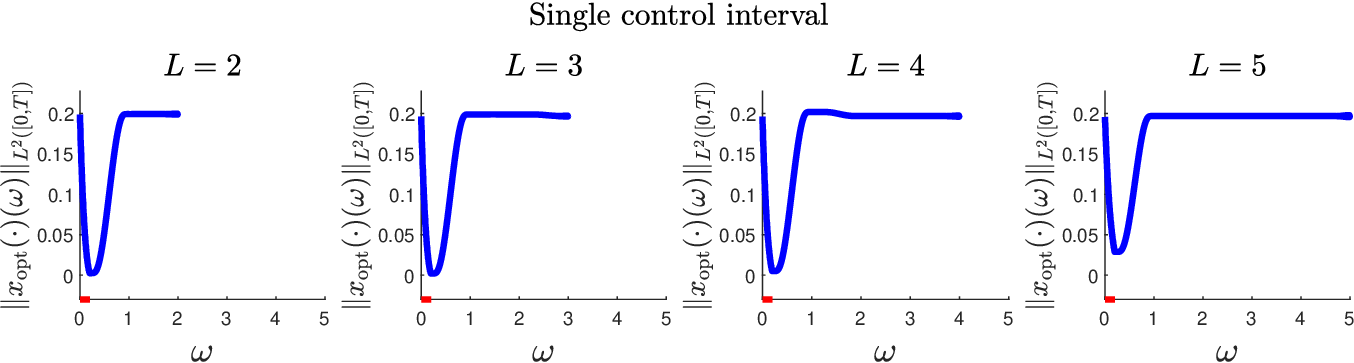}\\
\vspace{.2cm}
\includegraphics[trim={0 0 0 0},clip,width=\linewidth]{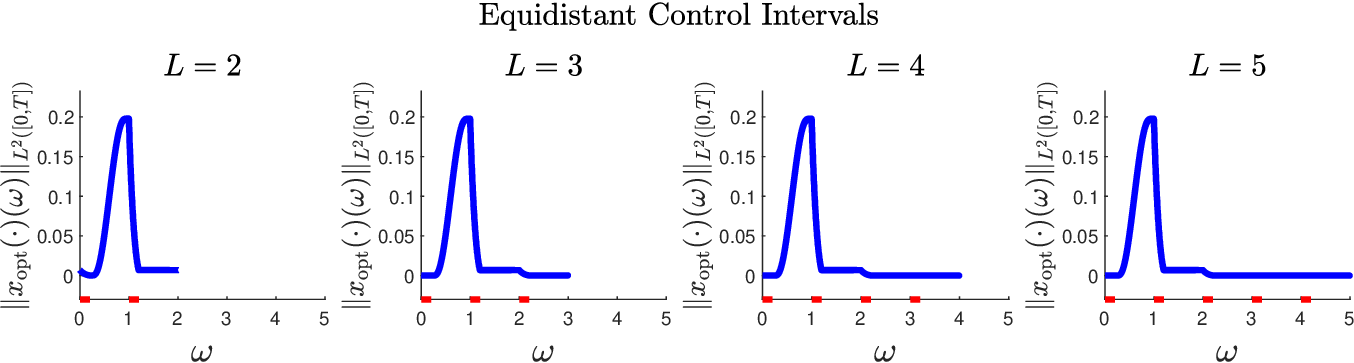}
  \caption{Spatial decay of the optimal state on different domains for a single control interval (top) and equidistant control intervals (bottom)}
  \label{L2NormTime2}
\end{figure}

\noindent\Cref{L2NormTime2} visualizes the \textit{spatial} decay of the optimal state trajectory for $T = 2.5$ and $\alpha = 0.125$. For this purpose and for fixed spatial coordinate $\omega \in [0,L]$, we computed the $L^2([0,T])$-norm of the optimal state $x_\mathrm{opt}(\cdot)(\omega)$ in the time domain. 
As this norm nearly remains constant over the remaining part of the spatial domain (upper row of~\Cref{L2NormTime2}) and as the width of this constant part increases with the domain size, this implies, that the spatial decay of the optimal state is slower for larger domain sizes. In case of equidistant intervals the decay does not change as the domain size increases (lower row of~\Cref{L2NormTime2}).

\begin{figure}[htb]
\centering
\includegraphics[width=0.9\linewidth]{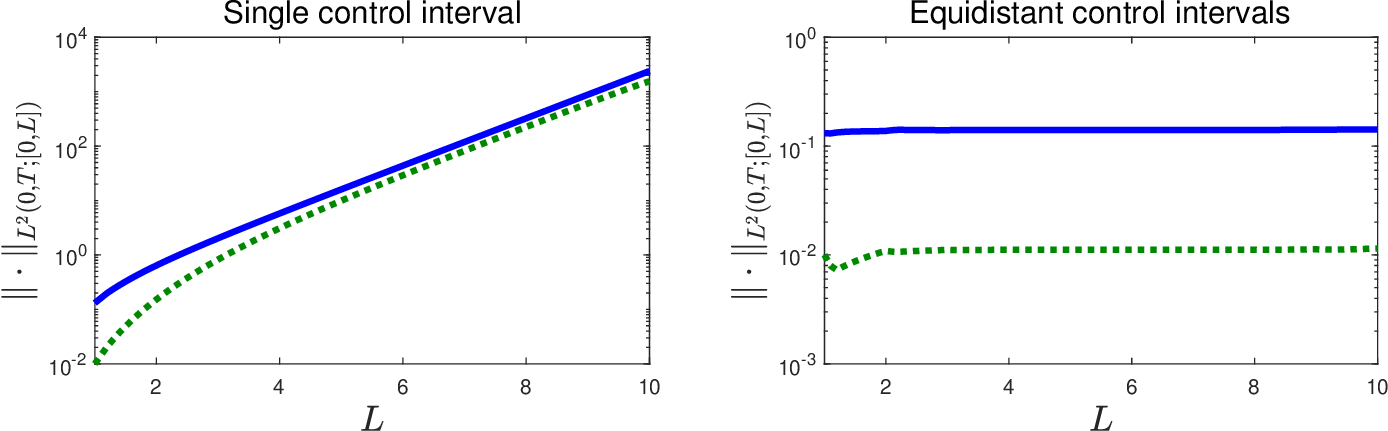}
  \caption{Relation between the domain size and the $L^2(0,T;[0,L])$-norm of the optimal state (blue) respectively costate (green, dotted) for parameters $T = 5$, $\alpha = 0.125$}
  \label{SpaceTimeNormPlot}
\end{figure}

\noindent In~\Cref{SpaceTimeNormPlot} the relation between the $L^2(0,T;[0,L])$-norm of the scaled optimal state and corresponding adjoint state and the domain size is shown. 
In the scenario of a single control interval the transport equation with periodic boundary conditions is not domain-uniformly stabilizable/detectable (see~\Cref{Thm: Stabilizability}). Therefore the (sufficient) assumptions of~\Cref{Thm: SensitivityResult} are not fulfilled, which means that the bound in equation~\eqref{SpatialDecayTheorem} does not necessarily hold.
This is confirmed by the plot on the left of \Cref{SpaceTimeNormPlot} the $L^2$-norm of the state/costate increases exponentially in the domain size $L$.
Evaluating condition (ii) from~\Cref{Thm: Stabilizability} for our proposed sequence of equidistant control intervals we find that there exist $M,k,K >0$ such that
\begin{equation*}
        k(a_n-b_m) - K\sum _{j=m+1}^{n-1} B_j = k(n-1-m)-0.2K(n-m-1) = (k-0.2K)(n-1-m) \leq M
\end{equation*}
for all $n\geq m \in \N$. This condition is for example fulfilled for $k=1$, $K=5$ and arbitrary $M>0$.
Therefore~\Cref{Thm: Stabilizability} ensures the domain-uniform stabilizability and detectability of the transport equation such that the assumptions of~\Cref{Thm: SensitivityResult} are fulfilled. Since the perturbation defined in~\eqref{Perturbation} is exponentially localized around $z=\varepsilon _\mathrm{center}$~\Cref{Thm: SensitivityResult} yields the same for the optimal state and costate. Therefore the norms in the right-hand side plot of~\Cref{SpaceTimeNormPlot} uniformly bounded.
\begin{figure}[H]
\centering
\includegraphics[trim={2cm 0 2cm 0},clip,width=0.95\linewidth]{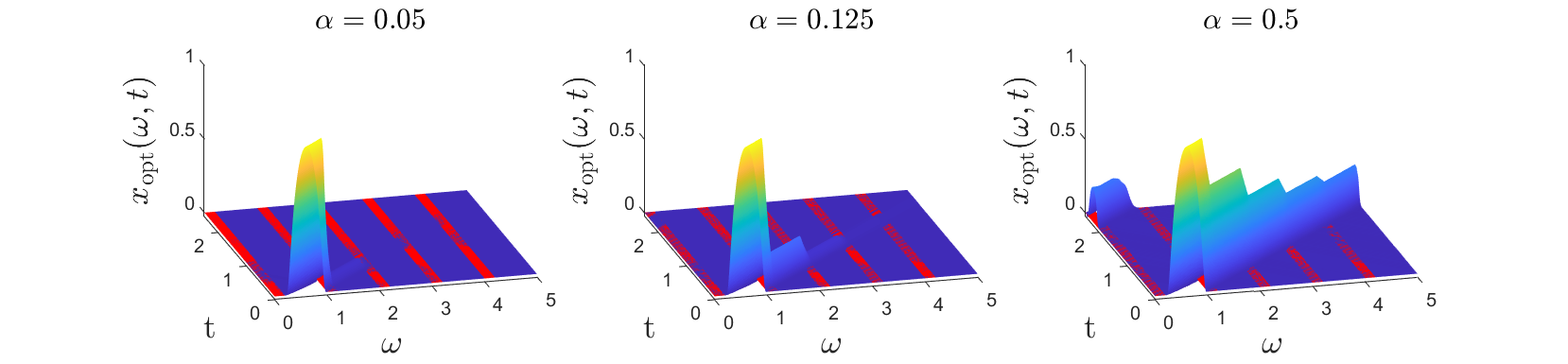}\\
\includegraphics[trim={2cm 0 2cm 0},clip,width=0.95\linewidth]{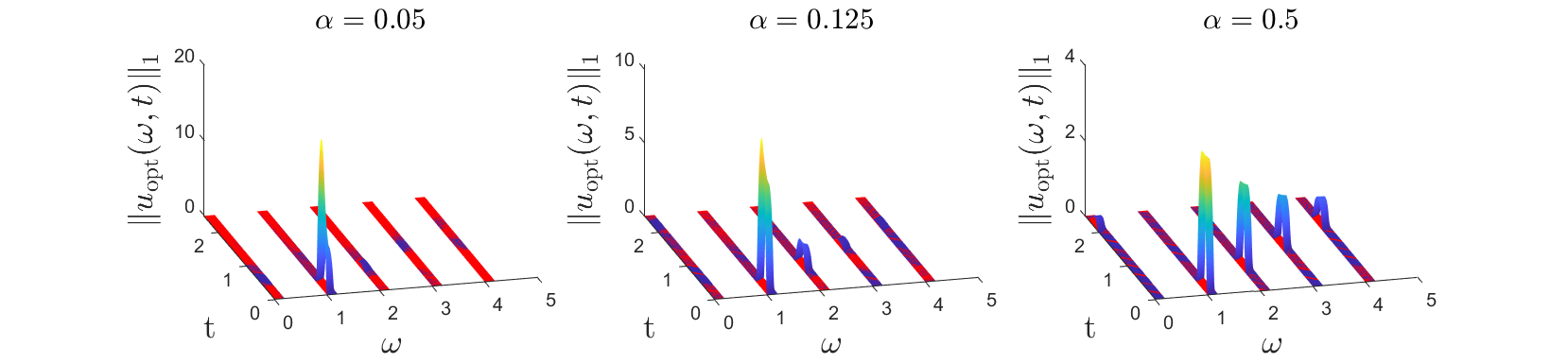}
  \caption{Optimal state (top) and control (bottom) for different control weights with equidistant control intervals (red) for $T = 5$}
  \label{3DPlotCOntrolWeight}
\end{figure}

\noindent Last, in~\Cref{3DPlotCOntrolWeight}, the influence of the control weight $\alpha$ on the solution of the OCP is visualized. A stronger weighting of the control means, that large amplitudes of the control signal come with a disproportionately high cost. Therefore increasing the parameter $\alpha$ leads to a slower decay of the state trajectory and a lower peak in the control signal, once the perturbation reaches the control domain.

\section{Conclusion and outlook}
In this work we considered linear-quadratic optimal control problems governed by general evolution equations. We showed that, under domain-uniform stabilizability and detectability assumptions on the involved operators, that the influence of spatially localized perturbations on the OCP's solution decays exponentially in space. We further characterized the domain-uniform stabilizability/detectability assumption for linear transport equations and provided numerical examples confirming the findings.

Future work considers the extension towards networks of partial differential equations.

\section*{Acknowledgements}
We thank Professor Martin Gugat (FAU Erlangen-Nürnberg) for fruitful discussions with regard to this work. Moreover, we are grateful to Professor Herbert Egger (JKU Linz) for valuable comments regarding the reformulation of the wave equation.

\printbibliography
\appendix
\section{Proof of~\Cref{Thm: Boundedness}}
\label{App: Boundedness}

During the remaining part of this section let $z \in C(0,T;X)^2$ and $r \in \left( L_1 (0,T;X) \times X \right)^2$ with
\begin{equation*}
    z = 
    \begin{pmatrix}
    x\\ \lambda    
    \end{pmatrix}
    = \M^{-1} r = \begin{pmatrix}
        C^*C & -\frac{\mathrm{d}}{\mathrm{d}t} - A ^*\\
        0 & E_T\\
        \frac{\mathrm{d}}{\mathrm{d}t} - A & -BQ^{-1}B^*\\
        E_0 & 0
    \end{pmatrix}^{-1}
    \begin{pmatrix}
        l_1 \\
        \lambda _T\\
        l_2 \\
        x_0
    \end{pmatrix}
\end{equation*}
or equivalently
\begin{equation}
\label{Pontryagin}
     \M z=\begin{pmatrix}
        C^*C & -\frac{\mathrm{d}}{\mathrm{d}t} - A ^*\\
        0 & E_T\\
        \frac{\mathrm{d}}{\mathrm{d}t} - A & -BQ^{-1}B^*\\
        E_0 & 0
    \end{pmatrix}
    \begin{pmatrix}
        x\\ \lambda 
    \end{pmatrix} = r = \begin{pmatrix}
        l_1 \\
        \lambda _T\\
        l_2 \\
        x_0
    \end{pmatrix}
\end{equation}
Since we want to derive a bound $c>0$ such that
\begin{equation*}
    \forall\, T>0 \, \forall \, \Omega \in \mathcal{O}: \n{\M^{-1}}_{L(W^{1\lor 2},W^{2\land \infty})}\leq c.
\end{equation*}
the aim of this section is to find an estimate of the form (leaving out the indices $\Omega$ and $T$)
\begin{equation}
    \n{z}_{2 \land \infty}^2 = \n{x}_{2\land \infty}^2 + \n{\lambda}_{2\land \infty}^2 \leq c \left( \n{l_1}_{1\lor 2}^2 + \n{\lambda _T}_X^2 + \n{l_2}_{1\lor 2}^2 + \n{x_0}_X^2 \right) = c\, \n{r}_{1\lor 2}^2.
\end{equation}
We emphasize that the constant $c>0$ must be independent of $T$ and $\Omega \in \mathcal{O}$.
Since a certain amount of technical obstacles has to be overcome in order to reach this goal, we have structured the corresponding proof in five smaller steps:
First for any $t\in [0,T]$ we will analyze the systems
\begin{equation}
\label{PhiSystem}
    \mathcolor{black}{\forall s \in [0,t]: -\dot{\varphi}(s) = (A_\Omega ^C + K_\Omega ^C C_\Omega)^*\varphi (s)} , \platz \varphi (t) = x_\Omega ^T(t)
\end{equation}
and
\begin{equation}
\label{PsiSystem}
    \mathcolor{black}{\forall s \in [t,T]:\dot{\psi}(s) = (A_\Omega  + B_\Omega K_\Omega ^B)\psi(s)} , \platz \psi (t) = \lambda_\Omega ^T(t)
\end{equation}
corresponding to the stabilizability and detectability conditions in~\Cref{Thm: Boundedness}. We will find exponential bounds for the solutions of these systems. Second we will rewrite the quantities $\n{x_\Omega ^T(t)}_{X_\Omega}^2$ and $\n{\lambda_\Omega ^T (t)}_{X_\Omega}^2$ such that they depend on the solutions $\varphi$ and $\psi$ of~\eqref{PhiSystem} and~\eqref{PsiSystem}. In the third step we will use the bounds on $\varphi$ and $\psi$ which we found in the first step to estimate $\n{x_\Omega ^T(t)}_{X_\Omega}^2$ and $\n{\lambda_\Omega ^T (t)}_{X_\Omega}^2$ from above. These bounds can then be integrated to find bounds on $\n{x_\Omega ^T}_{L^2(0,T;X_\Omega)}^2$ and $\n{\lambda_\Omega ^T}_{L^2(0,T;X_\Omega)}^2$. This will give us an estimate of the form (leaving out $\Omega$ and $T$)
\begin{equation*}
        \n{x}_{2\land \infty}^2 + \n{\lambda}_{2\land \infty}^2 \!\leq \!\tilde{c} \left( \n{Cx}_{L^2(0,T;Y)}^2 \!+ \!\n{R^{-*}B^*\lambda}_{L^2(0,T;U)}^2 \!+ \!\n{l_1}_{1\lor 2}^2 \! + \! \n{\lambda _T}_X^2 \!+ \!\n{l_2}_{1\lor 2}^2 \!+\! \n{x_0}_X^2\right).
\end{equation*}
In the fourth step we will find an estimate for
$\n{C_\Omega x_\Omega^T}_{L^2(0,T;Y_\Omega)}^2 + \n{R_\Omega^{-*}B_\Omega^*\lambda _\Omega^T}_{L^2(0,T;U_\Omega)}^2$ which allows us to prove the desired result in the fifth and last step.

\begin{lemma}[Step 1]
\label{Step1Lemma}
Let the stabilizability/detectability condition in~\Cref{Thm: Boundedness} be fulfilled. Let $\varphi : [0,t] \rightarrow D(A_\Omega ^*)$ be the mild solution of~\eqref{PhiSystem} on the interval $[0,t]$ and $\psi: [t,T] \rightarrow D(A _\Omega^T)$ be the mild solution of~\eqref{PsiSystem} on the interval $[t,T]$.
Then there exist constants $M_\varphi, k_\varphi >0, M_\psi, k_\psi >0$ such that for all $\Omega \in \mathcal{O}$ and $T>0$ the estimates
\begin{equation}
\label{xEstimate}
    \hspace{-0.2cm} \forall  v \in L^2(0,t;X_\Omega): \int_0^t |\left \langle v(s), \varphi (s)\right\rangle _{X_\Omega} | \mathrm{d} s \leq \n{x_\Omega ^T(t)}_{X_\Omega} \frac{M_\varphi}{\sqrt{k_\varphi}} \sqrt{\int _0^t \n{v(s)}_{X_\Omega}^2e^{-k_\varphi (t-s)}\mathrm{d}s}
\end{equation}
and
\begin{equation}
\label{lambdaEstimate}
    \hspace{-0.2cm} \forall  w \in L^2(t,T;X_\Omega): \int_t^T |\left \langle w(s), \psi (s)\right\rangle _{X_\Omega} | \mathrm{d} s \leq \n{\lambda_\Omega ^T(t)}_{X_\Omega} \frac{M_\psi}{\sqrt{k_\psi}} \sqrt{\int _0^t \n{w(s)}_{X_\Omega}^2e^{-k_\psi (s-t)}\mathrm{d}s}
\end{equation}
hold true for $\varphi$ and $\psi$.
\end{lemma}
\begin{proof}
Due to domain-uniform stabilizability there exist constants $M_\varphi, k_\varphi >0$ such that the semigroup family $\left(\left(\mathcal{T}_\Omega^\varphi(t)\right)_{t \geq 0}\right)_{\Omega \in \mathcal{O}}$ generated by the operators $A_\Omega + B_\Omega K_\Omega ^B$ fulfills the estimate
\begin{equation*}
    \forall \, T>0 \, \forall \, \Omega \in \mathcal{O}: \n{\mathcal{T}_\Omega^\varphi(t)}_{L(X_\Omega, X_\Omega)} \leq M_\varphi e^{-k_\varphi t}.
\end{equation*}
During the remainder of the proof we leave out the indices $\Omega$ and $T$ for readability. To prove~\eqref{xEstimate} we first estimate the solution $\varphi \in C(0,t;X)$ of~\eqref{PhiSystem} via
\begin{equation*}
    \forall s \in [0,t]: \n{\varphi (s)}_X = \n{\mathcal{T}^\varphi(t-s)x(t)}_X \leq \n{\mathcal{T}^\varphi(t-s)}_{L(X,X)}\n{x(t)}_X \leq M_\varphi e^{-k_\varphi (t-s)} \n{x(t)}_X
\end{equation*}
where the last equation follows from the exponential stability of the semigroup $\left( \mathcal{T}^\varphi (t) \right)_{t\geq 0}$. This inequality leads to
\begin{flalign}
\label{Step1Estimate1}
    \int_0^t |\left \langle v(s), \varphi (s)\right\rangle _X | \mathrm{d} s \overset{\textrm{CSI}}{\leq} \int_0^t \n{v(s)}_X \n{\varphi (s)}_X  \mathrm{d} s \leq \n{x(t)}_X M_\varphi \int_0^t \n{v(s)}_X e^{-k_\varphi (t-s)}\mathrm{d} s
\end{flalign}
where the first inequality follows from the Cauchy-Schwarz-Inequality (CSI). The right hand side of~\eqref{Step1Estimate1} can be further estimated via
\begin{flalign}
\label{Step1Estimate2}
\begin{split}
    M_\varphi \int_0^t \n{v(s)}_X e^{-k_\varphi (t-s)}\mathrm{d} s &\overset{\textrm{\textcolor{white}{Hölder}}}{=} M_\varphi \int_0^t \n{v(s)}_X e^{-\frac{k_\varphi}{2} (t-s)} e^{-\frac{k_\varphi}{2} (t-s)}\mathrm{d} s\\
    &\overset{\textrm{\textcolor{white}{Hölder}}}{=} \n{M_\varphi  v e^{-\frac{k_\varphi}{2} (t-\cdot)} e^{-\frac{k_\varphi}{2} (t-\cdot)}}_{L^1(0,t;X)}\\
    &\overset{\textrm{Hölder}}{\leq} \n{M_\varphi  v e^{-\frac{k_\varphi}{2} (t-\cdot)} }_{L^2(0,t;X)} \n{e^{-\frac{k_\varphi}{2} (t-\cdot)}}_{L^2(0,t;X)}\\
   &\overset{\textrm{\textcolor{white}{Hölder}}}{=} M_\varphi \sqrt{\int_0^t \n{v(s)}_X^2 e^{-k_\varphi (t-s)}\mathrm{d} s} \sqrt{\int_0^t e^{-k_\varphi (t-s)}\mathrm{d} s}\\
   &\overset{\textrm{\textcolor{white}{Hölder}}}{=} M_\varphi \sqrt{\int_0^t \n{v(s)}_X^2 e^{-k_\varphi (t-s)}\mathrm{d} s} \sqrt{\frac{1}{k_\varphi} \left( 1 - e^{-k_\varphi t}\right)}\\
   &\overset{\textrm{\textcolor{white}{Hölder}}}{\leq} \frac{M_\varphi}{\sqrt{k_\varphi}} \sqrt{\int_0^t \n{v(s)}_X^2 e^{-k_\varphi (t-s)}\mathrm{d} s},
\end{split}
\end{flalign} 
where we used the Hölder inequality in the third step of~\eqref{Step1Estimate2}.
Plugging~\eqref{Step1Estimate2} into~\eqref{Step1Estimate1} yields~\eqref{xEstimate}. \eqref{lambdaEstimate} can be derived in an analogous way from domain-uniform detectability.
\end{proof}

To find an alternative representation for the terms $\n{x_\Omega ^T(t)}_X^2$ and $\n{\lambda_\Omega ^T(t)}_X^2$ we test the state equation from~\eqref{Pontryagin} with $\varphi$ and the adjoint equation with $\psi$. This leads to the following lemma.

\begin{lemma}[Step 2]
\label{Step2Lemma}
Let $T>0$ and $\Omega \in \mathcal{O}$ be arbitrary. Let the stabilizability/detectability condition in~\Cref{Thm: Boundedness} be fulfilled. If $\varphi : [0,t] \rightarrow D(A_\Omega^*)$ is the solution of~\eqref{PhiSystem} on the interval $[0,t]$ and $\psi: [t,T] \rightarrow D(A _\Omega)$ is the solution of~\eqref{PsiSystem} on the interval $[t,T]$ then we have the two equalities (leaving out the indices $\Omega$ and $T$)
\begin{equation}
    \label{Step2Aim1}
    \n{x(t)}_{X}^2 \!= \!\int _0^t  \left \langle l_2(s) , \varphi (s) \right\rangle _{X} \!-\!\left \langle K ^C C x(s) , \varphi (s)\right\rangle _{X} \!+\! \left \langle R ^{-*} B ^* \lambda (s) , R ^{-*} B ^* \varphi (s) \right\rangle _{X} \,\mathrm{d}s \!+ \! \left \langle x_0 , \varphi (0) \right\rangle _{X}
\end{equation}
and
\begin{equation}
    \label{Step2Aim2}
    \n{\lambda(t)}_{X}^2 \!= \!\int _0^t  \left \langle l_1(s) , \psi (s) \right\rangle _{X} \!-\!\left \langle \left(K ^{B}\right)^* B^* \lambda(s) , \psi (s)\right\rangle _{X}\! - \!\left \langle C x (s) ,C \psi (s) \right\rangle _{X} \,\mathrm{d}s \! + \! \left \langle \lambda , \psi (T) \right\rangle _{X}.
\end{equation}
\end{lemma}
\begin{proof}
For readability we leave out the indices $\Omega$ and $T$ during the proof. By testing the state equation from~\eqref{Pontryagin} with $\varphi (s)$ and subtracting $\left \langle K_\Omega ^C Cx(s), \varphi (s)\right\rangle$ on both sides of the resulting equation we find the equality
\begin{flalign*}
    &\left \langle \dot{x} (s), \varphi (s) \right\rangle _X - \left \langle (A+K ^C C)x (s) , \varphi (s) \right\rangle _X - \left \langle BQ^{-1} B^* \lambda (s), \varphi (s) \right\rangle _X\\ 
    = &\left \langle l_2(s), \varphi (s) \right\rangle _X - \left \langle K ^C C x(s), \varphi (s)\right\rangle _X.
\end{flalign*}
By integrating from $0$ to $t$ and rearranging the terms we find
\begin{flalign*}
    &\int _0^t \left \langle \dot{x} (s), \varphi (s) \right\rangle _X - \left \langle (A+K ^C C)x (s) , \varphi (s) \right\rangle _X \mathrm{d}s\\ 
    =& \int _0^t \left \langle\dot{x} (s), \varphi (s) \right\rangle _X - \left \langle x (s) , (A^* + C^* \left(K ^C\right)^*)\varphi (s) \right\rangle _X \mathrm{d}s\\
    =&\int _0^t \left \langle \dot{x} (s), \varphi (s) \right\rangle _X + \left \langle x (s) , \dot{\varphi} (s) \right\rangle _X \mathrm{d}s\\
    =& \int _0^t \frac{\mathrm{d}}{\mathrm{d}s}\left \langle x (s), \varphi (s) \right\rangle _X \mathrm{d}s\\
    =& \left \langle x (t), \varphi (t) \right\rangle _X - \left \langle x (0), \varphi (0) \right\rangle _X\\
    =& \n{x(t)}_X^2 - \left \langle x _0, \varphi (0) \right\rangle _X\\
    =& \int _0^t \left \langle l_2(s), \varphi (s) \right\rangle _X - \left \langle K ^C C x(s), \varphi (s)\right\rangle _X + \left \langle BQ^{-1} B^* \lambda (s), \varphi (s) \right\rangle _X \mathrm{d}s\\
    =&\int _0^t  \left \langle l_2(s) , \varphi (s) \right\rangle _X -\left \langle K ^C C x(s) , \varphi (s)\right\rangle _X + \left \langle R^{-*} B^* \lambda (s) , R^{-*} B^* \varphi (s) \right\rangle _X \,\mathrm{d}s.
\end{flalign*}
This shows~\eqref{Step2Aim1}.~\eqref{Step2Aim2} follows by testing the costate equation with $\psi$, subtracting the term $\left \langle \left( K ^B\right)^* B^* \lambda(s), \psi (s)\right\rangle$ on both sides and proceeding in the same way as before.
\end{proof}

By plugging the estimate from~\Cref{Step1Lemma} (Step 1) into the equality from~\Cref{Step2Lemma} we can find a first estimate on $\n{x_\Omega ^T}_{2\land \infty}^2 + \n{\lambda _\Omega ^T}_{2\land \infty}^2$.
This results in the following Lemma.

\begin{lemma}[Step 3]
\label{Step3Lemma}
Let the stabilizability/detectability condition in~\Cref{Thm: Boundedness} be fulfilled. Then there exists a constant $\tilde{c}>0$ such that the inequality 
\begin{flalign}
\label{Step3Aim}
\begin{split}
        \n{x_\Omega ^T}_{2\land \infty}^2 &+ \n{\lambda_\Omega ^T}_{2\land \infty}^2\\
        \leq &\,\tilde{c} \left( \n{C_\Omega x_\Omega ^T}_{L^2(0,T;Y_\Omega)}^2 \!+\! \n{R_\Omega^{-*}B_\Omega^*\lambda_\Omega ^T}_{L^2(0,T;U_\Omega)}^2 \!+\! \n{l_1}_{1\lor 2}^2 \!+\! \n{\lambda^\Omega _T}_{X_\Omega}^2 \!+\! \n{l_2}_{1\lor 2}^2 \!+\! \n{x^\Omega _0}_{X_\Omega}^2\right)
\end{split}
\end{flalign}
is fulfilled for all $T>0$ and $\Omega \in \mathcal{O}$.
\end{lemma}

\begin{proof}
The proof is organised into four parts. First we show an estimate of the form
\begin{equation}
    \label{Step31Aim}
    \n{x_\Omega ^T}_{C(0,T;X_\Omega)}^2 \leq c_1 \left( \n{C_\Omega x_\Omega ^T}_{L^2(0,T;Y_\Omega)}^2 + \n{R_\Omega^{-*}B_\Omega^*\lambda _\Omega ^T}_{L^2(0,T;U_\Omega)}^2+ \n{l_2}_{1\lor 2}^2 + \n{x^\Omega _0}_{X_\Omega}^2\right).
\end{equation}
Then, we show an analogous estimate on the norm
$\n{\lambda _\Omega^T}_{C(0,T;X_\Omega)}^2$.
In the third and fourth part we show estimates on the $L^2$-norms of $x_\Omega^T$ an $\lambda _\Omega ^T$. As before we leave out the indices $\Omega$ and $T$ during the proof for readability. 
\\We start with the first part.
By using the results of Step 1 and Step 2 we find the estimate
\begin{flalign*}
    \n{x(t)}_X^2
    \overset{\textrm{Step 2}}{=}& \int _0^t  \left \langle l_2(s) , \varphi (s) \right\rangle _X -\left \langle K ^C C x(s) , \varphi (s)\right\rangle _X + \left \langle R^{-*} B^* \lambda (s) , R^{-*} B^* \varphi (s) \right\rangle _X \,\mathrm{d}s + \left \langle x_0 , \varphi (0) \right\rangle _X\\
    \overset{\textrm{Step 1}}{\leq}& \n{x(t)}_X \frac{M_\varphi}{\sqrt{k_\varphi}} 
    \n{K ^C}_{L(Y,X)} \sqrt{\int _0^t \n{Cx(s)}^2e^{-k_\varphi (t-s)}\mathrm{d}s} \,\, + \,\, \n{x_0}_X \n{x(t)}_X M_\varphi \sqrt{e^{-k_\varphi t}}\\
    \,+ \,\,\,&\,\, \n{x(t)}_X \frac{M_\varphi}{\sqrt{k_\varphi}} \n{BR^{-1}}_{L(U,X)} \sqrt{\int _0^t \n{R^{-*}B^*\lambda (s)}_U^2 e^{-k_\varphi (t-s)}\mathrm{d}s}\\
    \,+ \,\,\,&\,\, \min \left\{ 
    \n{x(t)}_X M_\varphi \int _0^t \n{l_2(s)}_X  e^{-k_\varphi (t-s)}\mathrm{d}s, \n{x(t)}_X \frac{M_\varphi}{\sqrt{k_\varphi}}\sqrt{\int_0^t \n{l_2(s)}^2  e^{-k_\varphi (t-s)}\mathrm{d}s}
    \right\}.
\end{flalign*}
Dividing by $\n{x(t)}$ leads to the inequality
\begin{flalign*}
    \n{x(t)}_X^2 \leq & \frac{M_\varphi}{\sqrt{k_\varphi}}  \left( \n{K ^C} \sqrt{\int _0^t \n{Cx(s)}^2e^{-k_\varphi (t-s)}\mathrm{d}s} \n{x(t)}_X\right.\\
    &\,\,\,\,\,+ \,\,\,\,\, \left. \n{BR^{-1}}_U \sqrt{\int _0^t \n{R^{-*}B^*\lambda (s)}^2 e^{-k_\varphi (t-s)}\mathrm{d}s} \right)\\
    &\,\,\,\,\,+ \,\,\,\,\, \min \left\{  M_\varphi \alpha(t),  \frac{M_\varphi}{\sqrt{k_\varphi}}\beta(t)
    \right\} + \n{x_0}_X M_\varphi \sqrt{e^{-k_\varphi t}}\\
    &\leq \hat{c}_1\left( \sqrt{\int _0^t \n{Cx(s)}^2e^{-k_\varphi (t-s)}\mathrm{d}s} +  \sqrt{\int _0^t \n{R^{-*}B^*\lambda (s)}^2 e^{-k_\varphi (t-s)}\mathrm{d}s} \right)\\
    &+ \hat{c}_1 \left(\min \left\{ \alpha(t), \beta(t) \right\} + \n{x_0}_X\sqrt{e^{-k_\varphi t}}
    \right),
\end{flalign*}
where we define functions $\alpha,\beta$ by $\alpha(t):=\int _0^t \n{l_2(s)}_X  e^{-k_\varphi (t-s)}\mathrm{d}s$, $\beta(t):=\sqrt{\int_0^t \n{l_2(s)}^2  e^{-k_\varphi (t-s)}\mathrm{d}s}$ and a constant $\hat{c}_1:=M_\varphi \max \left\{ 1, \frac{1}{\sqrt{k_\varphi}}\right\}\max \left\{ 1, C_\psi, \frac{C_B}{\alpha}\right\}$. By taking the square on both sides of the inequality and using the estimate
\begin{equation*}
    \forall a,b,c,d \in \R: (a+b+c+d)^2 \leq 2((a+b)^2+(c+d)^2) \leq 4(a^2+b^2+c^2+d^2)
\end{equation*}
we find the estimate
\begin{flalign}
\label{xEstimate3}
\begin{split}
    \n{x(t)}_X &\leq c_1\left( \int _0^t \n{Cx(s)}^2e^{-k_\varphi (t-s)}\mathrm{d}s +  \int _0^t \n{R^{-*}B^*\lambda (s)}^2 e^{-k_\varphi (t-s)}\mathrm{d}s\right)\\
    &+ c_1\left(\min \left\{ \alpha(t)^2, \beta(t)^2 \right\} +\n{x_0}_X^2 e^{-k_\varphi t}
    \right)
\end{split}
\end{flalign}
with constant $c_1:=4M_\varphi \max \left\{ 1, \frac{1}{\sqrt{k_\varphi}}\right\}\max \left\{ 1, C_\psi, \frac{C_B}{\alpha}\right\}$. This implies
\begin{flalign*}
\begin{split}
    \n{x}_{C(0,T;X)}^2 &= \underset{t\in [0,T]}{\max} \n{x(t)}_X^2 \leq c_1 \left(\n{Cx}_{L^2(0,T;Y)}^2 + \n{R^{-*}B^* \lambda}_{L^2(0,T;U)}^2\right) \\
    &+ c_1 \left(\min \left\{\n{l_2}_{L^1{0,T;X}^2}, \n{l_2}_{L^2{0,T;X}^2} \right\} + \n{x_0}_X^2\right)
\end{split}
\end{flalign*}
and therefore~\eqref{Step31Aim}.\\\\

The aim of the second part is to show an estimate of the form
\begin{equation}
    \label{Step32Aim}
    \n{\lambda}_{C(0,T;X)}^2 \leq c_2 \left( \n{Cx}_{L^2(0,T;Y)}^2 + \n{R^{-*}B^*\lambda}_{L^2(0,T;U)}^2+ \n{l_1}_{1\lor 2}^2 + \n{\lambda _T}_X^2\right).
\end{equation}
Again using the results of Step 1 and Step 2 we find the estimate
\begin{flalign*}
    &\n{\lambda(t)}_X^2\\
    \overset{\textrm{Step 2}}{=} &\int _0^t  \left \langle l_1(s) , \psi (s) \right\rangle _X -\left \langle \left(K ^B\right)^* B^* \lambda(s) , \psi (s)\right\rangle _X - \left \langle C x (s) ,C \psi (s) \right\rangle _X \,\mathrm{d}s + \left \langle \lambda _T , \psi (T) \right\rangle _X\\
    \overset{\textrm{Step 1}}{\leq} &\n{\lambda(t)}_X \frac{M_\psi}{\sqrt{k_\psi}} 
    \n{K ^B}_{L(X,U)} \sqrt{\int _t^T \n{B^*\lambda (s)}_U^2e^{-k_\psi (s-t)}\mathrm{d}s} \,\, + \,\, \n{\lambda _T}_X \n{\lambda(t)}_X M_\psi \sqrt{e^{-k_\psi t}}\\
    \,\,\,\,\,+ \,\,\,&\,\, \n{\lambda(t)}_X \frac{M_\psi}{\sqrt{k_\psi}} \n{C}_{L(X,Y)} \sqrt{\int _t^T \n{Cx(s)}_X^2 e^{-k_\psi (s-t)}\mathrm{d}s}\\
    \,\,\,\,\,+ \,\,\,&\,\, \min \left\{ 
    \n{\lambda(t)}_X M_\psi \int _t^T \n{l_1(s)}_X  e^{-k_\psi (s-t)}\mathrm{d}s,  \n{\lambda(t)}_X\frac{M_\psi}{\sqrt{k_\psi}}\sqrt{\int_t^T \n{l_1(s)}^2  e^{-k_\psi (s-t)}\mathrm{d}s}
    \right\}.
\end{flalign*}
Dividing by $\n{\lambda(t)}_X$ leads to the inequality
\begin{flalign*}
    &\n{\lambda(t)}_X\\
    \leq \,&\frac{M_\psi}{\sqrt{k_\psi}} \left(
    \n{K ^B}_{L(X,U)} \sqrt{\int _t^T \n{B^*\lambda (s)}_U^2e^{-k_\psi (s-t)}\mathrm{d}s} + \n{C}_{L(X,Y)} \sqrt{\int _t^T \n{Cx(s)}_X^2 e^{-k_\psi (s-t)}\mathrm{d}s}\right)\\
    + \,&\,\, \min \left\{ 
    M_\psi \gamma (t),  \frac{M_\psi}{\sqrt{k_\psi}}\delta (t)
    \right\} \,\, + \,\,  \n{\lambda _T}_X M_\psi \sqrt{e^{-k_\psi t}}\\
    \leq \, &\hat{c}_2\left( \sqrt{\int _t^T \n{Cx(s)}^2e^{-k_\psi (s-t)}\mathrm{d}s} +  \sqrt{\int _t^T \n{R^{-*}B^*\lambda (s)}^2 e^{-k_\psi (s-t)}\mathrm{d}s}\right)\\
    + \,& \hat{c}_2\left(\min \left\{ \gamma(t), \delta(t) \right\} + \n{\lambda _T}_X\sqrt{e^{-k_\psi t}}
    \right),
\end{flalign*}
where we define a constant $\hat{c}_2:=M_\psi \max \left\{ 1, \frac{1}{\sqrt{k_\psi}}\right\}\max \left\{ 1, C_\varphi C_R, C_C\right\}$ and functions $\gamma(t):=\int _t^T \n{l_1(s)}_X  e^{-k_\psi (s-t)}\mathrm{d}s$, $\delta(t):=\sqrt{\int_t^T \n{l_1(s)}^2  e^{-k_\psi (s-t)}\mathrm{d}s}$. Proceeding in the same way as in the first part we find the estimate
\begin{flalign}
\label{lambdaEstimate3}
\begin{split}
    \n{\lambda(t)}_X &\leq c_2\left( \int _t^T \n{Cx(s)}^2e^{-k_\psi (t-s)}\mathrm{d}s +  \int _t^T \n{R^{-*}B^*\lambda (s)}^2 e^{-k_\psi (t-s)}\mathrm{d}s\right)\\
    &+ c_2\left(\min \left\{ \gamma(t)^2, \delta(t)^2 \right\} +\n{\lambda_T}_X^2 e^{-k_\psi t}
    \right)
\end{split}
\end{flalign}
with constant $c_1:=4M_\psi \max \left\{ 1, \frac{1}{\sqrt{k_\psi}}\right\}\max \left\{ 1, C_\varphi C_R, C_C\right\}$. This implies
\begin{flalign*}
    &\n{\lambda}_{C(0,T;X)}^2\\
    =& \underset{t\in [0,T]}{\max} \n{\lambda(t)}_X^2\\
    \leq&  \, c_2 \left(\n{Cx}_{L^2(0,T;Y)}^2 + \n{R^{-*}B^* \lambda}_{L^2(0,T;U)}^2 + \min \left\{\n{l_1}_{L^1{0,T;X}^2}, \n{l_1}_{L^2{0,T;X}^2} \right\} + \n{\lambda _T}_X^2\right)
\end{flalign*}
and therefore~\eqref{Step32Aim}.\\

\noindent For the $L^2$-norm of $x$ we prove an estimate of the form
\begin{equation}
    \label{Step33Aim}
    \n{x}_{L^2(0,T;X)}^2 \leq c_3 \left( \n{Cx}_{L^2(0,T;Y)}^2 + \n{R^{-*}B^*\lambda}_{L^2(0,T;U)}^2+ \n{l_2}_{1\lor 2}^2 + \n{x_0}_X^2\right).
\end{equation}
We first define for $k >0$
\begin{equation*}
    g_k: \R_{\geq 0} \rightarrow \R _{\geq 0}, \platz g_k(t):=e^{-kt}
\end{equation*}
and note for all $T>0$ that
\begin{flalign}
    \label{gnorm}
    \begin{split}
     \n{g_k}_{L^1(0,T;\R)} &= \int _0^T e^{-kt} \mathrm{dt} \leq \n{g_k}_{L^1(0,\infty;\R)}=\int _0^\infty e^{-kt} \mathrm{dt}=\left[-\frac{1}{k} e^{-kt} \right]_0^\infty = \frac{1}{k}\\
    \n{g_k}_{L^2(0,T;\R)}^2 &= \int _0^T e^{-2kt} \mathrm{dt} \leq \n{g_k}_{L^2(0,\infty;\R)}^2=\int _0^\infty e^{-2kt} \mathrm{dt}=\left[-\frac{1}{2k} e^{-2kt} \right]_0^\infty = \frac{1}{2k}.
    \end{split}
\end{flalign}
From~\eqref{xEstimate3} we conclude
\begin{flalign}
\label{xEstimate32}
    \begin{split}
        \n{x}_{L^2(0,T;X)}^2
        &= \int _0^T \n{x(t)}_X^2 dt\\
        &\!\!\!\overset{~\eqref{xEstimate3}}{\leq} c_1\ \left(\int _0^T \int _0^t \left(\n{Cx(s)}^2+ \n{R^{-*}B^*\lambda (s)}^2 \right) e^{-k_\varphi (t-s)}\mathrm{d}s \, \mathrm{d}t\right)\\
        &+ c_1\ \left(\int _0^T \min \left\{ \alpha(t)^2, \beta(t)^2 \right\}\mathrm{d}t
        +\int _0^T \n{x_0}_X^2 e^{-k_\varphi t}\mathrm{d}t\right).
    \end{split}
\end{flalign}
We will now estimate the different terms from the right hand side of this inequality individually. Defining the function
\begin{equation*}
    h_\varphi :[0,T] \rightarrow \R_{\geq 0}, \platz h_\varphi (t):= \n{Cx(t)}^2+ \n{R^{-*}B^*\lambda (t)}^2
\end{equation*}
we find the estimate
\begin{flalign*}
    \begin{split}
        \int _0^T \int _0^t \left(\n{Cx(s)}^2+ \n{R^{-*}B^*\lambda (s)}^2 \right) e^{-k_\varphi (t-s)}\mathrm{d}s \, \mathrm{d}t &= \n{h_\varphi * g_{k_\varphi }}_{L^1(0,T;\R_{\geq 0})}\\
        &\leq \n{h_\varphi}_{L^1(0,T;\R_{\geq 0})} \n{g_{k_\varphi}}_{L^1(0,T;\R _{\geq 0})}\\
        &= \frac{1}{k_\varphi} \int _0^T \n{Cx(t)}_Y^2+ \n{R^{-*}B^*\lambda (t)}_U^2\mathrm{d}t\\
        &= \frac{1}{k_\varphi} \left( \n{Cx}_{L^2(0,T;Y)}^2 + \n{R^{-*}B^*\lambda}_{L^2(0,T;U)}^2 \right)
    \end{split}
\end{flalign*}
by using Young's inequality and~\eqref{gnorm}. For the second term we find the estimates
\begin{flalign*}
    \begin{split}
       \int _0^T \alpha (t)^2 \mathrm{d}t &= \int _0^T \int_0^t \n{l_2(s)}_X^2  e^{-k_\varphi (t-s)}\mathrm{d}s \mathrm{d}t = \n{\left(\n{l_2(\cdot)}_X^2 * g_{k_\varphi}\right)}_{L^1(0,T;\R_{\geq 0})}\\
       &\leq \n{\n{l_2(\cdot)}_X^2}_{L^1(0,T;\R_{\geq 0})} \n{g_{k_\varphi}}_{L^1(0,T;\R _{\geq 0})} = 
       \frac{1}{k_\varphi}\n{l_2}_{L^2(0,T;X)}^2
    \end{split}
\end{flalign*}
and
\begin{flalign*}
    \begin{split}
       \int _0^T \beta (t)^2 \mathrm{d}t &= \int _0^T \left(\int_0^t \n{l_2(s)}_X  e^{-k_\varphi (t-s)}\mathrm{d}s\right)^2 \mathrm{d}t = \n{\left(\n{l_2(\cdot)}_X * g_{k_\varphi}\right)}_{L^2(0,T;\R_{\geq 0})}^2\\
       &\leq \n{\n{l_2(\cdot)}_X}_{L^1(0,T;\R_{\geq 0})}^2 \n{g_{k_\varphi}}_{L^2(0,T;\R _{\geq 0}^2)}^2 = 
       \frac{1}{2k_\varphi}\n{l_2}_{L^1(0,T;X)}^2
    \end{split}
\end{flalign*}
For the third term we find
\begin{equation*}
    \int _0^T \n{x_0}_X^2 e^{-k_\varphi t}\mathrm{d}t = \n{x_0}_X^2\int _0^T  g_{k_\varphi}(t)\mathrm{d}t = \n{x_0}_X^2 \n{g_{k_\varphi}}_{L^1(0,T;\R_{\geq 0})} = \frac{1}{k_\varphi} \n{x_0}_X^2.
\end{equation*}
Plugging these three estimates into~\eqref{xEstimate32} yields~\eqref{Step33Aim} with
\begin{equation*}
    c_3:= \frac{1}{k_\varphi} c_1
\end{equation*}
as the sought for constant.\\\\

By proceeding analogously as in the third part we can find the estimate
\begin{equation}
    \label{Step34Aim}
    \n{\lambda}_{L^2(0,T;X)}^2 \leq c_4 \left( \n{Cx}_{L^2(0,T;Y)}^2 + \n{R^{-*}B^*\lambda}_{L^2(0,T;U)}^2+ \n{l_1}_{1\lor 2}^2 + \n{\lambda _T}_X^2\right).
\end{equation}
with
\begin{equation*}
    c_4 = \frac{1}{k_\psi} c_2
\end{equation*}
in the fourth part of this proof.
Combining the four estimates~\eqref{Step31Aim},~\eqref{Step32Aim},~\eqref{Step33Aim} and~\eqref{Step34Aim} implies~\eqref{Step3Aim} with
\begin{equation*}
    \tilde{c}= \max \{1, \frac{1}{k_\varphi}, \frac{1}{k_\psi}\} \max \{c_1, c_2\}.
\end{equation*}
Note that $\tilde{c}$ only depends on the quantities $C_B$, $C_C$, $C_\varphi$, $C_\psi$, $M_\varphi$, $M_\psi$, $k_\varphi$, $k_\psi$, $C_R$ and $\alpha$. By assumption this means that the estimates hold for all $T>0$ and $\Omega \in \mathcal{O}$.
\end{proof}

By testing the state equation from~\eqref{Pontryagin} with $\lambda _\Omega ^T$ and the adjoint equation with $x _\Omega ^T$ we can further estimate the right hand side of the inequality which was shown in~\Cref{Step3Lemma}. This way we succeed in linking $\n{z _\Omega ^T}_{2\land \infty}$ to the components of $r _\Omega ^T$ in the following lemma.
For technical reasons we will need the following Lemma which was taken from~\cite{Schaller2019}.
\begin{lemma}[Generalized Cauchy-Schwarz Inequality]
\label{CSULemma}
Let $X$ be a Hilbert space and $T>0$. For all $v\in C(0,T;X)$ and for all $w \in L^1(0,T;X)$ we have the inequality 
\begin{equation*}
    \int _0^T \left \langle v(s), w(s)\right\rangle \mathrm{d}s \leq \n{v}_{2\land \infty} \n{w} _{1 \lor 2}.
\end{equation*}
\end{lemma}
With the help of~\Cref{CSULemma} we are now able to complete the fourth step in showing boundedness for the solution operator.
\begin{lemma}[Step 4]
\label{Step4Lemma}
Let $x _\Omega ^T$, $\lambda _\Omega ^T$, $l_1$, $l_2$, $\lambda ^\Omega _T$ and $x ^\Omega_0$ be as in~\eqref{Pontryagin}. Then the estimate
\begin{flalign}
    \label{Step4Aim}
    \begin{split}
     &\n{C _\Omega x _\Omega ^T}_{L^2(0,T;Y _\Omega)}^2 
        + \n{R _\Omega^{-*}B _\Omega^*\lambda _\Omega ^T}_{L^2(0,T;U _\Omega)}^2\\
        \leq  &\n{\lambda ^\Omega _T}_{X _\Omega}\n{x _\Omega ^T(T)}_{X _\Omega} + 
        \n{x^\Omega _0}_{X _\Omega}\n{\lambda _\Omega ^T (0)}_{X _\Omega} + \n{l_2}_{1 \lor 2} \n{\lambda _\Omega ^T}_{2\land \infty} + \n{l_1}_{1 \lor 2} \n{x _\Omega ^T}_{2\land \infty}
    \end{split}
\end{flalign}
is fulfilled for all $T>0$ and $\Omega \in \mathcal{O}$.
\end{lemma}
\begin{proof}
For readability we leave out the indices $\Omega$ and $T$. Testing the adjoint equation from~\eqref{Pontryagin} with $x$ yields
\begin{equation}
\label{TestedAdjointEquation}
    \forall s \in [0,T]: \left \langle C^* Cx(s), x(s)\right\rangle _X - \left \langle \dot{\lambda}(s), x(s)\right\rangle _X - \left \langle A^* \lambda (s), x(s)\right\rangle _X = \left \langle l_1(s), x(s)\right\rangle _X.
\end{equation}
Analogously we find
\begin{equation}
\label{TestedStateEquation}
    \forall s \in [0,T]: \left \langle \dot{x}(s), \lambda (s)\right\rangle _X - \left \langle A x(s), \lambda(s)\right\rangle _X - \left \langle B(R^*R)^{-1}B^*\lambda (s), \lambda (s)\right\rangle _X= \left \langle l_2(s), \lambda(s)\right\rangle _X.
\end{equation}
by testing the state equation from~\eqref{Pontryagin} with $\lambda$.
By subtracting~\eqref{TestedStateEquation} from~\eqref{TestedAdjointEquation} and integrating from $0$ to $T$ we find
\begin{flalign*}
    \n{Cx}_{L^2(0,T;Y)}^2 &+ \n{R^{-*}B^*\lambda}_{L^2(0,T;U)}^2\\
    =&\int _0^T \left \langle B(R^*R)^{-1}B^*\lambda (s), \lambda (s)\right\rangle _X + \left \langle C^* Cx(s), x(s)\right\rangle _X \mathrm{d}s\\ 
    = &\int _0^T \left \langle \dot{\lambda}(s), x(s)\right\rangle _X + \left \langle \dot{x}(s), \lambda (s)\right\rangle _X \mathrm{d}s + \int _0^T \left \langle l_1(s), x(s)\right\rangle _X - \left \langle l_2(s), \lambda(s)\right\rangle _X \mathrm{d}s\\
    = &\int _0^T \frac{\mathrm{d}}{\mathrm{d}t} \left \langle \lambda(s), x(s)\right\rangle _X \mathrm{d}s + \int _0^T \left \langle l_1(s), x(s)\right\rangle _X - \left \langle l_2(s), \lambda(s)\right\rangle _X \mathrm{d}s\\
    = &\left \langle \lambda _T, x(T)\right\rangle _X - \left \langle x_0, \lambda (0)\right\rangle _X + \int _0^T \left \langle l_1(s), x(s)\right\rangle _X - \left \langle l_2(s), \lambda(s)\right\rangle _X \mathrm{d}s\\
    \leq &\n{\lambda _T}_{X}\n{x(T)}_{X} + \n{x_0}_{X}\n{\lambda (0)}_{X} + \n{l_2}_{1 \lor 2} \n{\lambda}_{2\land \infty} + \n{l_1}_{1 \lor 2} \n{x}_{2\land \infty}.
\end{flalign*}
The last inequality holds due to~\Cref{CSULemma}.
\end{proof}

\begin{lemma}[Step 5]
\label{Step5Lemma}
Let the stabilizability/detectability condition in~\Cref{Thm: Boundedness} be fulfilled. Assume that there exists a $\tilde{c}>0$ such that the inequality
\begin{flalign}
    \label{Lemma5Ass}
    \begin{split}
    \n{z _\Omega ^T}_{2\land \infty}^2 &\leq  \tilde{c} \left( \n{\lambda ^\Omega _T}_{X_\Omega}\n{x_\Omega ^T(T)}_{X_\Omega} + \n{x^\Omega _0}_{X_\Omega}\n{\lambda _\Omega ^T (0)}_{X_\Omega}\right)\\
    &+ \tilde{c} \left(\n{l_2}_{1 \lor 2} \n{\lambda _\Omega ^T}_{2\land \infty} + \n{l_1}_{1 \lor 2} \n{x _\Omega ^T}_{2\land \infty} + \n{r _\Omega ^T}_{1\lor 2}^2\right)
    \end{split}
\end{flalign}
holds for all $T>0$ and $\Omega \in \mathcal{O}$. Then there also exists a $c>0$ such that
\begin{equation}
    \label{Step5Aim}
        \forall \, T>0 \, \forall \,\Omega \in \mathcal{O}: \n{z _\Omega ^T}_{2\land \infty}^2 \leq c \, \n{r _\Omega ^T}_{1\lor 2}^2.
\end{equation}
\end{lemma}
\begin{proof}
    By first using the inequalities
    \begin{equation*}
        \n{x(T)}_X \leq \n{x}_{\infty} \leq \n{x}_{2\land \infty} \,\,\,\, \textrm{and} \,\,\,\, \n{\lambda(0)}_X \leq \n{\lambda}_{\infty} \leq \n{\lambda}_{2\land \infty}
        \end{equation*}
    and then adding the positive term
    \begin{equation*}
        \left(\n{\lambda _T}_{X} + \n{l_1}_{1 \lor 2}\right) \n{\lambda}_{2\land \infty} +   \left(\n{x_0}_{X} + \n{l_2}_{1 \lor 2}\right)\n{x}_{2\land \infty}
    \end{equation*}
    to the right hand side of~\eqref{Lemma5Ass} we get the inequalities
    \begin{flalign}
    \label{Ineq51}
        \begin{split}
            \n{z}_{2\land \infty}^2  &\leq  \tilde{c} \left( \left(\n{\lambda _T}_{X} + \n{l_1}_{1 \lor 2}\right)\n{x}_{2\land \infty} + \left(\n{x_0}_{X} + \n{l_2}_{1 \lor 2}\right) \n{\lambda}_{2\land \infty} + \n{r}_{1\lor 2}^2\right)\\
            & \leq \tilde{c} \left( \left(\n{\lambda _T}_{X} + \n{l_1}_{1 \lor 2} + \n{x_0}_{X} + \n{l_2}_{1 \lor 2}\right) \left(\n{x}_{2\land \infty} + \n{\lambda}_{2\land \infty}\right) + \n{r}_{1\lor 2}^2\right).
        \end{split}
     \end{flalign}
    By applying the inequality $(a+b)^2\leq 2(a^2 + b^2)$ we find the estimate
    \begin{flalign*}
        &\left(\n{\lambda _T}_{X} + \n{l_1}_{1 \lor 2} + \n{x_0}_{X} + \n{l_2}_{1 \lor 2}\right)^2 \left(\n{x}_{2\land \infty} + \n{\lambda}_{2\land \infty}\right)^2\\
        &\qquad \qquad \leq  8 \left(\n{\lambda _T}_{X}^2 + \n{l_1}_{1 \lor 2}^2 + \n{x_0}_{X}^2 + \n{l_2}_{1 \lor 2}^2\right) \left(\n{x}_{2\land \infty}^2 + \n{\lambda}_{2\land \infty}^2\right)
    \end{flalign*}
    which implies
    \begin{flalign*}
        &\left(\n{\lambda _T}_{X} + \n{l_1}_{1 \lor 2} + \n{x_0}_{X} + \n{l_2}_{1 \lor 2}\right) \left(\n{x}_{2\land \infty} + \n{\lambda}_{2\land \infty}\right)\\
        & \qquad \leq  \sqrt{8} \sqrt{\left(\n{\lambda _T}_{X}^2 + \n{l_1}_{1 \lor 2}^2 + \n{x_0}_{X}^2 + \n{l_2}_{1 \lor 2}^2\right)} \sqrt{\left(\n{x}_{2\land \infty}^2 + \n{\lambda}_{2\land \infty}^2\right)}\\
        &\qquad =  \sqrt{8} \n{r}_{1\lor 2} \n{z}_{2 \land \infty }.
    \end{flalign*}
    Plugging this into~\eqref{Ineq51} and using the inequality
    \begin{equation*}
        \sqrt{8}\tilde{c} \n{r}_{1\lor 2} \n{z}_{2 \land \infty } \leq \frac{1}{2} \left( 8\tilde{c}^2 \n{r}_{1\lor 2}^2  + \n{z}_{2 \land \infty }^2 \right)
    \end{equation*}
    we find
    \begin{flalign*}
        \n{z}_{2\land \infty}^2 &= \tilde{c} \left( \sqrt{8}\n{r}_{1\lor 2} \n{z}_{2 \land \infty } + \n{r}_{1\lor 2}^2\right) \\
        &\leq \left( 4\tilde{c}^2 + \sqrt{8} \tilde{c}\right) \n{r}_{1\lor 2}^2 + \frac{1}{2}\n{z}_{2\land \infty}^2.
    \end{flalign*}
    This implies~\eqref{Step5Aim} with constant
    \begin{equation*}
        c:=2\left( 4 \tilde{c}^2 + \sqrt{8} \tilde{c} \right ) .
    \end{equation*}
\end{proof}
By combining~\Cref{Step3Lemma},~\Cref{Step4Lemma} and~\Cref{Step5Lemma} we can finally prove the desired result on the boundedness of the solution operator $\M^{-1}$.

\begin{thm}
    Let the stabilizability/detectability condition in~\Cref{Thm: Boundedness} be fulfilled. Then there exists a constant $c>0$ such that the norm of the solution operator $\M^{-1}$ can be estimated by
    \begin{equation*}
        \forall\, T>0 \, \forall \, \Omega \in \mathcal{O}: \n{\M^{-1}}_{L(W^{1\lor 2},W^{2\land \infty})}\leq c.
    \end{equation*}
\end{thm}
\begin{proof}
    Using the estimates in~\Cref{Step3Lemma} and~\Cref{Step4Lemma} we find the inequalities
\begin{flalign*}
    \n{z}_{2\land \infty}^2 &\overset{\textrm{\textcolor{white}{Lemma 11}}}{=} \n{x}_{2\land \infty}^2 + \n{\lambda}_{2\land \infty}^2\\
    &\overset{\textrm{Lem.~\ref{Step3Lemma}}}{\leq} \tilde{c} \left( \n{Cx}_{L_2(Y)}^2 + \n{R^{-*}B^* \lambda}_{L_2(U)}^2 + \n{r}_{1\lor 2}^2\right)\\
    &\overset{\textrm{Lem.~\ref{Step4Lemma}}}{\leq} \tilde{c} \left( \n{\lambda _T}_{X}\n{x(T)}_{X} \!+\! \n{x_0}_{X}\n{\lambda (0)}_{X} \!+\! \n{l_2}_{1 \lor 2} \n{\lambda}_{2\land \infty} \!+\! \n{l_1}_{1 \lor 2} \n{x}_{2\land \infty} \!+\! \n{r}_{1\lor 2}^2\right).
\end{flalign*}
This implies that the requirements of~\Cref{Step5Lemma} are fulfilled. Therefore we know that there exists a $c>0$ independent of $T>0$ and $\Omega \in \mathcal{O}$ such that
\begin{equation*}
        \n{z}_{2\land \infty}^2 \leq c \, \n{r}_{1\lor 2}^2.
\end{equation*}
This proves our claim.
\end{proof}

\newpage
\section{Proof of Lemmata~\ref{Lem: SolutionFormulaStateFeedback1}-\ref{Lem: IntervalCondition}}
\label{App: Lemmata}
In the following we use the notation $B_j =b_j-a_j$ and $A_j = a_{j+1}-b_j$.
\paragraph{Proof of~\Cref{Lem: SolutionFormulaStateFeedback1}}
\begin{proof}
    For $x^0 \in \mathrm{dom}(A)$ the fundamental theorem of calculus leads to
    \begin{flalign*}
        \frac{\partial}{\partial \omega} x(\omega, t)\! =\! &-\frac{1}{c} \left(P_{\Omega _L}(k)(\omega)\!-\!P_{\Omega _L}(k)(\omega \!-\! ct) \right)e^{-\frac{1}{c}\int_{\omega-ct}^{\omega}P_{\Omega _L}(k)(y)\mathrm{d}y} P_{\Omega _L}(x^0)(\omega-ct)\\
        &+ e^{-\frac{1}{c}\int_{\omega-ct}^{\omega}P_{\Omega _L}(k)(y)\mathrm{d}y} P_{\Omega _L}\left(\frac{\mathrm{d}}{\mathrm{d} \omega}x^{0}\right)(\omega-ct)
    \end{flalign*}
    and
    \begin{flalign*}
        &-\frac{\partial}{\partial t} x(\omega, t)\\
        = &P_{\Omega _L}(k)(\omega-ct)e^{-\frac{1}{c}\int_{\omega\!-\!ct}^{\omega}P_{\Omega _L}(k)(y)\mathrm{d}y} P_{\Omega _L}(x^0)(\omega-ct)\!+\!ce^{-\frac{1}{c}\int_{\omega\!-\!ct}^{\omega}P_{\Omega _L}(k)(y)\mathrm{d}y} P_{\Omega _L}\left(\frac{\mathrm{d}}{\mathrm{d} \omega}x^{0}\right)(\omega\!-\!ct).
    \end{flalign*}
    This shows~\eqref{TransportEquationwithStateFeedback}.
    Checking the initial and boundary conditions is straightforward. The solution formula can be extended to $x^0 \in L^2([0,L])$ via a standard density argument (see~\cite[Proposition II.1.5.]{Werner2006}).
\end{proof}

\paragraph{Proof of~\Cref{Lem: Necessary}}
\begin{proof}
    \textbf{(i)} Assume $a_1 >0$. Choose $0 < L_1 < a_1$. Then the transport equation is uncontrolled on $\Omega _{L_1}=[0,L_1]$. Since the semigroup corresponding to the uncontrolled transport equation with domain $\Omega _{L_1}$ is not exponentially stable (see~\eqref{Eq: NormPreservation}), this is a contradiction to domain-uniform stabilizability.\\
    \textbf{(ii)} Assume $|\Omega _c|=\mu _c <\infty$. Let $k_B$ be an arbitrary state feedback as in~\eqref{GeneralStateFeedback} with $\hat{k}_B:=\n{k_B}_{\infty}$. For $t_L = \frac{L}{c}$ we find
    \begin{equation*}
        \n{T^\varphi _{\Omega _L}(t_L)x^0}_{L^2([0,L]}^2=\int _0^L \n{e^{-\frac{1}{c}\int _{\omega-L}^\omega P_{\Omega _L}(k_B|_{\Omega _c^L})(y)\mathrm{d}y}P_{\Omega _L}(x^0)(\omega-ct)}^2 d\omega \geq e^{-2\frac{\mu _c}{c}\hat{k}_B}\n{x^0}_{L^2([0,L])}.
    \end{equation*}
    This shows, that the controlled transport equation is not domain-uniformly stabilizable with this norm.\\
\end{proof}

\paragraph{Proof of~\Cref{Lem: IntervalCondition}(i)}
\begin{proof}
    For $\omega =0$ the inequality in~\eqref{ExponentialInequality} is equivalent to
    \begin{flalign*}
        \exists M,k>0 \, \forall t \geq 0: e^{kt - \int _0^t e(\tau )\mathrm{d}\tau} \leq M.
    \end{flalign*}
    Since the exponential function is continuous and monotonously increasing this condition is equi\-valent to
    \begin{equation}
    \label{ExponentialInequalitySimplified}
        \exists \tilde{M},k>0 \, \forall t \geq 0: kt - \int _0^t e(\tau )\mathrm{d}\tau \leq \tilde{M}.
    \end{equation}
    Now note that the left hand side of this inequality can be expressed in the form
    \begin{equation}
    \label{SimplifiedRepresentation}
        kt - \int _0^t e(\tau )\mathrm{d}\tau = \left\{ \begin{array}{cc}
            kt - \sum _{j=1}^{n-1} k_jB_j - k_n(t-a_n), & t \in [a_n,b_n] \\
            kt - \sum _{j=1}^{n-1} k_jB_j, & t \in [b_{n-1}, a_n]
        \end{array}\right. .
    \end{equation}
    For the rest of the proof it will be our aim to show, that the conditions~\eqref{ExponentialInequalitySimplified} and~\eqref{IntervalCondition1} are equivalent. For necessity assume that for any constants $\tilde{M},k>0$ there exists a number $n(\tilde{M},k) \in \N$ such that
    \begin{equation*}
        ka_{n(\tilde{M},k)} - \sum _{j=1}^{{n(\tilde{M},k)}-1} k_j (b_j-a_j) = ka_{n(\tilde{M},k)} - \sum _{j=1}^{{n(\tilde{M},k)}-1} k_j B_j > \tilde{M}.
    \end{equation*}
    Now choose $t(\tilde{M},k):=a_{n(\tilde{M},k)}$. We find
    \begin{equation*}
        kt(\tilde{M},k) - \int _0^{t(\tilde{M},k)} e(\tau )\mathrm{d}\tau = ka_{n(\tilde{M},k)} - \sum _{j=1}^{n(\tilde{M},k)-1} k_jB_j > \tilde{M}.
    \end{equation*}
    Therefore condition~\eqref{ExponentialInequalitySimplified} is not fulfilled which shows necessity. \\
    \noindent For sufficiency we assume that~\eqref{IntervalCondition1} is fulfilled and perform a case distinction to show~\eqref{ExponentialInequalitySimplified}.\newline
    \textit{Case 1} ($t \in [b_{n-1}, a_n]$): In this case we have the inequality
    \begin{equation*}
        kt - \int _0^t e(\tau )\mathrm{d}\tau \overset{\eqref{SimplifiedRepresentation}}{=} kt - \sum _{j=1}^{n-1} k_jB_j \leq ka_n - \sum _{j=1}^{n-1} k_jB_j \overset{\eqref{IntervalCondition}}{\leq} \tilde{M}.
    \end{equation*}
    \textit{Case 2} ($t \in [a_n,b_n]$ and $k\leq k_n$): In this case we have the inequality
    \begin{flalign*}
        kt - \int _0^t e(\tau )\mathrm{d}\tau &\overset{\eqref{SimplifiedRepresentation}}{=} kt - \sum _{j=1}^{n-1} k_jB_j - k_n(t-a_n) = (k-k_n)t + k_n a_n - \sum _{j=1}^{n-1} k_jB_j\\
        &\leq (k-k_n)a_n + k_n a_n - \sum _{j=1}^{n-1} k_jB_j = ka_n - \sum _{j=1}^{n-1} k_jB_j \overset{\eqref{IntegralTerm}}{\leq} \tilde{M}.
    \end{flalign*}
    \textit{Case 3} ($t \in [a_n,b_n]$ and $k\geq k_n$): In this case we have the inequality
    \begin{flalign*}
        kt - \int _0^t e(\tau )\mathrm{d}\tau &\overset{\eqref{SimplifiedRepresentation}}{=} kt - \sum _{j=1}^{n-1} k_jB_j - k_n(t-a_n) = (k-k_n)t + k_n a_n - \sum _{j=1}^{n-1} k_jB_j\\
        &\leq (k-k_n)b_n + k_n a_n - \sum _{j=1}^{n-1} k_jB_j = kb_n - \sum _{j=1}^{n} k_jB_j \leq ka_{n+1} - \sum _{j=1}^{n} k_jB_j \overset{\eqref{IntegralTerm}}{\leq} \tilde{M}.
    \end{flalign*}
    Therefore~\eqref{ExponentialInequalitySimplified} is fulfilled which concludes the proof.
\end{proof}

\paragraph{Proof of~\Cref{Lem: IntervalCondition}(ii)}
\begin{proof}
    \textbf{$\mathbf{\implies}$:} Suppose that~\eqref{IntervalCondition} does not hold. Then for any $\tilde{M}>0$ there exist integers $n,m\in \N$, $n \geq m$ such that
    \begin{equation*}
        k(a_n-b_m) - \sum_{j=m+1}^{n-1} k_j B_j> \tilde{M}.
    \end{equation*}
    Define $\omega_0 := b_m$ and $t_0 := a_n-b_m$. Then the inequality
    \begin{equation*}
        kt- \int _{\omega_0}^{\omega_0+t_0} e(\tau) \,\mathrm{d}\tau = k (a_n-b_m)- \int _{b_m}^{a_n}e(\tau) \,\mathrm{d}\tau = k(a_n-b_m) - \sum_{j=m+1}^{n-1} k_j B_j > \tilde{M}
    \end{equation*}
    follows. Plugging both sides of the inequality into the exponential function yields
    \begin{equation*}
        e^{- \int _{\omega_0}^{\omega_0+t_0} e(\tau) \,\mathrm{d}\tau} > e^{\tilde{M}} e^{-kt_0}
    \end{equation*}
    due to monotony of the exponential function. Since this function is also unbounded we find, that~\eqref{ExponentialInequality} is not fulfilled.\\
    \textbf{$\mathbf{\impliedby}$:} In order to show this direction we assume that~\eqref{IntervalCondition} is fulfilled. Let $\omega>0$ be arbitrary but fixed. Define $e_\omega: \R _{>0} \rightarrow \R , \, e_\omega(t):=e(t+\omega)$.
    By substitution $r:=\tau - \omega$ we find
    \begin{equation*}
        \int _\omega^{\omega+t} e(\tau ) \mathrm{d}\tau = \int _0^{t} e(r + \omega) \mathrm{d}r = \int _0^{t} e_\omega(r) \mathrm{d}r.
    \end{equation*}
    For $\omega \in [b_m, b_{m+1})$, $m\in \N \cup \{0\}$ define $\tilde{b}_0 := 0$, $\tilde{a}_1:= \mathrm{max} \{0, a_{m+1}-\omega\}$, $\tilde{b}_1:= b_{m+1}-\omega$ and
    \begin{equation*}
        \forall n \geq 2: \tilde{a}_n:=a_{n+m}-\omega, \platz \tilde{b}_n:=b_{n+m}-\omega.
    \end{equation*}
    Furthermore define for $n\in \N$ the shifted sequence $(\tilde{k}_n)_{n\in \N}$ with $\tilde{k}_n:=k_{n+m}$. Note that for all $n\geq 2$ the equality
    \begin{equation*}
        \tilde{B}_n = \tilde{b}_n - \tilde{a}_n =  (b_{n+m} - \omega) - (a_{n+m} - \omega) = b_{n+m}-a_{n+m} = B_{n+m}
    \end{equation*}
    holds.
    Our aim will be to show, that condition~\eqref{IntervalCondition} holds for the shifted function $e_\omega$, the shifted sequence $(\tilde{k}_n)_{n\in \N}$ and the shifted intervals $\tilde{I}_n:=[\tilde{a}_n,\tilde{b}_n]$. For $n=1$ we have
    \begin{equation*}
        k \tilde {a_1} - \sum _{j=1}^{0} \tilde{k}_j \tilde{B}_j = k \tilde {a_1} = k \, \mathrm{max} \{0, a_{m+1}-\omega\} \leq k \, (a_{m+1}-b_m) \overset{\eqref{IntervalCondition}}{\leq} \tilde{M}.
    \end{equation*}
    For $n\geq 2$ we find
    \begin{flalign*}
        k \, \tilde {a_n} - \sum _{j=1}^{n-1} \tilde{k}_j \tilde{B}_j
        = \, &k \, (a_{n+m}-\omega) - \sum _{j=2}^{n-1} k_{j+m} B_{j+m} - k_{m+1}(b_{m+1} - \omega -\mathrm{max} \{0, a_{m+1}-\omega\})\\
        =\, &k \, (a_{n+m}-\omega) - \sum _{j=m+2}^{n+m-1} k_{j} B_{j} - k_{m+1}(b_{m+1} -\omega -\mathrm{max} \{0, a_{m+1}-\omega\})\\
        \leq \, & \left\{\begin{array}{cc}
            k \, (a_{n+m}-b_m) - \sum _{j=m+2}^{n+m-1} k_{j} B_{j} - k_{m+1}(b_{m+1}- a_{m+1}), &  \omega \leq a_{m+1}\\
            k \, (a_{n+m}-b_{m+1}) - \sum _{j=m+2}^{n+m-1} k_{j} B_{j}, & \omega > a_{m+1}
        \end{array}\right.\\
        = \, &\left\{\begin{array}{cc}
            k \, (a_{n+m}-b_m) - \sum _{j=m+1}^{n+m-1} k_{j} B_{j}, &  \omega \leq a_{m+1}\\
            k \, (a_{n+m}-b_{m+1}) - \sum _{j=m+2}^{n+m-1} k_{j} B_{j}, & \omega > a_{m+1}
        \end{array}\right.\overset{~\eqref{IntervalCondition}}{\leq} \tilde{M}.
    \end{flalign*}
    From (i) it now follows, that
    \begin{equation*}
        \forall t \geq 0: e^{-\int _\omega^{\omega+t} e(\tau ) \mathrm{d}\tau} = e^{-\int _0^{t} e_\omega(r) \mathrm{d}r} \overset{(i)}{\leq} \tilde{M}e^{-kt}.
    \end{equation*}
\end{proof}

\paragraph{Proof of~\Cref{Lem: IntervalCondition}(iii)}
\begin{proof}
     \textbf{$\mathbf{\implies}$:} Assume, that~\eqref{IntervalCondition} does not hold. According to (ii) this implies, that for any constant $M>0$ there exists $\omega_0 \in [0,L], t_0\geq 0$ such that
     \begin{equation*}
         e^{-\int _{\omega_0}^{\omega_0+t_0}e(\tau )\mathrm{d}\tau} > Me^{-kt_0} \implies e^{-\frac{1}{c}\int _{\omega_0}^{\omega_0+t_0}e(\tau )\mathrm{d}\tau} > M^{\frac{1}{c}}e^{-kt_0}.
     \end{equation*}
     Choosing $L=2(\omega_0+t_0)$ it follows that
     \begin{equation*}
          e^{-\frac{1}{c}\int _{\omega_0}^{\omega_0+c\frac{t_0}{c}}P_{\Omega _L}(e_L)(y )\mathrm{d}y} = e^{-\frac{1}{c}\int _{\omega_0}^{\omega_0+t_0}e(y)\mathrm{d}y} > M^{\frac{1}{c}}e^{-kt_0}.
     \end{equation*}
     Therefore~\eqref{ExponentialInequality} is also not fulfilled.\\
     \textbf{$\mathbf{\impliedby}$:} Assume that~\eqref{IntervalCondition} is fulfilled. Let $L>0$ $\omega\in [0,L]$ and $t\geq 0$ be arbitrary. Using Euclidean division we can find $m_E\in \N$ and $t_0 \in \left[0,\frac{L}{c}\right)$ such that $t=m\frac{L}{c}+t_0$. This allows us to rewrite the integral in~\eqref{ExponentialInequality}. By using the periodicity of $P_{\Omega _L}(e_L)$ we find the equalities
     \begin{flalign*}
         \int _{\omega-ct}^\omega P_{\Omega _L}(e_L)(y)\mathrm{d}y &= \int _{\omega-m_EL-ct_0}^\omega P_{\Omega _L}(e_L)(y)\mathrm{d}y\\
         &= \int _{\omega-m_EL-ct_0}^{\omega-m_EL} P_{\Omega _L}(e_L)(y)\mathrm{d}y + \int _{\omega-m_EL}^{\omega} P_{\Omega _L}(e_L)(y)\mathrm{d}y\\
         &= \int _{\omega-ct_0}^{\omega} P_{\Omega _L}(e_L)(y)\mathrm{d}y + m_E \int _{0}^{L} e_L(y)\mathrm{d}y.
     \end{flalign*}
     We consider two different cases.\\
     \textbf{Case 1 ($\omega-ct_0<0$):} In this case we have
     \begin{flalign*}
         \int _{\omega-ct}^\omega P_{\Omega _L}(e_L)(y)\mathrm{d}y &= \int _{\omega-ct_0}^{\omega} P_{\Omega _L}(e_L)(y)\mathrm{d}y + m_E\int _{0}^{L} e_L(y)\mathrm{d}y\\
         &= \int _{0}^{\omega} e_L(y)\mathrm{d}y + \int _{\omega-ct_0+L}^{L} e_L(y)\mathrm{d}y + m_E \int _{0}^{L} e_L(y)\mathrm{d}y\\
         &= \int _{0}^{\omega} e(y)\mathrm{d}y + \int _{\omega-ct_0+L}^{L} e(y)\mathrm{d}y + m_E \int _{0}^{L} e(y)\mathrm{d}y.
     \end{flalign*}
     Using (ii) we find, that there exists a constant $\hat{M}>0$ which is independent of $\omega$, $t$ and $L$ such that
     \begin{flalign*}
         e^{-\int _{\omega-ct}^\omega P_{\Omega _L}(e_L)(y)\mathrm{d}y} \!=\! e^{-\left(\int _{0}^{\omega} e(y)\mathrm{d}y + \int _{\omega-ct_0+L}^{L} e(y)\mathrm{d}y + m_E \int _{0}^{L} e(y)\mathrm{d}y\right)}\!\leq \! \hat{M} e^{-k(m_EL+\omega-\omega+ct_0)}\! =\! \hat{M} e^{-k(m_EL+ct_0)}
     \end{flalign*}
     \textbf{Case 2 ($\omega-ct_0\geq 0$):} In this case we have
     \begin{flalign*}
         \int _{\omega-ct}^\omega P_{\Omega _L}(e_L)(y)\mathrm{d}y &= \int _{\omega-ct_0}^{\omega} P_{\Omega _L}(e_L)(y)\mathrm{d}y + m_E \int _{0}^{L} e_L(y)\mathrm{d}y\\
         &= \int _{\omega-ct_0}^{\omega} e_L(y)\mathrm{d}y  + m_E \int _{0}^{L} e_L(y)\mathrm{d}y= \int _{\omega-ct_0}^{\omega} e(y)\mathrm{d}y + m_E \int _{0}^{L} e(y)\mathrm{d}y.
     \end{flalign*}
     Using (ii) we find, that there exists a constant $\hat{M}>0$ which is independent of $\omega$, $t$ and $L$ such that
     \begin{flalign*}
         e^{-\int _{\omega-ct}^\omega P_{\Omega _L}(e_L)(y)\mathrm{d}y} &= e^{-\left(\int _{\omega-ct_0}^{\omega} e(y)\mathrm{d}y + m_E \int _{0}^{L} e(y)\mathrm{d}y\right)}
         \leq \hat{M} e^{-k(m_EL+ct_0)}
     \end{flalign*}
     For both cases we can conclude
     \begin{equation*}
         e^{-\frac{1}{c}\int _{\omega-ct}^\omega P_{\Omega _L}(e_L)(y)\mathrm{d}y} \leq \hat{M}^{\frac{1}{c}}e^{-\frac{k}{c}(m_EL+ct_0)} = \hat{M}^{\frac{1}{c}}e^{-kt}.
     \end{equation*}
     Choosing $M=\hat{M}^{\frac{1}{c}}$ we find that~\eqref{ExpInequalityPeriodified} is fulfilled.
\end{proof}

\newpage

\section{Proof of~\Cref{Lem: SolNonConstantTransport}}
\label{App: LemmataNonConstant}

\paragraph{Proof of~\Cref{Lem: SolNonConstantTransport}(i)}
\begin{proof}
We first prove uniqueness of a solution to \eqref{ODETransportForward}. Note that for initial values $p_0$ at which $c$ is Lipschitz-continuous, local uniqueness of the solution follows from the well known Theorem of Picard-Lindelöf \cite[Theorem 8.14]{Amann.2006}. Therefore w.l.o.g.\ let $p_0$ be a point at which $c$ is not Lipschitz-continuous. Assume that there are two solutions $p_1(t)$ and $p_2(t)$ such that
\begin{equation*}
    \exists \, t_1 >0 \, \forall \, t \in (0,t_1): p_1(t) \neq p_2(t).
\end{equation*}
Since $\forall x \in \R: c(x) > c_\mathrm{min} >0$ we find $\forall \, t \in (0,t_1): \min \{p_1(t), p_2(t)\} >p_0$. Let $r>p_0$ be such that $c$ is Lipschitz-continuous on the interval $(p_0,r]$ with Lipschitz constant $L\geq 0$. Since $c(x) < c_\mathrm{max} < \infty$ we find $0<t_2 \leq t_1$ such that $\forall \, t \in (0,t_2): p_1(t),p_2(t) \in (p_0,r]$. 
Thus, for all $t \in (0,t_2]$,
\begin{flalign*}
    \n{p_1(t)\!-\!p_2(t)}\!=\!\!\! \int _0^t \n{c(p_1(\tau))\!-\!c(p_2(\tau))} d\tau
    \!=\!\!\! \underset{(0,t]}{\int} \n{c(p_1(\tau))\!-\!c(p_2(\tau))} d\tau\!\leq\! \int _0^t L\n{p_1(\tau)\!-\!p_2(\tau)} d\tau.
\end{flalign*}
Now Gronwall's Lemma~\cite[Lemma 2.4]{Logemann.2014} yields $p_1(t) = p_2(t)$ for all $t \in (0,t_2)$ and therefore a contradiction. An analogous argumentation yields the result for \eqref{ODETransportBackward}.
\end{proof}

\paragraph{Proof of~\Cref{Lem: SolNonConstantTransport}(ii)}
\begin{proof}
For arbitrary $t\in [0,T]$ separation of variables yields
\begin{flalign}
\label{SVCalc}
\begin{split}
     &-\int _{q_0}^{q(t,q_0)} \frac{1}{P_{\Omega _L}(c)(y)}\mathrm{d}y = \int _0^t 1\mathrm{d}\tau = t = \int _{T-t}^T 1\mathrm{d}\tau\\
     = &\int _{p(T-t,p_0)}^{p(T,p_0)} \frac{1}{P_{\Omega _L}(c)(y)}\mathrm{d}y = \int _{p(T-t,p_0)}^{q_0} \frac{1}{P_{\Omega _L}(c)(y)}\mathrm{d}y = -\int _{q_0}^{p(T-t,p_0)} \frac{1}{P_{\Omega _L}(c)(y)}\mathrm{d}y.
     \end{split}
\end{flalign}
Due to regularity of $c$ the solution $q(t,q_0)$ of~\eqref{ODETransportBackward} is unique. Using separation of variables we find that $q(t,q_0)$ is the unique solution of solving the equality
\begin{equation*}
    -\int _{q_0}^{x_q} \frac{1}{P_{\Omega _L}(c)(y)}\mathrm{d}y = \int _0^t 1\mathrm{d}\tau
\end{equation*}
for $x_q$. From~\eqref{SVCalc} we know, that $p(T-t,p_0)$ also solves this equality. Uniqueness of the solution implies
\begin{equation*}
    q(t,q_0)=p(T-t,p_0)
\end{equation*}
and therefore the claim.
\end{proof}

\paragraph{Proof of~\Cref{Lem: SolNonConstantTransport}(iii)}
\begin{proof}
    Let $t>0$ be arbitrary but fixed. By definition of the derivative we have
    \begin{equation*}
        \frac{\partial}{\partial q_0}q(t,q_0)= \underset{h\uparrow 0}{\lim}\frac{q(t,q_0+h)-q(t,q_0)}{h}.
    \end{equation*}
    Due to positivity of $c$ the solution $q(t,q_0)$ is strictly monotonously falling. Therefore $q(t,q_0)<q_0$. Let $h \in (q(t,q_0)-q_0,0)$. Since $q(t,q_0)$ is continuous in $t$ there exists $T_h \in (0,t)$ such that $q(T_h,q_0)=q_0+h$. Using separation of variables we find
    \begin{equation*}
        T_h = \int _0^{T_h} 1\mathrm{d}\tau = -\int _{q_0}^{q_0+h} \frac{1}{P_{\Omega _L}(c)(y)}\mathrm{d}y.
    \end{equation*}
    The cocycle property of autonomous differential equations implies
    \begin{equation}
    \label{CocycleEq}
        q(t,q_0+h)=q(t,q(T_h,q_0))=q(t+T_h,q_0) = q\left(t-\int _{q_0}^{q_0+h} \frac{1}{P_{\Omega _L}(c)(y)}\mathrm{d}y,q_0\right).
    \end{equation}
Let $(h_k)_{k\in \N} \in (q(t,q_0)-q_0,0)^\N$ be a sequence with $\underset{k \rightarrow \infty}{\lim}h_k=0$. For the left-side limit value of the differential quotient we find the equalities
\begin{flalign*}
    \frac{\partial}{\partial q_0}q(t,q_0)&= \underset{k\rightarrow \infty}{\lim}\frac{q(t,q_0+h_k)-q(t,q_0)}{h_k} \overset{\eqref{CocycleEq}}{=} \underset{k\rightarrow \infty}{\lim}\frac{q(t-\int_{q_0}^{q_0+h_k} \frac{1}{P_{\Omega _L}(c)(y)}\mathrm{d}y,q_0)-q(t,q_0)}{h_k}\\
    &= \frac{\partial}{\partial p}q\left(t-\int_{q_0}^{p} \frac{1}{P_{\Omega _L}(c)(y)}\mathrm{d}y,q_0\right)|_{p=q_0} = -\frac{1}{c(p)}\dot{q}\left(t-\int_{q_0}^{p} \frac{1}{P_{\Omega _L}(c)(y)}\mathrm{d}y,q_0\right)|_{p=q_0}\\
    &= -\frac{1}{P_{\Omega _L}(c)(q_0)}\dot{q}\left(t,q_0\right).
\end{flalign*}
The case of the right-hand limit value can be treated analogously.
\end{proof}

\paragraph{Proof of~\Cref{Lem: SolNonConstantTransport}(iv)}
\begin{proof}
    For $x^0\in \mathrm{dom}(A_L)$ we find
    \begin{flalign*}
        \frac{\partial}{\partial \omega} P_{\Omega _L}(x^0)(q(t,\omega)) &= \left(\frac{\partial}{\partial \omega}q(t,\omega)\right) P_{\Omega _L}\left(\frac{\mathrm{d}}{\mathrm{d} \omega}x^{0}\right)(q(t,\omega))\\
        &= -\frac{1}{P_{\Omega _L}(c)(\omega)}\dot{q}\left(t,\omega\right) P_{\Omega _L}\left(\frac{\mathrm{d}}{\mathrm{d} \omega}x^{0}\right)(q(t,\omega))
    \end{flalign*}
    and
    \begin{equation*}
        \frac{\partial}{\partial t} P_{\Omega _L}(x^0)(q(t,\omega)) = \dot{q}(t,\omega)P_{\Omega _L}\left(\frac{\mathrm{d}}{\mathrm{d} \omega}x^{0}\right)(q(t,\omega))
    \end{equation*}
    using the chain rule and (ii). Therefore we have
    \begin{equation*}
        \frac{\partial}{\partial t} x(\omega ,t) = -P_{\Omega _L}(c)(\omega ) \frac{\partial}{\partial \omega} x(\omega ,t )
    \end{equation*}
    and
    \begin{equation*}
        x(\omega ,0) = P_{\Omega _L}(x^0)(q(0,\omega)) = P_{\Omega _L}(x^0)(\omega) = x^0(\omega)
    \end{equation*}
    for almost all $\omega \in [0,L]$.
    Finally we use separation of variables to find
    \begin{equation*}
        \int _0^t 1 \mathrm{d}t = -\int _0^{q(t,0)} \frac{1}{P_{\Omega _L}(c)(y)}\mathrm{d}y = -\int _L^{q(t,L)} \frac{1}{P_{\Omega _L}(c)(y)}\mathrm{d}y = -\int _0^{q(t,L)-L} \frac{1}{P_{\Omega _L}(c)(y)}\mathrm{d}y
    \end{equation*}
    and therefore via uniqueness of the solution $q(t,L)-L=q(t,0)$. This implies
    \begin{equation*}
        x(0,t)=P_{\Omega _L}(x^0)(q(t,0))=P_{\Omega _L}(x^0)(q(t,L)-L)=P_{\Omega _L}(x^0)(q(t,L))=x(L,t).
    \end{equation*}
    Again the solution formula can be extended to $x^0 \in L^2([0,L])$ via (see~\cite[Proposition II.1.5.]{Werner2006}).
\end{proof}

\paragraph{Proof of~\Cref{Lem: SolNonConstantTransport}(v)}
\begin{proof}
    For $x^0 \in \mathrm{dom}(A_L)$ using the fundamental theorem of calculus, the chain rule of differentiation and (ii) leads to
    \begin{flalign*}
        \frac{\partial}{\partial \omega} x(\omega, t)
        = &-\left( \frac{P_{\Omega _L}(k)(\omega)}{P_{\Omega _L}(c)(\omega )} - \frac{\partial}{\partial \omega}q (t,\omega ) \frac{P_{\Omega _L}(k)(q (t,\omega))}{P_{\Omega _L}(c)(q (t,\omega))}\right)e^{-\int_{q(t,\omega)}^{\omega}\frac{P_{\Omega _L}(k)(y)}{P_{\Omega _L}(c)(y)}\mathrm{d}y} P_{\Omega _L}(x^0)(q(t,\omega))\\
        &+\frac{\partial}{\partial \omega}q (t,\omega ) e^{-\int_{q(t,\omega)}^{\omega}\frac{P_{\Omega _L}(k)(y)}{P_{\Omega _L}(c)(y)}\mathrm{d}y} P_{\Omega _L}\left(\frac{\mathrm{d}}{\mathrm{d} \omega}x^{0}\right)(q(t,\omega))\\
        =& - \frac{(kx)(\omega, t)}{c(\omega )} - \frac{\dot{q}(t,\omega)}{c(\omega)} e^{-\int_{q(t,\omega)}^{\omega}\frac{P_{\Omega _L}(k)(y)}{P_{\Omega _L}(c)(y)}\mathrm{d}y} \left(\frac{P_{\Omega _L}(k)(q (t,\omega))}{P_{\Omega _L}(c)(q (t,\omega))} P_{\Omega _L}(x^0)(q(t,\omega))\right)\\
        &-\frac{\dot{q}(t,\omega)}{c(\omega)} e^{-\int_{q(t,\omega)}^{\omega}\frac{P_{\Omega _L}(k)(y)}{P_{\Omega _L}(c)(y)}\mathrm{d}y}\left(P_{\Omega _L}\frac{\mathrm{d}}{\mathrm{d} \omega}x^{0}\right)(q(t,\omega))
    \end{flalign*}
    and
    \begin{flalign*}
        \frac{\partial}{\partial t} x(\omega, t) 
        \!=&\dot{q}(t,\omega) e^{-\int_{q(t,\omega)}^{\omega}\frac{P_{\Omega _L}(k)(y)}{P_{\Omega _L}(c)(y)}\mathrm{d}y}\! \left(\frac{P_{\Omega _L}(k)(q (t,\omega))}{P_{\Omega _L}(c)(q (t,\omega))} P_{\Omega _L}(x^0)(q(t,\omega)) \!+ \!P_{\Omega _L}\left(\frac{\mathrm{d}}{\mathrm{d} \omega}x^{0}\right)(q(t,\omega))\right).
    \end{flalign*}
    This shows~\eqref{TransportEquationwithStateFeedback}. The boundary and initial conditions as well as the extension to $x^0 \in L^2([0,L])$ can be shown in analogy to the proof of~\Cref{Lem: SolutionFormulaStateFeedback1}.
\end{proof}

\section{Proof of~\Cref{Prop: Equivalence}}
\label{App: WaveEq}
\vspace{-0.5cm}
\begin{proof}
For readability we often leave out the arguments $\omega$ and/or $t$ in the following.\\
    \textbf{$\Rightarrow$}: Assume that $x\in C^1([0,T],H^2(\Omega _L)\cap H_0^1(\Omega _L))$ is a classical solution of~\eqref{eq: ClosedLoop}. Define $v \in C^1([0,T],H^1(\Omega _L))$ as the solution of
    \begin{equation}
    \label{WaveODE}
        \frac{\partial}{\partial t}v=-k\,\chi _{\Omega_L^c}v-c^2\frac{\partial}{\partial \omega}x, \quad  v(\omega,0) = -\int _0^\omega x_{\Omega _L}^1(s)+k\,\chi _{\Omega_L^c}(s) x_{\Omega _L}^0(s)\mathrm{d}s.
    \end{equation}
    Note that for each $\omega \in \Omega _L$~\eqref{WaveODE} is a linear ordinary differential equation with stabilizing part $-kv$. Since $x\in C^1([0,T],H^2(\Omega _L)\cap H_0^1(\Omega _L))$ its spatial derivative is continuous and therefore bounded on the compact set $[0,T]\times \overline{\Omega _L}$. Therefore, the solution $v$ of~\eqref{WaveODE} is well- defined on $[0,T]\times \overline{\Omega _L}$. Rearranging the state equation in~\eqref{eq: ClosedLoop} we find
    \begin{equation*}
    c^2\frac{\partial ^2}{\partial \omega ^2}x(\omega,t)=\frac{\partial^2}{\partial t^2}{x}(\omega,t)+ \chi _{\Omega_L^c}(\omega) \left(2k \frac{\partial}{\partial t}x(\omega,t) +k^2x(\omega,t)\right)
    \end{equation*}
    Taking the time integral over this equation and using the fundamental theorem of calculus (FTC) we find
    \begin{flalign}
    \label{ProofStep3}
    \begin{split}
        c^2\int _0^t \frac{\partial^2}{\partial \omega^2}x(\tau)\mathrm{d}\tau &\overset{\eqref{eq: ClosedLoop}}{=}\int _0^t \frac{\partial^2}{\partial \tau^2}x(\tau) + 2k \,\chi _{\Omega_L^c}\frac{\partial}{\partial \tau} x(\tau) + k^2 \,\chi _{\Omega_L^c} x(\tau)\mathrm{d}\tau\\
        &\overset{\mathrm{FTC}}{=} \frac{\partial}{\partial t}x(t)-\frac{\partial}{\partial t}x(0)+2k\,\chi _{\Omega_L^c}(x(t)-x(0))+\int _0^t k^2 \,\chi _{\Omega_L^c} x(\tau)\mathrm{d}\tau,
        \end{split}
    \end{flalign}
    where we omitted the spatial argument $\omega$ for brevity.
    Integrating the differential equation~\eqref{WaveODE} over time yields
    \begin{equation}
    \label{ProofStep3.5}
        \int _0^t \frac{\partial}{\partial \tau}v(\tau)\mathrm{d}\tau  \overset{\mathrm{FTC}}{=} v(t)-v(0) \overset{\eqref{WaveODE}}{=} -k \,\chi _{\Omega_L^c}\int _0^t v(\tau)\mathrm{d}\tau-c^2\int _0^t \frac{\partial}{\partial \omega} x(\tau) \mathrm{d}\tau.
    \end{equation}
    Taking the spatial derivative of this equality and using the initial condition from~\eqref{WaveODE} we find
    \begin{flalign}
    \label{ProofStep4}
    \begin{split}
        \frac{\partial}{\partial \omega}v(t) - \frac{\partial}{\partial \omega}v(0) &\overset{\eqref{WaveODE}}{=} \frac{\partial}{\partial \omega}v(t)+\frac{\partial}{\partial t}x(0)+k\,\chi _{\Omega_L^c}x(0)\\
        &\overset{\eqref{ProofStep3.5}}{=}  -k\,\chi _{\Omega_L^c}\int _0^t\frac{\partial}{\partial \omega}v(\tau)\mathrm{d}\tau-c^2\int _0^t \frac{\partial^2}{\partial \omega ^2}x(\tau)\mathrm{d}\tau
        \end{split}
    \end{flalign}
    almost everywhere since $\chi _{\Omega_L^c}$ is piecewise constant.
    Rearranging the terms in~\eqref{ProofStep4} leads to
    \begin{equation}
        \label{ProofStep4.5}
        c^2\int _0^t \frac{\partial^2}{\partial \omega ^2}x(\tau)\mathrm{d}\tau = -\frac{\partial}{\partial \omega}v(t)-\frac{\partial}{\partial t}x(0)-k\,\chi _{\Omega_L^c}x(0)-k\,\chi _{\Omega_L^c}\int _0^t\frac{\partial}{\partial \omega}v(\tau)\mathrm{d}\tau.
    \end{equation}
    In equations~\eqref{ProofStep4.5} and~\eqref{ProofStep2} the same term incorporating the second spatial derivative of $x$ appears on the right hand side. Therefore, the left hand sides are equal as well. This observation implies
    \begin{flalign*}
        &\frac{\partial}{\partial t}x(t)-\frac{\partial}{\partial t}x(0)+2k\,\chi _{\Omega_L^c}(x(t)-x(0))+\int _0^t k^2 \,\chi _{\Omega_L^c} x(\tau)\mathrm{d}\tau\\
        = &-\frac{\partial}{\partial \omega}v(t)-\frac{\partial}{\partial t}x(0)-k\,\chi _{\Omega_L^c}x(0)-k\,\chi _{\Omega_L^c}\int _0^t\frac{\partial}{\partial \omega}v(\tau)\mathrm{d}\tau.
    \end{flalign*}
    By removing the terms which appear on both sides of this equation and using the fundamental theorem of calculus again we find
    \begin{flalign*}
         &\frac{\partial}{\partial t}x(t)+k\,\chi _{\Omega_L^c}x(t)+k\chi _{\Omega_L^c}(x(t)-x(0))+k^2\int _0^t x(\tau)\mathrm{d}\tau\\
         \overset{\mathrm{FTC}}{=}& \frac{\partial}{\partial t}x(t)+k\,\chi _{\Omega_L^c}x(t)+k \,\chi _{\Omega_L^c}\int _0^t \frac{\partial}{\partial \tau} x(\tau)+ k\,\chi _{\Omega_L^c}x(\tau)\mathrm{d}\tau\\ 
         =& -\frac{\partial}{\partial \omega}v(t)-k\,\chi _{\Omega_L^c}\int _0^t \frac{\partial}{\partial \omega}v(\tau)\mathrm{d}\tau.
    \end{flalign*}
    This implies for $h(t):=\frac{\partial}{\partial t}x(t)+k\,\chi _{\Omega_L^c}x(t)+\frac{\partial}{\partial \omega}v(t)$ that
    \begin{equation*}
        h(t)+k\,\chi _{\Omega_L^c}\int _0^t h(\tau)\mathrm{d}\tau=0 \implies \frac{\partial}{\partial t} h(t) = -k\,\chi _{\Omega_L^c}h(t)
    \end{equation*}
    Together with the initial condition 
        $h(0)=\frac{\partial}{\partial t}x(0)+xu(0)+\frac{\partial}{\partial \omega}v(0)=0$
    we find
    \begin{equation*}
        h \equiv 0 \implies \frac{\partial}{\partial t}x(t)=-\frac{\partial}{\partial \omega}v(t)-k\,\chi _{\Omega_L^c}x(t).
    \end{equation*}
    Therefore $(x,v)$ is a classical solution of~\eqref{eq: ClosedLoop}.
    \\
    \textbf{$\Leftarrow$}: Assume that $(\xi _1, \xi _2) \in C^2([0,T],H^2(\Omega _L)\cap H_0^1(\Omega _L))\times C^1([0,T],H^1(\Omega _L))$ is a classical solution of~\eqref{WaveEquation1}. Taking the time derivative of the first and the spatial derivative of the second equation in~\eqref{WaveEquation1} we find
    \begin{equation}
    \label{ProofStep1}
        \frac{\partial^2}{\partial t^2}{\xi _1}=\frac{\partial^2}{\partial t \partial \omega}{\xi _2} - k\,\chi _{\Omega_L^c}\frac{\partial}{\partial t}{\xi _1} \quad \textrm{and} \quad c^2\frac{\partial^2}{\partial \omega^2}{\xi _1}=\frac{\partial^2}{\partial t \partial \omega}{\xi _2} - k\,\chi _{\Omega_L^c}\frac{\partial}{\partial \omega}{\xi _2}
    \end{equation}
    Subtracting the second from the first equation in~\eqref{ProofStep1} leads to
    \begin{equation}
        \label{ProofStep2}
        \frac{\partial^2}{\partial t^2}{\xi _1}-c^2\frac{\partial^2}{\partial \omega^2}{\xi _1} = -k\,\chi _{\Omega_L^c}\left(\frac{\partial}{\partial t}{\xi _1} + \frac{\partial}{\partial \omega}{\xi _2}\right) \overset{\eqref{WaveEquation1}}{=}-2k\frac{\partial}{\partial t}{\xi _1}-k^2\xi _1 .
    \end{equation}
    Therefore $\xi _1$ is a classical solution of~\eqref{eq: ClosedLoop}.
\end{proof}

\end{document}